\documentclass[a4paper,12pt,reqno]{amsart}
\usepackage[top=25truemm,bottom=25truemm,left=20truemm,right=20truemm]{geometry}
\usepackage[dvipsnames]{xcolor}
\usepackage[dvipdfmx]{graphicx}
\usepackage{float}

\usepackage{subcaption}

\usepackage{latexsym}

\newcommand{\thm}{\sc{Theorem}\ }
\newcommand{\prop}{\sc{Proposition}\ }
\newcommand{\cor}{\sc{Corollary}\ }

\newcommand{\lem}{\sc{Lemma}\ }
\newcommand{\rem}{\sc{Remark}\ }

\usepackage[dvipdfmx]{hyperref}
\hypersetup{colorlinks=true,citecolor=Magenta,linkcolor=Magenta}

\newcommand{\R}{\mathbb{R}}
\newcommand{\bs}[1]{\boldsymbol{#1}}

\usepackage{color}
\begin{document}
\begin{center}
  {\bf \Large  Dual pairs of generic conformally flat hypersurfaces}\\[4mm]
\end{center}
\begin{center}
{ \large Yoshihiko Suyama\footnote{Fukuoka University: E-mail; suyama@fukuoka-u.ac.jp}}\\[4mm]
\end{center}

{\bf Abstract}. \ We study generic conformally flat (local-)hypersurfaces in the Euclidean 4-space $\R^4$. Such a hypersurface $f$ has the dual (hypersurface) $f^*$ in $\R^4$, which is also generic and conformally flat.
By repeating the composite action of inversion and the dual transformation on a hypersurface $f$, infinitely many non-equivalent generic conformally flat hypersurfaces are obtained from a single $f$.
The dual $f^*$ is defined by the total differential of the embedding expressed in terms of the original $f$. However, the exact formula in $\R^4$ of $f^*$ is not obvious, because of difficulty of integrating the total differential. 
Therefore, for the study of generic conformally flat hypersurfaces it is important to clarify in some way an explicit correspondence between the dual pair in $\R^4$.  The aim of this paper is to clarify that correspondence between the dual pair $f$ and $f^*$, which we do by making approximate discrete hypersurfaces of $f^*$ for all positive integers $n$ as maps from $3$-dimensional nets in $f$ to $\R^4$.
The approximations are constructed from the dual invariants of a generic conformally flat hypersurface $f$, of which dual invariants are defined as the maps from $f$ to $\R^4$, and as $n$ tends to $\infty$, the approximations induce maps between the corresponding curvature surfaces of $f$ and $f^*$.\\



{\bf 1. Introduction}

We study generic conformally flat (local-)hypersurfaces in the Euclidean 4-space $\R^4$. Here, we say that a hypersurface is generic if its principal curvatures are distinct from each other at each point, and a curvature surface of a hypersurface means a surface woven of the curvature lines for two principal curvatures of the hypersurface. Any generic conformally flat hypersurface $f$ has the dual (hypersurface) $f^*$ in $\R^4$, which is also generic and conformally flat (\cite{hsuy}). In this paper, for a given hypersurface $f$ and every positive integer $n$, we construct two kind of discrete hypersurfsces of $f^*$ as maps from a $3$-dimensional net in $f$ to $\R^4$, which induce a map from $f$ to $f^*$ as $n$ tends to $\infty$. 

Let a map $f: U(\subset \R^3)\rightarrow \R^4$ be a generic conformally flat hypersurface, where $U$ is a connected open set. Let $\kappa_i \ (i=1,2,3)$ be the principal curvatures of $f$. In this paper, for the sake of simplicity we assume that $\kappa_3$ is the middle principal curvature: $\kappa_1>\kappa_3>\kappa_2$ or $\kappa_1<\kappa_3<\kappa_2$. Then, there is a principal curvature line coordinate system $(x,y,z)$ taken in the order of $\kappa_i \ (i=1,2,3)$ and a function $\varphi(x,y,z)$ on $U$ such that the metric $g$,
$$g:=\cos^2\varphi(dx)^2+\sin^2\varphi (dy)^2+(dz)^2, \eqno{(1.1)}$$
is conformally flat on $U$ and the first fundamental form $I_f$ of $f$ is expressed as $$I_f=P^2g:=P^2(\cos^2\varphi(dx)^2+\sin^2\varphi (dy)^2+(dz)^2) \eqno{(1.2)}$$
with a function $P=P(x,y,z)> 0$ on $U$, where $\varphi(x,y,z)$ takes values in a connected open set of interval $(-\pi,\pi)$ satisfying $\sin\varphi\cos\varphi\neq 0$. Thus, we can recognize that the hypersurface $f(x,y,z)$ is a one parameter family of $(x,y)$-curvature surfaces with parameter $z$. 
We determine the unit normal vector field $N$ of $f$ so that $(\kappa_1-\kappa_2)\sin\varphi\cos\varphi(=P^{-1})>0$. Then, the principal curvatures $\kappa_i \ (i=1,2)$ are denoted as
$$\kappa_1=P^{-1}\tan\varphi+\kappa_3, \ \ \ \kappa_2=-P^{-1}\cot\varphi+\kappa_3 \eqno{(1.3)}$$
by $\varphi=\varphi(x,y,z)$, $P=P(x,y,z)$ and $\kappa_3=\kappa_3(x,y,z)$, where $P^{-1}(x,y,z):=(1/P)(x,y,z)$ (see Section 2 and \cite[Theorem 8]{ct}, \cite{hs3}). We call the conformally flat metric $g$ on $U$ in (1.2) the Guichard net\footnote{We referred the canonical principal Guichard net of $f$ to simply ``the Guichard net" of $f$ (see Section 2).} of $f.$  

Conversely, for a conformally flat metric $g$ (or a function $\varphi=\varphi(x,y,z)$) in (1.1) on $U$, there is a generic conformally flat hypersurface $f(x,y,z)$ with the Guichard net $g$ uniquely up to a conformal transformation of $\R^4$, if $U$ is simply connected (\cite{he1}, \cite{he2}, \cite[Section 2]{hs1}). Hence, the existence of a generic conformally flat hypersurface is {\it theoretically} equivalent to that of a conformally flat metric $g$ (or a function $\varphi=\varphi(x,y,z)$) on $U$ in (1.1). Here, the term ``theoretically" means that 
we have to find the functions $P=P(x,y,z)$ and $\kappa_3=\kappa_3(x,y,z)$ from $\varphi=\varphi(x,y,z)$ by the Gauss and the Codazzi equations to realize a hypersurface in $\R^4$ (see Section 4 and \cite{hs3}, \cite[Proposition 2.3]{su4}). The conformally flat metric $g$ on $U$ in (1.1) is also called a Guichard net.

Now, for a generic conformally flat hypersurface $f(x,y,z)$ on $U$, let us denote its dual in $\R^4$ by $f^*(x,y,z)$. Then, the system $(x,y,z)$ of $f^*$ is also the principal curvature line coordinates associated with the Guichard net $g^*$ of $f^*$, and further among the principal curvatures $\kappa_i^*$ of $f^*$ the third curvature $\kappa_3^*$ is middle (see Section 2). 
With $q\in \R^4$, let $\iota_q(x):=(x-q)/\parallel x-q\parallel^2+{ q}$ be the inversion in a unit sphere with center $q$. 
Then, a $4$-dimensional set of non-equivalent generic conformally flat hypersurfaces is induced from a single $f(x,y,z)$ by $(\iota_q\circ f)^*(x,y,z)$ with $q$. Moreover, an infinite dimensional set of non-equivalent generic conformally flat hypersurfaces is obtained from $f$ by repeating the above transformation to $f$ as $(\iota_{q_2}\circ(\iota_{q_1}\circ f)^*)^*$, $(\iota_{q_3}\circ(\iota_{q_2}\circ(\iota_{q_1}\circ f)^*)^*)^*$ and so on (see Section 2 and \cite{hsuy}). Therefore, it is important for the study of generic conformally flat hypersurfaces to clarify in some way an explicit correspondence in $\bs R^4$ between $f$ and $f^*$, since obtaining the exact formula of the dual is difficult as mentioned below.

Let $S$ be the Schouten $(1,1)$-tensor of $f$. The $S$ commutes with the shape operator of $f$ and symmetric with respect to the first fundamental form $I_f$. Hence, the principal curvature vectors $\partial/\partial x$, $\partial/\partial y$ and $\partial/\partial z$ are the eigen-vectors of $S$. Let $\sigma_i \ (i=1,2,3)$ be the eigen-values of $S$ in the order of $(x,y,z)$. 
The dual $f^*$ of $f$ is defined as $df^*:=df\circ S=\sigma_1f_xdx+\sigma_2f_ydy+\sigma_3f_zdz$ by the total differential, and hence $f^*$ is determined uniquely up to a parallel translation and non-degenerate on $U^*:=\{(x,y,z)\in U|\det S\neq 0\}$ (i.e., $f^*$ is generic and conformally flat only on $U^*$). The following facts for $f^*$ follow directly from the definition (see Section 2 and \cite{hsuy}): $(f^*)^*=f$ (i.e., $df=df^*\circ S^*$) holds; $f$ and $f^*$ have the common Schouten $(0,2)$-tensor; the corresponding principal curvature directions in $\R^4$ of $f$ and $f^*$ are parallel to each other. However, the exact formula of $f^*$ is not apparent from the definition, because of difficulty of integrating $df^*$. Our aim in this paper is to study an exact correspondence in $\R^4$ between the dual pair, by making approximate discrete hypersurfaces of $f^*$ as maps from certain $3$-dimensional nets in $f$ to $\R^4$.

Now, we explain the results in each section.
In Section 2, we review the several results mentioned above. 
In particular, at the point $p_0:=(a_0,b_0,c_0)\in U\setminus U^*$, we study the meaning of the fact that the corresponding principal curvature directions in $\R^4$ of $f$ and $f^*$ are parallel on $U^*$ (Proposition 2.1).
In Section 3, for $f(x,y,z)$ we define a map $Inv^{D}(f)_3: U\rightarrow \R^4$ by
$$Inv^{D}(f)_3(x,y,z):=(f^*+\sigma_3f+\kappa_3N)(x,y,z). \eqno{(1.4)}$$
The map $Inv^{D}(f)_3$ is an invariant of $f$ on $U$ with respect to the dual transformation. In fact, we shall verify the following facts: the dual transformation $*: f\leftrightarrow f^*$ induces a dual action $*$ on $Inv^D(f)_3$, and then $Inv^D(f)_3$ is invariant with respect to the action (Lemma 3.4). Furthermore, the two maps $Inv^D(f)_i:=f^*+\sigma_if+\kappa_iN \ (i=1,2)$ are dual invariants of $f$ on $U^*$.

Now, among generic conformally flat hypersurfaces, the hypersurfaces with conformal product-type and cyclic Guichard net were explicitly constructed and completely classified (\cite{he1}, \cite{hs1}, \cite{la}--\cite{su3}): the Guichard nets in (1.1) of these hypersurfaces are characterized by the second derivative condition $\varphi_{xz}=\varphi_{yz}\equiv 0$ of $\varphi$.  
In Section 4, using $Inv^{D}(f)_3$ we obtain an exact dual of every normal form in these hypersurfaces. 
Here, the normal forms refer to the hypersurfaces having the simplest structures among hypersurfaces with conformal product-type and cyclic Guichard net, of which all hypersurfaces are obtained by conformal transformations of the normal forms. The structure of the normal forms falls into three classes (see Section 4). Then, the dual $f^*$ of every normal form $f$ will be a normal form with the same structure as $f$. However, $f$ and $f^*$ are not conformally equivalent except when some $\kappa_i$ vanishes identically. Even the hypersurfaces with conformal product-type and cyclic Guichard net would be difficult to find these explicit duals in any form if they are not of the normal forms. Because those duals $f^*$ are new hypersurfaces not included in such types (see Section 2).

Now, let us return our argument to any generic conformally flat hypersurface $f(x,y,z)$ on $U$. In Section 5, using the dual invariants $Inv^{D}(f)_i$, two kinds of approximations of $f^*(x,y,z)$ are constructed for each positive integer $n$. Let $Q=[x_0,x_e]\times[y_0,y_e]\times[z_0,z_e]$ be a compact cube in $U^*$ and $x_e-x_0=y_e-y_0=z_e-z_0=a>0$. 
For each positive integer $n$, let $x_0<x_1<\cdots <x_n=x_e$, $y_0<y_1<\cdots <y_n=y_e$ and $z_0<z_1<\cdots <z_n=z_e$ be the divisions of the three sides of $Q$ of equal length $\delta_n:=a/n$. 
Let $L_n$ be the lattice in $Q$ created by the coordinate lines passing through all points $(x_i,y_j,z_k)$. Then, note that $f(L_n)$ is a $3$-dimensional net in $f$. 
We construct two kinds of $\R^4$-valued functions $(f^*)^{\bar\delta_n}(x,y,z)$ and $(f^*)^{\underline{\delta}_n}(x,y,z)$ on $L_n$ so that both sequences $\{(f^*)^{\bar\delta_n}\}_n$ and $\{(f^*)^{\underline{\delta}_n}\}_n$ converge uniformly to $f^*$ on $Q$ as $n$ tends to $\infty$. 
These approximations $(f^*)^{\bar\delta_n}(x,y,z)$ and $(f^*)^{\underline{\delta}_n}(x,y,z)$ are constructed from $f(x,y,z)$ on $L_n$ in the method outlined below, by focusing respectively on the principal curvature $x$- and $y$-lines of each $(x,y)$-curvature surface $f^*(x,y,z_k)$ with fixed $z=z_k$.

Now, for $n$ and $z=z_k$ fixed, let $S^{\bar\delta_n}(x,y,z_k)$ and $S^{{\underline\delta}_n}(x,y,z_k)$ be the $(x,y)$-surfaces $(f^*)^{\bar\delta_n}(x,y,z_k)$ and $(f^*)^{{\underline\delta}_n}(x,y,z_k)$, respectively. We only describe the idea of constructing the approximate $(x,y)$-surface $S^{\bar\delta_n}(x,y,z_k)$, under the condition $S^{\bar\delta_n}(x_0,y_0,z_k)={\bf 0}$. 
Then, we use the fact that the surface $f^*(x,y,z_k)$ is a one parameter family of the principal curvature $x$-lines with parameter $y$. 
Firstly, for $i=0,1,\cdots,n-1$ and $j=0,1,\cdots,n$, we define an $x$-curve $\bs{u}^{\delta_n}_{i}(x,y_j,z_k)$ from $[x_i,x_{i+1}]$ to $\R^4$ as a kind of parallel curve of $f(x,y_j,z_k)-f(x_i,y_j,z_k)$ as follows: 
$$\bs{u}^{\delta_n}_i(x,y_j,z_k):=-\sigma_3(x_i,y_j,z_k)\left[f(x,y_j,z_k)\right]^x_{x_i}-\kappa_3(x_i,y_j,z_k)\left[N(x,y_j,z_k)\right]^x_{x_i}, \eqno{(1.5)}$$where $\left[f(x,y_j,z_k)\right]^x_{x_i}=f(x,y_j,z_k)-f(x_i,y_j,z_k)$. (Usually, for a plane curve $l(x)$, the parallel curve ${\bar l}(x)$ of $l(x)$ means a curve ${\bar l}(x):=l(x)+cn(x)$, where $c$ and $n(x)$ are a constant and the unit normal vector of $l(x)$, respectively.)
We have $\bs{u}^{\delta_n}_i(x_i,y_j,z_k)=\bs{0}$. Note that, if $Inv^{D}(f)_3(x,y_j,z_k)=\bs{0}$ would be satisfied, then we have $\left[f^*(x,y_j,z_k)\right]_{x=x_i}^x=-[(\sigma_3f+\kappa_3N)(x,y_j,z_k)]_{x=x_i}^x$. Next, for $j=0,1,\cdots,n$, we define an $x$-curve $\bs{u}^{\delta_n}(x,y_j,z_k)$ from $[x_0,x_e]$ to $\R^4$ by 
$$\bs{u}^{\delta_n}(x,y_j,z_k):=\Sigma_{l=0}^{i-1}\bs{u}^{\delta_n}_l(x_{l+1},y_j,z_k)+\bs{u}^{\delta_n}_i(x,y_j,z_k) \ \ \ \ for \ x\in [x_i,x_{i+1}], \ \ 0\leq i\leq n-1, \eqno{(1.6)}$$ 
where $\Sigma_{l=0}^{i-1}\bs{u}^{\delta_n}_l(x_{l+1},y_j,z_k)=:{\bf 0}$ for $i=0$. The curve $\bs{u}^{\delta_n}(x,y_j,z_k)$ will approximate the principal $x$-line $\left[f^*(x,y_j,z_k)-f^*(x_0,y_j,z_k)\right]$ on $[x_0,x_e]$. In particular, the tangent vectors of both curves $\bs{u}^{\delta_n}(x,y_j,z_k)$ and $f(x,y_j,z_k)$ will be parallel on $[x_0,x_e]$ (even at the division points $x=x_i$), as is the case for curves $f^*(x,y_j,z_k)$ and $f(x,y_j,z_k)$. 
Now, the $(x,y)$-surface $S^{\bar\delta_n}(x,y,z_k)$ is created by translating these approximate $x$-curves $\bs{u}^{\delta_n}(x,y_j,z_k) \ (j=0,1,\cdots, n)$ by the appropriate constant vectors $\bs{c}^{\delta_n}_{j,k}$, that is, $S^{\bar\delta_n}(x,y_j,z_k):=\bs{c}^{\delta_n}_{j,k}+\bs{u}^{\delta_n}(x,y_j,z_k)$ for $j\in\{0,1,\cdots,n\}$. As mentioned above, the approximation $S^{\bar\delta_n}(x,y,z_k)$ is constructed by focusing on $x$-curves of the surface $f^*(x,y,z_k)$. 

In Section 6, for every $(i,j)$ an approximate $y$-curve of $f^*(x,y,z)$ connecting two points $S^{\bar\delta_n}(x_i,y_j,z_k)$ and $S^{\bar\delta_n}(x_i,y_{j+1},z_k)$ is constructed, that is, the approximation $S^{\bar\delta_n}(x,y,z_k)$ is a discrete surface defined on $L_n\cap (\R^2_{(x,y)}\times\{z_k\})$.
 
Now, although it is difficult to find the specific integral of $df^*$, the relation between the two $(x,y)$-curvature surfaces $f(x,y,z_k)$ and $f^*(x,y,z_k)$ has been clear from the formulae of $\bs{u}^{\delta_n}_i(x,y_j,z_k)$ and $\bs{u}^{\delta_n}(x,y_j,z_k)$ in the approximate discrete surfaces $S^{\bar\delta_n}(x,y,z_k)$. In particular, we emphasize that these formulae give a more interesting relation between two surfaces $f(x,y,z_k)$ and $f^*(x,y,z_k)$ even if the exact integral of $df^*$ is found. 

On the other hand, by focusing on $y$-curves of the $(x,y)$-surface $f^*(x,y,z_k)$, the approximate discrete surface $S^{\underline{\delta}_n}(x,y,z_k)$ is constructed from $Inv^{D}(f)_3(x,y,z)$ in the same way. Furthermore, using the invariants $Inv^D(f)_2$ and $Inv^D(f)_1$, respectively, we construct the $z$-curves of $(f^*)^{\bar\delta_n}(x_i,y_j,z)$ and $(f^*)^{\underline{\delta}_n}(x_i,y_j,z)$ to approximate the curvature surfaces $f^*(x,y,z)$ with $y=y_j$ and $x=x_i$ fixed. 
Then, both $(f^*)^{\bar\delta_n}(x,y,z)$ and $(f^*)^{\underline{\delta}_n}(x,y,z)$ on $L_n$ are approximate discrete hypersurfaces of $f^*(x,y,z)$.

In the above construction of the approximate discrete hypersurfaces, the curvature surfaces of the hypersurface $f^*$ are an important object. The geometry on curvature surfaces of generic conformally flat hypersurfaces in $\R^4$ were studied in a series of papers (\cite{bhs}, \cite{ms}, \cite{su4}).  \\

{\bf 2. Preliminaries}

In this section, for later use, we summarize the known facts about generic conformally flat hypersurfaces $f$ and their duals $f^*$. Then. the fact in (2.3) below says that the corresponding principal curvature vectors in $\R^4$ of a dual pair $f$ and $f^*$ are parallel at all regular points of $f^*$. This seems to suggest that any degenerated point of $f^*$ is a cusp of a curvature line of $f^*$ (cf. \cite[Theorem 3.8]{ms}). We give a sufficient condition under which it is true (Proposition 2.1).

Now, let $f$ be a generic conformally flat hypersurface in $\R^4$ defined on a connected open set $U$ in $\R^3$. Let $\kappa_i \ (i=1,2,3)$ be the principal curvatures of $f$ and $X_{\alpha}$, $X_{\beta}$, $X_{\gamma}$ be the unit principal curvature vectors in $\R^4$ taken in the order of $\kappa_i \ (i=1,2,3)$. We assume that $\kappa_3$ is middle among $\kappa_i$, as in Section 1. Let $\alpha$, $\beta$, $\gamma$ be the dual differential 1-forms for the unit vectors on $U$ corresponding to $X_{\alpha}$, $X_{\beta}$, $X_{\gamma}$. Then, $f$ is conformally flat if and only if the conformal fundamental forms 
$$\theta^1:=\sqrt{(\kappa_1-\kappa_3)(\kappa_1-\kappa_2)}~\alpha, \ \ 
\theta^2:=\sqrt{(\kappa_1-\kappa_2)(\kappa_3-\kappa_2)}~\beta, \ \ 
\theta^3:=\sqrt{(\kappa_3-\kappa_2)(\kappa_1-\kappa_3)}~\gamma$$
are closed (\cite[Section 2]{hs1}): $\theta^i$ are determined uniquely up to signs as well as $X_{\alpha}$, $X_{\beta}$ and $X_{\gamma}$, and they are invariant of $f$ with respect to a conformal transformation of $\R^4$. A principal curvature line coordinate system $(x,y,z)$ of $f$ is determined such that $\theta^1=dx$, $\theta^2=dy$ and $\theta^3=dz$: the coordinates $x$, $y$ and $z$ are determined up to signs and constant of integration.  Thus, the first fundamental form of $f(x,y,z)$ is expressed as in (1.2):
$$I_f=\alpha^2+\beta^2+\gamma^2= P^2g:=P^2(\cos^2\varphi dx^2+\sin^2\varphi dy^2+dz^2), $$
where
$$P:=\frac{1}{\sqrt{(\kappa_3-\kappa_2)(\kappa_1-\kappa_3)}}(>0), \ \ \cos^2\varphi=\frac{\kappa_3-\kappa_2}{\kappa_1-\kappa_2}, \ \ 
\sin^2\varphi=\frac{\kappa_1-\kappa_3}{\kappa_1-\kappa_2}.\eqno{(2.1)} $$
The above conformally flat metric $g=\cos^2\varphi dx^2+\sin^2\varphi dy^2+dz^2$ on $U$ is called the {\it canonical} Guichard net of $f$, where the term ``canonical" refers to the conditions $\theta_1=dx$, $\theta_2=dy$ and $\theta_3=dz$. If we drop the condition of ``{\it canonical}", then there is a one parameter family $\{f_c\} \ (c\in \R\setminus\{0\})$ of non-equivalent generic conformally flat hypersurfaces with the Guichard net $g$ so that $f_1=f$. 
In this paper, we restrict our argument to generic conformally flat hypersurfaces $f$ with the canonical Guichard net $g$. However, it does not inherently limit hypersurfaces, as mentioned below.

We explain the above facts. For a generic conformally flat hypersurface $f$, the Wang's M$\rm\ddot{o}$bius curvature $w:=(\kappa_2-\kappa_3)/(\kappa_1-\kappa_2)$ is another conformal invariant: in our case, $w=-\cos^2\varphi$ holds. Two hypersurfaces are conformally equivalent if and only if they have the same conformal fundamental forms $\theta^i$ and the M$\rm\ddot{o}$bius curvature $w$ (\cite{cpw}). For a generic conformally flat hypersurface $f(x,y,z)$ with the canonical Guichard net $g$, there is an associated family $\{f_c(x,y,z)\} \ (c\in \R\setminus\{0\})$ of generic conformally flat hypersurfaces: $f_1=f$; two hypersurfaces $f_c$ and $f_{-c}$ are conformally equivalent, but $f_c$ and $f_{c'}$ for $c>c'>0$ are not conformally equivalent. Each hypersurface $f_c(x,y,z)$ is characterized by the conditions $\theta^i_c=c\theta^i$ and $w_c=w$. Hence, all $f_c(x,y,z)$ have the same Guichard net $g$ as $f(x,y,z)$. Then, the Guichard net $g$ of $f_c(x,y,z) \ (c\neq \pm1)$ is not canonical. However, the $g$ of $f_c(x,y,z) \ (c\neq \pm1)$ can be replaced with the canonical Guichard net $g_c$, which is determined by a slight deformation $\varphi_c(x_c,y_c,z_c)$ of $\varphi(x,y,z)$. In fact, the canonical principal coordinates $(x_c,y_c,z_c)$ of $f_c$ are given as $x_c:=cx$, $y_c:=cy$ and $z_c:=cz$ by $\theta^i_c=c\theta^i$. Then, since we have $w_c=w$ under the correspondence $(x_c,y_c,z_c)\leftrightarrow(x,y,z)$, we define $\varphi_c$ by $\varphi_c(x_c,y_c,z_c):=\varphi(x_c/c,y_c/c,z_c/c)=\varphi(x,y,z)$. 
The function $\varphi_c(x_c,y_c,z_c)$ determines the canonical Guichard net $g_c:=\cos^2\varphi_c(dx_c)^2+\sin^2\varphi_c(dy_c)^2+(dz_c)^2$ of $f_c(x_c/c,y_c/c,z_c/c)$.

From now on, we assume that a hypersurface $f$ has the canonical Guichard net $g$, and call the canonical Guichard net of $f$ simply ``the Guichard net" of $f$. 
For a given Guichard net $g$ on $U$, there is a generic conformally flat hypersurface in $\R^4$ with the Guichard net $g$ uniquely up to a conformal transformation, if $U$ is simply connected (\cite[Section 2]{hs1}).   
Note that, even if we change $\varphi(x,y,z)$ in (1.1) into a function $\tilde\varphi(x,y,z)$ so that $\cos\tilde\varphi(x,y,z)=\pm\cos\varphi(x,y,z)$ and $\sin\tilde\varphi(x,y,z)=\pm\sin\varphi(x,y,z)$, the function $\tilde\varphi$ determines the same Guichard net as $\varphi$, by definition. Similarly, the coordinate change $(x,y)\rightarrow (\tilde x:=y,\tilde y:=x)$ determines the same Guichard net: in this case, we need to replace $\varphi$ by $\tilde\varphi=\varphi\pm\pi/2$.  
The other function $\hat\varphi$ which defines the same Guichard net as $\varphi$ is given by the coordinate change $(x,y,z)\rightarrow (\hat x:=\pm x+c_1,\hat y:=\pm y+c_2,\hat z:=\pm z+c_3)$ with a constant vector $(c_1,c_2,c_3)$ as mentioned above: hence $\hat\varphi(\hat x,\hat y,\hat z):=\varphi(\pm{\hat x}\mp c_1,\pm{\hat y}\mp c_2,\pm{\hat z}\mp c_3)~(=\varphi(x,y,z))$. Our argument proceeds by fixing a canonical coordinate system $(x,y,z)$ for a given hypersurface $f$. We can determine the unit normal vector field $N$ of $f$ such that the principal curvatures $\kappa_i \ ( i=1.2)$ are given as (1.3). In fact, in the case of $\sin\varphi\cos\varphi>0$, we determine $N$ by (2.1) so that $\kappa_1>\kappa_2$. If $\sin\varphi\cos\varphi<0$, then we replace that $N$ with $-N$.

Now, we define the unit principal curvature vectors in $\R^4$ of $f(x,y,z)$ by $X_{\alpha}:=(P\cos\varphi)^{-1}\partial f/\partial x$, $X_{\beta}:=(P\sin\varphi)^{-1}\partial f/\partial y$ and $X_{\gamma}:=P^{-1}\partial f/\partial z$, where $P^{-1}(x,y,z):=(1/P)(x,y,z)$. Let $S$ be the Schouten $(1,1)$-tensor of $f$. Then, the principal curvature vectors $\partial/\partial x$, $\partial/\partial y$ and $\partial/\partial z$ are the eigen-vectors of $S$, as mentioned in Section 1. We denote by $\sigma_i \ (i=1,2,3)$ the eigen-values of $S$ taken in the order of the vectors, and then $\sigma_i$ are given by 
$$\sigma_1=\textstyle{\frac{1}{2}}(\kappa_1\kappa_2-\kappa_2\kappa_3+\kappa_3\kappa_1), \ \ 
\sigma_2=\textstyle{\frac{1}{2}}(\kappa_1\kappa_2+\kappa_2\kappa_3-\kappa_3\kappa_1), \ \ 
\sigma_3=\textstyle{\frac{1}{2}}(-\kappa_1\kappa_2+\kappa_2\kappa_3+\kappa_3\kappa_1) \eqno{(2.2)}$$
(\cite{hsuy}). Particularly, we have $2\sigma_3=P^{-2}+\kappa_3^2\ (>0)$ by (1.3). 

The dual $f^*(x,y,z)$ in $\R^4$ of $f(x,y,z)$ is defined by the total differential $df^*:=df\circ S=\sigma_1f_xdx+\sigma_2f_ydy+\sigma_3f_zdz$, and then $f^*(x,y,z)$ is also a generic conformally flat hypersurface on $U^*:=\{(x,y,z)\in U|\det S(x,y,z)\neq 0\}$ (\cite{hsuy}): the set $U\setminus U^*$ is the degenerate points of $f^*(x,y,z)$.
In the following, the matters concerning $f^*(x,y,z)$ will be discussed on $U^*$ unless otherwise noted.
The unit principal curvature vectors $X^*_{\alpha},\ X^*_{\beta}$ and $X^*_{\gamma}$ in $\R^4$ of $f^*(x,y,z)$ are given as $X^*_{\alpha}=X_{\alpha},\ X^*_{\beta}=X_{\beta}$ and $X^*_{\gamma}=X_{\gamma}$ by 
$$X^*_{\alpha}=f^*_x/(\sigma_1P\cos\varphi)=X_{\alpha}, \ \ X^*_{\beta}=f^*_y/(\sigma_2P\sin\varphi)=X_{\beta}, \ \ X^*_{\gamma}=f^*_z/(\sigma_3P)=X_{\gamma}. \eqno{(2.3)}$$
The first fundamental form $I_{f^*}$ is determined by 
$$I_{f^*}=(\sigma_1\alpha)^2+(\sigma_2\beta)^2+(\sigma_3\gamma)^2=P^2\left((\sigma_1)^2\cos^2\varphi(dx)^2+(\sigma_2)^2\sin^2\varphi(dx)^2+(\sigma_3)^2(dz)^2\right) $$
by $df^*=df\circ S$. Then, $(\sigma_1/\sigma_3)^2\cos^2\varphi+(\sigma_2/\sigma_3)^2\sin^2\varphi=1$ holds by (2.1) and (2.2). In particular, we define the unit normal vector $N^*$ of $f^*$ by $N^*:=-N$ (not $N^*=N$), and then the principal curvatures $\kappa_i^*$ of $f^*$ are given by $\kappa_i^*=-\kappa_i/\sigma_i$: 
$$ N^*:=-N, \ \ \ \kappa^*_i=-\kappa_i/\sigma_i. \eqno{(2.4)}$$
Then, since $\kappa_i\sigma_j-\kappa_j\sigma_i=(\kappa_j-\kappa_i)\sigma_k$ holds for any permutation $\{i,j,k\}$ of $\{1,2,3\}$, the third principal curvature $\kappa^*_3$ is middle among $\kappa_i^* \ (i=1,2,3)$ by 
$$\left((\kappa_1/\sigma_1)-(\kappa_3/\sigma_3)\right) \left((\kappa_3/\sigma_3)-(\kappa_2/\sigma_2)\right)=(\kappa_1-\kappa_3)(\kappa_3-\kappa_2)/\sigma_3^2>0.$$

In consequence, the metric $g^*:=(\sigma_1/\sigma_3)^2\cos^2\varphi(dx)^2+(\sigma_2/\sigma_3)^2\sin^2\varphi(dx)^2+(dz)^2$ is the Guichard net of $f^*(x,y,z)$ and we can denote the first fundamental form $I_{f^*}$ of $f^*(x,y,z)$ by 
$$I_{f^*}=(P^*)^2g^*=(P^*)^2\left(\cos^2\varphi^*(dx)^2+\sin^2\varphi^*(dy)^2+dz^2\right), \eqno{(2.5)}$$ 
where
$$P^*=\sigma_3P(>0), \ \ \cos\varphi^*=(\sigma_1/\sigma_3)\cos\varphi, \ \ \sin\varphi^*=(\sigma_2/\sigma_3)\sin\varphi. \eqno{(2.6)}$$
Note that the Guichard net $g^*$ of $f^*(x,y,z)$ is canonical. In fact, for the conformal fundamental forms $(\theta^*)^i$ of $f^*$, we have $(\theta^*)^1={\it sign}(\sigma_1)\cdot\theta^1$, $(\theta^*)^2={\it sign}(\sigma_2)\cdot\theta^2$ and $(\theta^*)^3=\theta^3$ by replacing $(\kappa_i,\alpha,\beta,\gamma)$ with $(\kappa_i^*,\sigma_1\alpha,\sigma_2\beta,\sigma_3\gamma)$ in the definition of $\theta^i$. Note that the formulae (2.5) and (2.6) are differentiable on $U$ by $\sigma_3>0$ on $U$.

The eigen-values $\sigma^*_i$ of the Schouten (1,1)-tensor $S^*$ of $f^*$ are given as $\sigma^*_i=1/\sigma_i$ by (2.2) and (2.4).  
In particular, both hypersurfaces $f(x,y,z)$ and $f^*(x,y,z)$ have the common Schouten $(0,2)$-tensor and $df=df^*\circ S^*$ holds: we express the equation $df=df^*\circ S^*$ as $f=(f^*)^*$, ignoring parallel translations of $f$.   
Furthermore, using $P^{-1}\tan\varphi=\kappa_1-\kappa_3$ and $-P^{-1}\cot\varphi=\kappa_2-\kappa_3$ in (1.3), we have
$$(P^*)^{-1}\tan\varphi^*+\kappa_3^*=\kappa_1^*, \ \ \ -(P^*)^{-1}\cot\varphi^*+\kappa_3^*=\kappa_2^* \eqno{(2.7)}$$
by (2.2), (2.4) and (2.6). Note that the formulae in (2.7) for $f^*$ follow from $N^*=-N$ and they are the same as the formulae in (1.3) for $f$. We call the above correspondence $*: f(x,y,z) \leftrightarrow f^*(x,y,z)$, including the coordinates $(x,y,z)$, the dual transformation between $f$ and $f^*$. Then, note that we have determined as $N^*:=-N$.

Next, we study the meaning of the conditions $X_{\alpha}=X_{\alpha}^*$ and $X_{\beta}=X_{\beta}^*$ on $U^*$, at a degenerate point $p_0\in U\setminus U^*$ of the dual $f^*(x,y,z)$. Suppose that $f(x,y,z)$ satisfies $({\it det}~ S)(p_0)=0$ at a point $p_0:=(a_0,b_0,c_0)\in U$. Then, since $2\sigma_3=P^{-2}+\kappa_3^2>0$ on $U$, we have $(\sigma_1\sigma_2)(p_0)=0$. Furthermore, we have either $\sigma_1(p_0)=0$ or $\sigma_2(p_0)=0$. In fact, if $\sigma_1(p_0)=\sigma_2(p_0)=0$, then we have $(\kappa_1\kappa_2)(p_0)=\kappa_3(p_0)=0$ by (2.2).  
In the following proposition, let $p_0:=(a_0,b_0,c_0)\in U$. \\

{\prop 2.1}. \ {\it Let $f(x,y,z)$ be a generic conformally flat hypersurface with the first fundamental form $I_f=P^2g$ defined on a connected open set $U\subset \R^3$, where $g$ is a Guichard net in (1.1) defined from a function $\varphi(x,y,z)$. Let $f^*(x,y,z)$ be the dual of $f$. The functions $\kappa_i$ and $\sigma_i$ be the principal curvatures and the eigen-values of the Schouten tensor $S$ for $f$, respectively. Then, we have the following facts (1) and (2):

(1) If $\sigma_1(p_0)=0$ and $(\sigma_1)_x(p_0)\neq 0$ hold at $p_0\in U$, then the $x$-curve $f^*(x,b_0,c_0)$ in $\R^4$ has a cusp at $x=a_0$. Here, $(\sigma_1)_x(p_0)\neq 0$ is equivalent to $(P^{-1})_x(p_0)\neq (\kappa_3\varphi_x)(p_0)$. 

(2) If $\sigma_2(p_0)=0$ and $(\sigma_2)_y(p_0)\neq 0$ hold at $p_0\in U$, then the $y$-curve $f^*(a_0,y,c_0)$ in $\R^4$ has a cusp at $y=b_0$. Here, $(\sigma_1)_y(p_0)\neq 0$ is equivalent to $(P^{-1})_y(p_0)\neq(\kappa_3\varphi_y)(p_0)$.}\\

{\it Proof.} \ Since $U$ is a connected open set and the first fundamental form $I_f$ of $f$ is positive-definite on $U$, we can assume that the functions $P\cos\varphi$ and $P\sin\varphi$ are positive on $U$. Then, by $X_{\alpha}=(P\cos\varphi)^{-1}f_x$ and $X_{\beta}=(P\sin\varphi)^{-1}f_y$, the unit vectors $X_{\alpha}$ and $X_{\beta}$ on $U$ have the same orientations as $f_x$ and $f_y$, respectively. We only verify the fact (1), since (2) is verified in the same way. 

Now, suppose that $\sigma_1(p_0)=0$ and $(\sigma_1)_x(p_0)\neq 0$ hold at $p_0=(a_0,b_0,c_0)\in U$. When $x$ varies, the function $\sigma_1(x,b_0,c_0)$ of $x$ changes sign at $x=a_0$ and $X_{\alpha}(x,b_0,c_0)$ for $x\neq a_0$ is the tangent vector of the $x$-curve $f^*(x,b_0,c_0)$ by $f^*_x(x,b_0,c_0)=(\sigma_1f_x)(x,b_0,c_0)$. Then, we have $f^*_x(p_0)=\bf 0$ and $X_{\alpha}(p_0)$ is also tangent to the curve $f^*(x,y_0,z_0)$ at $x=a_0$ by the continuity of $X_{\alpha}(x,b_0,c_0)$ at $x=a_0$.
Furthermore, when $x$ varies, the curve $f^*(x,y_0,z_0)$ changes its orientation at $x=x_0$ against $X_{\alpha}(p_0)$. These facts imply that the point $f^*(x_0,y_0,z_0)$ is a cusp of the $x$-curve $f^*(x,y_0,z_0)$. 

Next, $(\sigma_1)_x(p_0)\neq 0$ and $[(P^{-1})_x-\kappa_3\varphi_x](p_0)\neq 0$ are equivalent. In fact, we have the following equation by (1.3), (2.2) and Lemma 2.2 below: 
$$2(\sigma_1)_x=[-P^{-2}+\kappa_3^2+2\kappa_3P^{-1}\tan\varphi]_x=-2(P\cos^2\varphi)^{-1}[(P^{-1})_x-\kappa_3\varphi_x].$$
\hspace{\fill}$\Box$\\

In general, the condition $(P^{-1})_x(p_0)\neq (\kappa_3\varphi_x)(p_0)$ or $(P^{-1})_y(p_0)\neq (\kappa_3\varphi_y)(p_0)$ at the point $p_0$ in Proposition 2.1 is satisfied, since we have $(P^{-1})_x=-(\kappa_3)_x\cot\varphi$ and $(P^{-1})_y=(\kappa_3)_y\tan\varphi$ on $U$, as mentioned in Lemma 2.2 below (these formulae have been used in the above proof).  
 
Now, the dual pair $f$ and $f^*$ are conformally equivalent if and only if $\kappa_i\equiv 0$ holds on $U$ for some $\kappa_i$. In fact, the M$\rm\ddot{o}$bius curvature $w^*=-\cos^2\varphi^*$ of $f^*$ is determined from $f$ by $\cos\varphi^*=(\sigma_1/\sigma_3)\cos\varphi$ in (2.6). Hence, $f$ and $f^*$ are conformally equivalent if and only if $\sigma_1=\pm \sigma_3$ holds on $U$. Furthermore, the equation $\sigma_1=\pm \sigma_3$ on $U$ implies that $\kappa_1\kappa_2\kappa_3\equiv 0$ holds on $U$, by (2.2). Suppose that $\kappa_1\kappa_2\kappa_3\equiv 0$ holds on $U$. Then, for each $\kappa_i$, the set $U'_i$ of points $p\in U$ such that $\kappa_i(p)=0$ is open and closed in $U$, because we have $(\kappa_j\kappa_k)(p)\neq 0$ for $p\in U'_i$, where $i\neq j\neq k\neq i$.

In the case $\kappa_1\equiv 0$, the vector $X_{\alpha}$ is independent of $x$ and all $x$-curves $f(x,y,z)$ with fixed $(y,z)$ are lines in $\R^4$ (cf. \cite{su3}, \cite[(2.2.3)]{su4}): the similar fact holds true in both cases of $\kappa_2\equiv 0$ and $\kappa_3\equiv 0$.

In the following Lemmata 2.2, 2.3, and 2.4, we state the other known facts about $f$ and $f^*$: \\

{\lem 2.2}. (cf. \cite[Proposition 2.2]{su4}) \ {\it We have the following equations on $U$:
$$(\kappa_3)_x=-(P^{-1})_x\tan\varphi, \ \ (\kappa_3)_y=(P^{-1})_y\cot\varphi, \ \ (\kappa_3)_z=-P^{-1}\varphi_z.$$ }\\

{\lem \sc{and} \sc{Definition} 2.3}. \ \ {\it We define a function $K=K(x,y,z)$ for $f(x,y,z)$ by  $$K:=-(P^{-1})_z+\kappa_3\varphi_z.$$
Then, the function $K^*=K^*(x,y,z)$ for $f^*(x,y,z)$ is given by $K^*=-K/\sigma_3$ on $U$.
Furthermore, we have the following equations on $U$: 
$$(\varphi^*)_x=\varphi_x-(\kappa_1/\sigma_3)(P^{-1})_x, \ \ 
(\varphi^*)_y=\varphi_y-(\kappa_2/\sigma_3)(P^{-1})_y, \ \ (\varphi^*)_z=-\varphi_z+(\kappa_3/\sigma_3)K.$$ }\\

{\it{Proof}.} \ We have $(\sigma_3)_z=-P^{-1}K$ by $\sigma_3=(P^{-2}+\kappa_3^2)/2$ and $(\kappa_3)_z=-P^{-1}\varphi_z$ in Lemma 2.2. Then, we have $[(P^*)^{-1}]_z=[P^{-1}/\sigma_3]_z=[P^{-2}K+(P^{-1})_z\sigma_3]/\sigma_3^2.$ Hence, if the formula for $(\varphi^*)_z$ holds true, then we obtain $K^*=-K/\sigma_3$ by the definition $K^*:=-[(P^*)^{-1}]_z+\kappa_3^*(\varphi^*)_z.$ 

The three formulae for $(\varphi^*)_x$, $(\varphi^*)_y$ and $(\varphi^*)_z$ have been verified in \cite[(5.13) and Corollary 5.2]{su4}. In that paper, the unit normal vector $N^*$ of $f^*$ is taken as $N^*=N$, but these formulae are independent of the choice of $N^*=\pm N$. 
\hspace{\fill}$\Box$\\

The function $K=K(x,y,z)$ for $f$ was introduced in \cite{su3}: for hypersurfaces with the Guichard net $g$ determined from $\varphi$ satisfying $\varphi_{xz}=\varphi_{yz}\equiv 0$, the function $K$ has a significant geometrical meaning (see Lemma 4.1 in Section 4).  \\

{\lem 2.4}. \ \ {\it We have the following equations on $U$: 
$$K_x=\kappa_2\varphi_{zx}, \ \ K_y=\kappa_1\varphi_{zy}, \ \ (\sigma_3)_x=\kappa_2(\kappa_3)_x, \ \ (\sigma_3)_y=\kappa_1(\kappa_3)_y, \ \ (\sigma_3)_z=-P^{-1}K.$$
In particular, if $\varphi_{xz}=\varphi_{yz}\equiv 0$, then $K$ depends only on $z$. If $\kappa_3$ depends only on $z$, then so does $\sigma_3$. If $K=(\kappa_3)_x=(\kappa_3)_y\equiv 0$ holds, then $\sigma_3$ is a positive constant.
}\\

{\it{Proof}.} \ The equations follow directly from Lemma 2.2, $2\sigma_3=P^{-2}+\kappa_3^2$ and the following equations in (\cite{hs3}, \cite[(2.3.1)]{su4}): 
$$\begin{array}{lll}
(P^{-1})_{xy}-(P^{-1})_y\varphi_x\cot\varphi+(P^{-1})_x\varphi_y\tan\varphi=0, \\[3mm]
(P^{-1})_{zx}+(P^{-1})_x\varphi_z\tan\varphi-P^{-1}\varphi_{zx}\cot\varphi=0,\\[3mm]
(P^{-1})_{zy}-(P^{-1})_y\varphi_z\cot\varphi+P^{-1}\varphi_{zy}\tan\varphi=0.
\end{array} \eqno{(2.8)}$$
\hspace{\fill}$\Box$\\

Next, let $\langle {\bf a},{\bf b}\rangle$ and $\parallel{\bf a}\parallel$ be the canonical inner product and the standard norm for vectors of $\R^4$, respectively. With $q\in \R^4$, let $\iota_q(x):=(x-q)/\parallel x-q\parallel^2+q$ be the inversion in a unit sphere with center $q$.  
We denote by $(\varphi^q,P^q,\kappa_i^q, \sigma_i^q, X_{\alpha}^q, \cdots)$ the geometric quantities for the hypersurface $\iota_q\circ f$ corresponding to $(\varphi,P,\kappa_i, \sigma_i, X_{\alpha},\cdots)$ for $f$. 
Then, since $f$ and $\iota_q\circ f$ are conformally equivalent, we have $\varphi^q=\varphi$ and obtain the following equations by direct calculation:
$$\begin{array}{ll} 
P^q=\displaystyle\frac{P}{\parallel f-q\parallel^{2}}, \ \ \kappa_i^q=\parallel f-q\parallel^2\kappa_i+2\langle N,f-q\rangle, \ \ K^q=\parallel f-q\parallel^2K+2\langle A,f-q\rangle,\\[3mm]
\sigma_i^q=\parallel f-q\parallel^4\sigma_i+2\kappa_i\parallel f-q\parallel^2\langle N,f-q\rangle+2\langle N,f-q\rangle^2,
\end{array} \eqno{(2.9)}$$
where $A:=-X_{\gamma}+\varphi_zN$, and for each vector $Z=X_{\alpha}, X_{\beta}, X_{\gamma}$ or $N$ of $f$ we have
$$Z^q=Z-2\langle Z,f-q\rangle~(f-q)/\parallel f-q\parallel^2. \eqno{(2.10)}$$

The vector $A=A(x,y,x)$ above and a vector $Kf+A$ for $f$ were defined in \cite{su3}: $Kf+A$ has an important geometric meaning for the hypersurfaces with the Guichard net determined from $\varphi$ satisfying $\varphi_{xz}=\varphi_{yz}\equiv 0$ (see Proposition 4.2 in Section 4). In the next section, more important application of $Kf+A$ to all generic conformally flat hypersurfaces will be proposed (see Theorem 3.2): $Kf+A$ is an inversion invariant for such a hypersurface, and our dual invariant of the hypersurface is induced from $Kf+A$.

Now, the M$\rm\ddot{o}$bius curvature of $(\iota_q\circ f)^*$ is given by $(w^q)^*=-(\sigma_1^q/\sigma_3^q)^2\cos^2 \varphi$, as mentioned before. Hence, by (2.9) we obtain a $4$-dimensional set of non-equivalent generic conformally flat hypersurfaces from a single $f$ by $(\iota_q\circ f)^*$ with $q\in \R^4$ (more precisely, see \cite[Corollary 3.6]{hsuy}). Moreover, we obtain an infinite dimensional set of non-equivalent generic conformally flat hypersurfaces from $f$ by $(\iota_{q_2}\circ (\iota_{q_1}\circ f)^*)^*$, $(\iota_{q_3}\circ(\iota_{q_2}\circ (\iota_{q_1}\circ f)^*)^*)^*$ and so on. 
Therefore, for the study of generic conformally flat hypersurfaces it is important to clarify in some way the correspondence between the images in $\R^4$ of each dual pair, although it is difficult to obtain the exact expression $f^*(x,y,z)$ from $df^*=df\circ S$. In the next section $3$, for a hypersurface $f$ we define three dual invariants $Inv^D(f)_i \ (i=1,2,3)$ and our argument will proceed on the basis of these invariants $Inv^D(f)_i$: in Section $4$ we obtain an exact dual for every normal form in the hypersurfaces with conformal product-type and cyclic Guichard net; in Sections $5$ and $6$, with respect to an arbitrarily given hypersurface, we construct two approximate discrete hypersurfaces of its dual for each positive integral $n$, that is, they induce the correspondence between the dual pair as $n$ tends to $\infty$. \\

{\bf 3. Two kinds of invariants determined for a hypersurface}

Let $f(x,y,z)$ be a generic conformally flat hypersurface in $\R^4$ with the first fundamental form $I_f=P^2g$ defined on a connected open set $U\subset \R^3$, where $g=\cos^2\varphi(dx)^2+\sin^2\varphi(dy)^2+(dz)^2$ is a Guichard net in (1.1). Let $\iota_q(x)=(x-q)/\parallel x-q\parallel^2+q$ be the inversion of $\R^4$, discussed in the previous section.  
In this section, we define two kinds of invariants for each hypersurface $f(x,y,z)$, which are determined as the maps from $U$ or $U^*$ into $\R^4$ and invariant with respect to the action of the inversions $\iota_q$ and the dual transformation to $f$, respectively.

Now, let $K=-(P^{-1})_z+\kappa_3\varphi_z$ be the function determined for $f(x,y,z)$, in Lemma 2.3, and $(X_{\alpha},X_{\beta},X_{\gamma},N)$ be the unit orthonormal frame field of $f$, given in the previous section. \\

{\lem 3.1}. \ \ {\it Let $A=-X_{\gamma}+\varphi_zN$ be the $\R^4$-valued function in (2.9), defined for $f(x,y,z)$ on $U$. 
Then, the following equations hold on $U$: 
$$[Kf+A]_x=\varphi_{zx}(\kappa_2f+N), \ \ [Kf+A]_y=\varphi_{zy}(\kappa_1f+N),$$
$$[Kf+A]_z=\varphi_{zz}(\kappa_3f+N)-(P^{-1}f)_{zz}-P^{-1}\varphi_z^2f.$$
}

{\it{Proof}.} \ Let $\nabla'$ be the canonical connection of $\R^4$: $\nabla'_{\partial/\partial x}=d/dx$. Now, we have $\kappa_1=P^{-1}\tan\varphi+\kappa_3$ and $\kappa_2=-P^{-1}\cot\varphi+\kappa_3$ in (1.3). Then, the first two equations in the lemma are derived directly from Lemma 2.4 and the following equations in \cite[(2.2.3)]{su4}:
$$\begin{array}{ll}
(X_{\gamma})_x=\nabla'_{\partial/\partial x}X_{\gamma}=(P^{-2}/\cos\varphi)(P\cos\varphi)_zf_x, \\[3mm]
 (X_{\gamma})_y=\nabla'_{\partial/\partial y}X_{\gamma}=(P^{-2}/\sin\varphi)(P\sin\varphi)_zf_y. 
\end{array} \eqno{(3.1)}$$
The third equation follows from $(\kappa_3)_z=-P^{-1}\varphi_z$ in Lemma 2.2 and $X_{\gamma}=P^{-1}f_z$.
\hspace{\fill}$\Box$\\

Since the formula for $[Kf+A]_z$ in Lemma 3.1 is complicated, we do not use it in our argument.  
Now, the $\R^4$-valued function  $A^q=-X_{\gamma}^q+\varphi^q_zN^q$ for $\iota_q\circ f$ corresponds to $A$ for $f$, where $\varphi^q=\varphi$.  
The following theorem follows directly from (2.9), (2.10) and $\varphi^q=\varphi$.\\

{\thm 3.2}. \ \ {\it With $q\in \R^4$, the hypersurfaces $f(x,y,z)$ and $(\iota_q\circ f)(x,y,z)$ satisfy the following relation:
$$K^q(\iota_q\circ f-q)+A^q=K(f-q)+A.$$}\\

We replace $\iota_q(x)$ with the inversion $\hat\iota_q(x):=(x-q)/\parallel x-q\parallel^2$. Note that $(\hat\iota_q\circ f)^*=(\iota_q\circ f)^*$ holds by the definition $df^*=df\circ S$. Then, the equation in Theorem 3.2 is rewritten as $K^q(\hat\iota_q\circ f)+A^q=(Kf+A)-Kq$. 
When we change $f$ into $\hat\iota_{q_1}\circ f$ in Theorem 3.2, we have 
$$K^{q_1,q_2}(\hat\iota_{q_2}\circ(\hat\iota_{q_1}\circ f))+A^{q_1,q_2}=
K^{q_1}(\hat\iota_{q_1}\circ f)+A^{q_1}-K^{q_1}q_2=(Kf+A)-Kq_1-K^{q_1}q_2,$$
where $K^{q_1,q_2}:=(K^{q_1})^{q_2}$ and $A^{q_1,q_2}:=(A^{q_1})^{q_2}$. 
Hence, the map $Kf+A$ determined for $f$ is invariant with respect to the repeated action of inversion $\hat\iota_q$ on $f$, in the above sense.  
We call $Kf+A$ the inversion invariant of $f$, and denote $Inv^{I}(f)(x,y,z):=(Kf+A)(x,y,z)$. The inversion invariant $Inv^{I}(f)(x,y,z)$ of $f$ is a map from $U$ into $\R^4$.

Next, after reviewing the action of the dual operator $*$ on geometric quantities of $f$ and $f^*$, we define the dual invariant of $f$ from the inversion invariant $Inv^I(f)$.  
Let $\kappa^*$ be a function or a vector field determined for $f^*$ corresponding to $\kappa$ for $f$, as in the previous section.

For a hypersurface $f(x,y,z)$, we have denoted its dual by $f^*(x,y,z)$. 
Hence, the coordinates $(x,y,z)$ of the dual $f^*$ are regarded as $(x^*,y^*,z^*)=(x,y,z)$. Then, since $\varphi^*(x,y,z)=\varphi^*(x^*,y^*,z^*)=(\varphi(x,y,z))^*$, we have $(\varphi^*)_z=(\varphi_z)^*$. In fact, for $h\neq 0$ we have 
$$(\varphi^*(x,y,z+h)-\varphi^*(x,y,z))/h=(\varphi^*(x^*,y^*,z^*+h^*)-\varphi^*(x^*,y^*,z^*))/h^*$$
$$=[(\varphi(x,y,z+h)-\varphi(x,y,z))/h]^*.$$ 
In the same way, we have $((P^{-1})_z)^*=((P^{-1})^*)_z=((P^*)^{-1})_z$, because $(P^{-1})^*=(P^*)^{-1}$ holds by $(P^{-1})^*=(1/P)^*=(1/P^*)=(P^*)^{-1}$. In consequence, we can denote $\varphi^*_z:=(\varphi^*)_z=(\varphi_z)^*$ and $(P^{-1})^*_z:=((P^{-1})_z)^*=((P^{-1})^*)_z=((P^*)^{-1})_z$.

Next, for $\kappa_i^*=-\kappa_i/\sigma_i$, we obtain $(\kappa_i^*)^*=-\kappa_i^*/\sigma_i^*$. In fact, we have $(\kappa_i^*)^*=\kappa_i$, $\kappa_i^*=-\kappa_i/\sigma_i$ and $\sigma_i^*=1/\sigma_i$. 
In the same way, for $K=-(P^{-1})_z+\kappa_3\varphi_z$ and $\varphi^*_z=-\varphi_z+(\kappa_3/\sigma_3)K$ in Lemma 2.3, we obtain $K^*=-[(P^{-1})_z]^*+\kappa_3^*(\varphi_z)^*$ and $(\varphi^*_z)^*=-(\varphi_z)^*+(\kappa_3^*/\sigma_3^*)K^*$.
In fact, we have $K^*=-K/\sigma_3$, $(\varphi^*_z)^*=\varphi_z$ and   
$$-((P^{-1})_z)^*+\kappa_3^*(\varphi_z)^*=-((P^*)^{-1})_z+\kappa_3^*\varphi^*_z=((\sigma_3)_zP^{-1}-\kappa_3^2K)/\sigma_3^{2}+K/\sigma_3=-K/\sigma_3,$$ 
$$-(\varphi_z)^*+(\kappa^*_3/\sigma^*_3)K^*=-\varphi^*_z+(\kappa_3/\sigma_3)K=\varphi_z$$ 
by $2\sigma_3=P^{-2}+\kappa_3^2$ and Lemma 2.4.

Furthermore, the vector field $A^*=-X_{\gamma}^*+\varphi_z^*N^*$ and the inversion invariant $K^*f^*+A^*$ of $f^*$ corresponds on $A=-X_{\gamma}+\varphi_zN$ and $Kf+A$ of $f$, respectively. Hence, we can recognize the action of the dual operator $*$ on $A$ and $Kf+A$ as follows: 
$$A^*=(-X_{\gamma}+\varphi_zN)^*:=-X_{\gamma}^*+\varphi^*_zN^*, \ \ \  (Kf+A)^*:=K^*f^*+A^*. \eqno{(3.2)}$$ 
Then, the following desired equations also hold:
$$(A^*)^*=-(X_{\gamma}^*)^*+(\varphi^*_z)^*(N^*)^*=-X_{\gamma}+\varphi_zN=A,$$
$$(K^*f^*+A^*)^*=(K^*)^*(f^*)^*+(A^*)^*=Kf+A. $$

Now, the dual operator $*$ also acts on the $\R^4$-valued functions $(Kf+A)\pm(K^*f^*+A^*)$ on $U$ determined from $f$ and $f^*$ as follows:
$$\begin{array}{ll}
\left\{(Kf+A)\pm(K^*f^*+A^*)\right\}^*:=(Kf+A)^*\pm(K^*f^*+A^*)^* \\[3mm]
=\pm\left\{(Kf+A)\pm(K^*f^*+A^*)\right\}.
\end{array}  \eqno{(3.3)}$$
The equations (3.3) mean that  
two $\R^4$-valued functions $(Kf+A)\pm(K^*f^*+A^*)$ are invariant under the dual action $*$. In particular, we have
$$(Kf+A)-(K^*f^*+A^*)=(K/\sigma_3)(f^*+\sigma_3f+\kappa_3N). \eqno{(3.4)}$$
by direct calculation. When $K\equiv 0$, we have $(Kf+A)-(K^*f^*+A^*)={\bf 0}$. However, even in this case, the $\R^4$-valued function $f^*+\sigma_3f+\kappa_3N$ on $U$ is well defined and has an important implication for our argument, as shown in Lemma 3.4 and Theorem 3.5 below. \\

{\rem 3.3}. \ (1) The equation $\varphi^q(x,y,z)=\varphi(x,y,z)$ holds, but we have $(\varphi^q)^*\neq \varphi^*$. In fact, $\varphi^*$ is determined by $\cos\varphi^*=(\sigma_1/\sigma_3)\cos\varphi$ and $\sin\varphi^*=(\sigma_2/\sigma_3)\cos\varphi$, on the other hand $(\varphi^q)^*$ is determined by $\cos(\varphi^q)^*=(\sigma_1^q/\sigma_3^q)\cos\varphi$ and $\sin(\varphi^q)^*=(\sigma_2^q/\sigma_3^q)\cos\varphi$, from (2.6) and (2.9). In particular, we have $(\varphi^q)^*_z=-\varphi_z+(\kappa_3^q/\sigma^q_3)K^q$ by Lemma 2.3. 
As shown in the previous section, we have $(\varphi^{q_1})^*\neq (\varphi^{q_2})^*$ in general if $q_1\neq q_2$.

(2) As the result of (1), the dual operator $*$ is meaningful only for the action between a fixed dual pair $f$ and $f^*$, their functions or vector fields. 
Hence, although the equation $K^q(\hat\iota_q\circ f)+A^q=K(f-q)+A$ holds in Theorem 3.2, we have $[K^q(\hat\iota_q\circ f)+A^q]^*\neq [K(f-q)+A]^*$. Here, we have $[K^q(\hat\iota_q\circ f)+A^q]^*=(K^q)^*(\iota_q\circ f)^*+(A^q)^*$ and $[K(f-q)+A]^*=K^*f^*+A^*$, by $(\hat\iota_q\circ f)^*=(\iota_q\circ f)^*$ and $(f-q)^*=f^*$, and then $(K^q)^*(\iota_q\circ f)^*+(A^q)^*$ and $K^*f^*+A^*$ are the inversion invariants of $(\iota_q\circ f)^*$ and $f^*$, respectively.\\

Now, we have the following lemma 3.4 and theorem 3.5: \\

{\lem 3.4}. \ \ {\it The following equations hold on $U$: 
$$(f^*+\sigma_3f+\kappa_3N)^*=(1/\sigma_3)(f^*+\sigma_3f+\kappa_3N),$$
$$\left\{(1/\sqrt{\sigma_3}(f^*+\sigma_3f+\kappa_3N)\right\}^*=(1/\sqrt{\sigma_3})(f^*+\sigma_3f+\kappa_3N).$$
}\\

{\it{Proof}.} \ Firstly, note that $f^*+\sigma_3f+\kappa_3N$ is a linear combination of $\R^4$-valued functions determined from $f$ and $f^*$. Now, we have 
$$(f^*+\sigma_3f+\kappa_3N)^*=(f^*)^*+\sigma_3^*f^*+\kappa_3^*N^*=f+(1/\sigma_3)f^*+(\kappa_3/\sigma_3)N,$$ 
which shows the first equation. The second equation follows from $2\sigma_3=P^{-2}+\kappa_3^2>0$, $(\sigma_3)^*=(1/\sigma_3)$ and the first equation. 
\hspace{\fill}$\Box$\\

The $\R^4$-valued function $f^*+\sigma_3f+\kappa_3N$ on $U$ is a dual invariant of $f$, in the sense of Lemma 3.4. We denote $Inv^D(f)_3:=f^*+\sigma_3f+\kappa_3N$ on $U$. Furthermore, the $\R^4$-valued functions $f^*+\sigma_if+\kappa_iN \ (i=1,2)$ are also dual invariants of $f$. However, the definition domain of these invariants is restricted to $U^*$. In fact, only for the points $(x,y,z)$ such that $\sigma_i(x,y,z)\neq 0$, we have 
$$[f^*+\sigma_if+\kappa_iN]^*=f+(1/\sigma_i)f^*+(\kappa_i/\sigma_i)N=(1/\sigma_i)(f^*+\sigma_if+\kappa_iN).$$
We denote $Inv^D(f)_i:=f^*+\sigma_if+\kappa_iN \ (i=1,2)$ on $U^*$.

The following Theorem 3.5 follows from (2.2), Lemmata 2.2--2.4 and the definition of $f^*$. \\

{\thm 3.5}. \ \ {\it We have the following equations on $U$ for $Inv^D(f)_3=f^*+\sigma_3f+\kappa_3N$: 
$$[f^*+\sigma_3f+\kappa_3N]_x=(\kappa_3)_x(\kappa_2f+N), \ \ \ 
[f^*+\sigma_3f+\kappa_3N]_y=(\kappa_3)_y(\kappa_1f+N). \leqno{(1)}$$
$$ \left[ f^* + \sigma_3 f+\kappa_3 N \right]_z=-P^{-1}(Kf+A).$$
$$ \left[ (1/\sqrt{\sigma_3}) (f^* + \sigma_3 f+\kappa_3 N) \right]_z
=-(2P\sqrt{\sigma_3})^{-1}\left\{(Kf+A)+(K^*f^*+A^*)\right\}, \leqno{(2)}$$
$$\left[(1/\sigma_3)(f^*+\sigma_3f+\kappa_3N)\right]_z=-(P^*)^{-1}(K^*f^*+A^*).$$}\\

The last equation of Theorem 3.5-(1) gives the relation between the two kinds of invariants $Inv^D(f)_3=f^* + \sigma_3 f+\kappa_3 N$ and $Inv^I(f)=Kf+A$ of $f$, and the first equation of (2) implies that the $\R^4$-valued function $(Kf+A)+(K^*f^*+A^*)$ in (3.3) is induced from the dual invariant $Inv^D(f)_3$.
 
When we notice the relation between the equations in Lemma 3.1 and Theorem 3.5-(1), we have the following corollary.  
In the corollary and the proof, let $\varphi$, $P$, $\kappa_i$, $\sigma_i$, $X_{\alpha}$, $X_{\beta}$, $X_{\gamma}$, $N$ and $A$ be the geometric quantities of $f$ (or $f-q$). The notation $\bs{a}\parallel\bs{b}$ means that the two vectors $\bs{a}$ and $\bs{b}$ in $\R^4$ are parallel to each other.  \\

{\cor 3.6}. \ \ {\it For any $q\in \R^4$, we have the following facts:
$$\left[Inv^D(f-q)_3\right]_x \ \parallel \ \left[Inv^D(\hat\iota_q\circ f)_3 \right]_x \ \parallel \ \left[Inv^I(f-q)\right]_x,$$
$$\left[Inv^D(f-q)_3\right]_y \ \parallel \ \left[Inv^D(\hat\iota_q\circ f)_3 \right]_y \ \parallel \ \left[Inv^I(f-q)\right]_y,$$
$$\left[ Inv^D(f-q)_3 \right]_z \ \parallel \ \left[Inv^D(\hat\iota_q\circ f)_3 \right]_z \ \parallel \ Inv^I(f-q). \ \ \ $$}\\

{\it{Proof}.} \ Firstly, we note that $(f-q)^*=f^*$ and $(\hat\iota_q\circ f)^*=(\iota_q\circ f)^*$ hold. Now, we have the following equation by (2.9), (2.10) and Theorem 3.5-(1): 
$$\left[(f-q)^*+\sigma_3(f-q)+\kappa_3 N\right]_x=(\kappa_3)_x(\kappa_2(f-q)+N),$$ 
$$\left[(f-q)^*+\sigma_3(f-q)+\kappa_3N\right]_y=(\kappa_3)_y(\kappa_1(f-q)+N),$$
$$\left[ (f-q)^* + \sigma_3 (f-q)+\kappa_3 N) \right]_z=-P^{-1}(K(f-q)+A),$$
$$\left[(\hat\iota_q\circ f)^* + \sigma_3^q (\hat\iota_q\circ f)+\kappa_3^q N^q \right]_x=(\kappa_3^q)_x(\kappa_2(f-q)+N),$$  
$$\left[(\hat\iota_q\circ f)^* + \sigma_3^q (\hat\iota_q\circ f)+\kappa_3^q N^q \right]_y=(\kappa_3^q)_y(\kappa_1(f-q)+N),$$
$$\left[(\hat\iota_q\circ f)^* + \sigma_3^q (\hat\iota_q\circ f)+\kappa_3^q N^q \right]_z=-\parallel f-q\parallel^2 P^{-1}(K(f-q)+A).$$
\hspace{\fill}$\Box$\\

{\bf 4. Normal forms in certain classes of hypersurfaces, and their duals}

Let the Guichard net $g=\cos^2\varphi(dx)^2+\sin^2\varphi(dy)^2+(dz)^2$ be defined by a function $\varphi(x,y,z)$ on a connected domain $U\subset \R^3$. In this section, we assume that $\varphi(x,y,z)$ satisfies $\varphi_{zx}=\varphi_{zy}\equiv 0$. Then, $\varphi$ is expressed as $\varphi(x,y,z)=\lambda(x,y)+\mu(z)$, and $\lambda$ and $\mu$ satisfy the equations 
$$\lambda_{xx}-\lambda_{yy}=(c_1^2/2)\sin(2\lambda), \ \ \ \mu\equiv 0 \ \ or \ \ \mu_z^2=c_2^2-c_1^2\sin^2\mu,  \eqno{(4.1)}$$
respectively, where $c_1$ and $c_2$ are non-zero constants (\cite{he1}, \cite[Section 4.1]{hs1}): in the case $\lambda_{xx}-\lambda_{yy}=-(c_1^2/2)\sin(2\lambda)$, we change the coordinates $(x,y)$ into $(y,x)$; if $\mu_z\equiv 0$, then $\mu(z)\equiv 0$ by the definition.  
If a hypersurface $f(x,y,z)$ has the Guichard net $g$ determined by $\varphi(x,y,z)=\lambda(x,y)$ ($\mu\equiv 0$), $f$ is called conformal product-type (\cite{he1}). If a hypersurface $f(x,y,z)$ has the Guichard net $g$ determined by $\varphi(x,y,z)=\lambda(x,y)+\mu(z)$ ($\mu_z\neq 0$), $f$ is called having a cyclic Guichard net (\cite{hs1}). All the hypersurfaces with conformal product-type and cyclic Guichard net were characterized and completely classified (\cite{he1}, \cite{hs1}, \cite{spt}). Then, regardless of whether $\mu\equiv 0$ or $\mu_z\neq 0$, such a hypersurface is conformally equivalent to a hypersurface with one of the following three structures (a), (b) and (c):

(a) Every $(x,y)$-curvature surfaces of $f(x,y,z)$ with fixed $z$ lies in a hyperplane $H(z)$, and $H(z)$ translates in parallel about $z$. In other words, the unit normal vector ${\bs B}(z)$ of $H(z)$ is constant: ${\bs B}(z)=:{\bs B}$. 

(b) Every $(x,y)$-curvature surfaces of $f(x,y,z)$ with fixed $z$ lies in a linear hyperplane $H(z)$, and the unit normal vector ${\bs B}(z)$ of $H(z)$ moves on a unit circle in a linear plane $L$.    

(c) All $(x,y)$-curvature surfaces of $f(x,y,z)$ with fixed $z$ lie on concentric $3$-spheres $H(z)$ centered at the origin.

\noindent
We shall say that a hypersurface $f$ having the structure (a), (b) or (c) is the normal form of its conformal equivalence class. By the way, the hypersurface with the structure (a), (b) or (c) is created as an evolution in $z$ from a suitable linear Weingarten surface in a Euclidean $3$-space $H(z_0)$ with fixed $z_0$, in a hyperbolic $3$-space $H(z_0)^+$ of some constant curvature $-1/r^{2}$ or in a standard $3$-sphere $H(z_0)$ of some constant curvature $1/r^2$, where $H(z_0)^+$ is an upper half hyperplane $L^{\perp}\times \{t {\bs B}(z_0)^{\perp}\}_{t>0}$ with respect to a plane $L^{\perp}$ perpendicular to $L$ and a unit vector ${\bs B}(z_0)^{\perp}$ in the plane $L$ perpendicular to ${\bs B}(z_0)$. However, we shall not go into such properties of surfaces that create the hypersurfaces.

Now, let a Guichard net $g$ on $U$ be defined from $\varphi(x,y,z)=\lambda(x,y)+\mu(z)$ in (4.1), and let $f(x,y,z)$ be a normal form with the first fundamental form $I_f=P^2g$. Let $(X_{\alpha}, X_{\beta}, X_{\gamma},N)$ be the unit orthonormal frame field of $f$, given in Section 2.
As shown later (Lemma 4.1 and Proposition 4.2), we have the following facts for the normal forms $f(x,y,z)$, using the inversion invariant $Inv^I(f)=Kf+A$ of $f$.
In the case $\mu\equiv 0$: ${\bs B}(z)=X_{\gamma}(x,y,z)$ in the structure (a) and (b); for $f(x,y,z)$ with the structure (a) or (c), every $z$-curve $f(x,y,z)$ with fixed $(x,y)$ is a line of $\R^4$; for $f(x,y,z)$ with the structure (b), every $z$-curve $f(x,y,z)$ with fixed $(x,y)$ is a circle of $\R^4$. In the case $\mu_z\neq 0$: ${\bs B}(z)=-(A/\sqrt{1+\mu_z^2})(x,y,z)$ in the structure (a) and (b); for all $f(x,y,z)$, every $z$-curve $f(x,y,z)$ with fixed $(x,y)$ is intrinsically circular but not extrinsically circular.
 
In this section, we obtain an exact dual $f^*$ of the normal form $f$ with conformal product-type and cyclic Guichard net, which is given as a formula related with the above (a)-(c) structures of $f$, respectively.  
That is, let $f$ be a normal form. Then, we have the following facts for the dual $f^*$ of $f$ (Theorem 4.3): $\varphi^*_{xz}=\varphi^*_{yz}\equiv 0$ hold, and hence $\varphi^*$ is also expressed as $\varphi^*(x,y,z)=\lambda^*(x,y)+\mu^*(z)$; $\mu^*\equiv 0$ holds if $\mu\equiv 0$, and $\mu^*_z\neq 0$ holds if $\mu_z\neq0$; the dual $f^*$ is a normal form with the same structure (a), (b) or (c) as $f$. However, $f$ and $f^*$ are not conformally equivalent if any principal curvature $\kappa_i$ of $f$ does not vanish identically, as mentioned in Section 2.

In Proposition 4.2 later, we also study the {\it properties} of the function $P$ and the third principal curvature $\kappa_3$. Here, the term ``the property of $P$" implies a sufficient condition on $P$ under which the aforementioned results on $f^*$ are obtained (not necessarily the exact formula of $P$). That is, after acknowledging the existence of hypersurfaces $f$ in $\R^4$ with the structure (a), (b) or (c), we study the properties of $P$ and $\kappa_3$: these proofs are quite different from those in the papers mentioned at the beginning of this section.

In order to study the properties of $P$ and $\kappa_3$, we firstly recall the equations for $P$, which are the condition that, for a Guichard net $g=\cos^2\varphi(dx)^2+\sin^2\varphi(dy)^2+(dz)^2$ (not necessarily $\varphi(x,y,z)=\lambda(x,y)+\mu(z)$), there is a generic conformally flat hypersurface $f$ in $\R^4$ with the first fundamental form $I_f=P^2g$ (cf. \cite{hs3}, \cite[Proposition 2.3]{su4}).  Now, let $P^{-1}$ be a solution to the following differential equation determined by $\varphi$,
$$\begin{array}{ll}
P^{-1}=-\left[(P^{-1})_{xx}+2\varphi_x\tan\varphi(P^{-1})_x\right]-\left[(P^{-1})_{yy}-2\varphi_y\cot\varphi(P^{-1})_y\right] \\[2mm]
\hspace{2cm}+(P^{-1})_{zz}+2\psi_{zz}P^{-1}, 
\end{array} \eqno{(4.2)}$$
where $\psi_{zz}:=(\varphi_{xx}-\varphi_{yy}-\varphi_{zz}\cos2\varphi)/\sin2\varphi$. For a solution $P^{-1}$ to (4.2), we define the functions $\kappa_3=\kappa_3(x,y,z)$ and $\zeta=\zeta(x,y,z)$ by 
$$\begin{array}{ll}
\kappa_3:=\left(\tan\varphi (P^{-1})_{xx}-\varphi_x\frac{\cos2\varphi}{\cos^2\varphi}(P^{-1})_x\right)-\left(\cot\varphi(P^{-1})_{yy}-\varphi_y\frac{\cos2\varphi}{\sin^2\varphi}(P^{-1})_y\right)  \\[3mm]
\hspace{2cm} +\left(P^{-1}\varphi_{zz}-(P^{-1})_z\varphi_z\right),
\end{array} \eqno{(4.3)}$$
$$\zeta:=P^{-1}\left[2(P^{-1})_{zz}+(-1+\varphi_z^2+2\psi_{zz})P^{-1}\right]
-\left(\frac{[(P^{-1})_x]^2}{\cos^2\varphi}+\frac{[(P^{-1})_y]^2}{\sin^2\varphi}+[(P^{-1})_z]^2\right). \eqno{(4.4)}$$
Then, the existence condition for $f$ in $\R^4$ with the metric $I_f=P^2g$ is that a non-zero solution $P^{-1}$ to (4.2) satisfies the following four equations: the three equations in (2.8) and $\zeta=\kappa_3^2$.

Now, the derivative $\nabla'_{X_{\gamma}}X_{\gamma}$ is expressed as
$$\nabla'_{X_{\gamma}}X_{\gamma}=aX_{\alpha}+bX_{\beta}+\kappa_3N, \ \ \ {\rm where} \ \  a=(P^{-1})_x/\cos\varphi, \ \ b=(P^{-1})_y/\sin\varphi \eqno{(4.5)}$$ (cf. \cite{su3}, \cite[(2.2.3) and (2.2.4)]{su4}).

In the following Lemma, the vector ${\bs B}(z)$ is the unit normal vector of the hyperplane $H(z)$ in the structures (a) and (b). Note that we have $A=-X_{\gamma}$ if $\mu\equiv 0$, since $A=-X_{\gamma}+\mu_zN$. \\

{\lem 4.1}. \ \ {\it Let a Guichard net $g$ on a connected open set $U\subset \R^3$ be defined by a function $\varphi(x,y,z)=\lambda(x,y)+\mu(z)$ in (4.1). Suppose that $f(x,y,z)$ is a hypersurface in $\R^4$ with the first fundamental form $I_f=P^2g$.
Then, the functions $a(x,y,z)$ and $b(x,y,z)$ in (4.5) are independent of $z$: $a(x,y,z)=:a(x,y)$ and $b(x,y,z)=:b(x,y)$, and $K(x,y,z)$ depends only on $z$: $K(x,y,z)=:K(z)$. In particular, if $K(z_0)=0$ at $z=z_0$, then the $(x,y)$-surface $f(x,y,z_0)$ lies in a hyperplane. If $K(z_0)\neq 0$ at $z=z_0$, then the $(x,y)$-surface $f(x,y,z_0)$ lies in a 3-sphere. 
Furthermore, when $f(x,y,z)$ is a normal form, we have the following facts: 

(1) In the case that $f(x,y,z)$ has the structure (a): Regardless of whether $\mu\equiv 0$ and $\mu_z\neq 0$, we have $K(z)\equiv 0$. Furthermore, the following facts hold. 

If $\mu\equiv 0$, then $X_{\gamma}(x,y,z)$ is constant and we have ${\bs B}=X_{\gamma}(x,y,z)$.

If $\mu_z\neq 0$, then $A(x,y,z)/\sqrt{1+\mu_z^2}$ is constant and we have ${\bs B}=-A(x,y,z)/\sqrt{1+\mu_z^2}$. 

(2) In the case that $f(x,y,z)$ has the structure (b): Regardless of whether $\mu\equiv 0$ and $\mu_z\neq 0$, we have $K(z)\equiv 0$. Furthermore, the following facts hold.

If $\mu\equiv 0$, then $X_{\gamma}(x,y,z)$ depends only on $z$ and we have ${\bs B}(z)=X_{\gamma}(x,y,z)$.

If $\mu_z\neq 0$, then $A(x,y,z)$ depends only on $z$ and we have ${\bs B}(z)=-A(x,y,z)/\sqrt{1+\mu_z^2}$. 

(3) In the case that $f(x,y,z)$ has the structure (c): Regardless of whether $\mu\equiv 0$ and $\mu_z\neq 0$, we have $K(z)\neq 0$ for all $(x,y,z)\in U$. Furthermore, the following facts hold.

If $\mu\equiv 0$, then we have $f(x,y,z)=X_{\gamma}(x,y,z)/K(z)$ (in Proposition 4.2 below, we show that $X_{\gamma}(x,y,z)$ is independent of $z$).

If $\mu_z\neq 0$, then we have $f(x,y,z)=-A(x,y,z)/K(z)$.  }\\

{\it{Proof}.} \ Suppose that a hypersurface $f(x,y,z)$ in $\R^4$ has the first fundamental form $I_f=P^2g$ and its Guichard net $g$ is defined by $\varphi(x,y,z)=\lambda(x,y)+\mu(z)$ on $U$: we have $\varphi_{xz}=\varphi_{yz}\equiv 0$. Then, in the equation $\nabla'_{X_{\gamma}}X_{\gamma}=aX_{\alpha}+bX_{\beta}+\kappa_3N$ of (4.5), the functions $a$ and $b$ do not depend on $z$. In fact, we have 
$$a_z=[(P^{-1})_x/\cos\varphi]_z=(1/\cos\varphi)\left((P^{-1})_{xz}+(P^{-1})_x\mu_z\tan\varphi\right)\equiv 0,$$
$$b_z=[(P^{-1})_y/\sin\varphi]_z=(1/\sin\varphi)\left((P^{-1})_{yz}-(P^{-1})_x\mu_z\cot\varphi\right)\equiv 0$$
by (2.8) and $\varphi_{zx}=\varphi_{zy}\equiv 0$.

Next, for $K=-(P^{-1})_z+\kappa_3\mu_z$ in Lemma 2.3 and $A=-X_{\gamma}+\mu_zN$, we have $K_x=K_y=0$ and $[Kf+A]_x=[Kf+A]_y={\bf 0}$ by Lemma 2.4 and Lemma 3.1.  Hence, $K(x,y,z)$ depends only on $z$: $K(x,y,z)=:K(z)$, and there is a vector $\bs{a}(z)$ depending only on $z$ such that $(Kf+A)(x,y,z)=\bs{a}(z)$. If $K(z_0)\neq 0$ at $z=z_0$, then we have $f(x,y,z_0)=-A(x,y,z_0)/K(z_0)+\bs{a}(z_0)/K(z_0)$. Hence, the $(x,y)$-surface of $f(x,y,z_0)$ lies in a $3$-sphere of radius $(\sqrt{1+\mu_z^2}/|K|)(z_0)$ and center $\bs{a}(z_0)/K(z_0)$ by the definition of $A$. If $K(z_0)=0$ at $z=z_0$, then $A(x,y,z_0)$ does not depend on $x$ and $y$ by $A_x(x,y,z_0)=A_y(x,y,z_0)={\bf 0}$ and ${\bs a}(z_0)=A(x,y,z_0)$ is perpendicular to $X_{\alpha}(x,y,z_0)$ and $X_{\beta}(x,y,z_0)$ by the definition of $A(x,y,z)$. Hence, 
the $(x,y)$-surface $f(x,y,z_0)$ lies in a hyperplane perpendicular to ${\bs a}(z_0)=A(x,y,z_0)$. Note that either $K(z)=0$ or $K(z)\neq 0$ holds for each point $(x,y,z)\in U$.

In consequence, by the definitions of the structures (a), (b), (c) and $A=-X_{\gamma}$ when $\mu\equiv 0$, we have verified the facts (1) and (2). 
In fact, for the hypersurface $f(x,y,z)$ with the structure (a) or (b), we have $K\equiv 0$ and ${\bs B}(z)=-(A/\sqrt{1+\mu_z^2})(z)$ by $\langle f_z,X_{\gamma}\rangle>0$. For the hypersurface with the structure (c), we have $K(z)\neq 0$ for all $z$ and $f(x,y,z)=-A(x,y,z)/K(z)$, since the above vector $\bs{a}(z)/K(z)$ is a zero vector by the condition of the structure (c).
\hspace{\fill}$\Box$\\

Let $\nabla$ be the connection of the normal form $f$ induced from $\nabla'$ of $\R^4$. Then, by (4.5) and $a_z=b_z\equiv 0$ in Lemma 4.1, we have $[(\nabla_{X_{\gamma}})^2X_{\gamma}](x,y,z)=-(a^2+b^2)(x,y)X_{\gamma}(x,y,z)$, which shows that all $z$-curves $f(x,y,z)$ with fixed $(x,y)$ are intrinsically circular. Then, if $(\kappa_3)_z\neq 0$, then these $z$-curves are not circles in $\R^4$. When $\mu_z\neq 0$, we have $(\kappa_3)_z\neq 0$ by $(\kappa_3)_z=-P^{-1}\mu_z$ in Lemma 2.2. 

By Lemma 4.1, we have the following facts for $inv^I(f)$ of the normal form $f$: if $f$ has the structure (a), 
then $Inv^I(f)(x,y,z)=-\sqrt{1+\mu_z^2}~{\bf B}$ with a unit constant vector ${\bf B}$; if $f$ has the structure (b), then $Inv^I(f)(x,y,z)=-\sqrt{1+\mu_z^2}~{\bf B}(z)$, where the vector ${\bf B}(z)$ of $z$ moves on a unit circle of a plane; if $f$ has the structure (c), then $Inv^I(f)(x,y,z)\equiv {\bf 0}$.

Now, in the following Proposition 4.2, we study the properties of $P$ and $\kappa_3$ for each normal form $f$ and identify $\varphi$ for $f$ that determines its Guichard net $g$. To do so, we have to study how the movement of ${\bs B}(z)$ in the structure (a) or (b) of $f$ relates to these $P$ and $\kappa_3$. Then, an exact form of $f$ associated with the structures (a)-(c) is also obtained. 
Let $c_1$ and $c_2$ be the non-zero constants in (4.1).  \\

{\prop 4.2}. \ \ {\it Let a Guichard net $g$ be defined by a function $\varphi(x,y,z)=\lambda(x,y)+\mu(z)$ in (4.1) on a connected open set $U\subset \R^3$. Suppose that $f(x,y,z)$ is a normal form in $\R^4$ with the first fundamental form $I_f=P^2g$. 
 Then, $f(x,y,z)$ satisfies the following facts according to its structure (a), (b) or (c).

(1) \underline{In the case that $f(x,y,z)$ has the structure (a)}: 

If $\mu\equiv 0$, then  
$\kappa_3(x,y,z)\equiv 0$, $c_1^2=1$ and $P$ is constant: $P(x,y,z)=:P(>0)$. In particular, when we take $X_{\gamma}=:(0,0,0,1)$, 
there is a surface ${\tilde f}(x,y)$ in a hyperplane $H(z_0)$ perpendicular to $X_{\gamma}$ such that 
$f$ is expressed as $f(x,y,z)=({\tilde f}(x,y),Pz)\subset H(z_0)\times \{tX_{\gamma}\}_{t\in \R}$.

If $\mu_z\neq 0$, then $P$ and $\kappa_3$ depend only on $z$:  $P(x,y,z)=:P(z)$ and $\kappa_3(x,y,z)=:\kappa_3(z)$, and they are expressed as $P^{-1}(z)=C\cos\mu(z)(>0)$ and $\kappa_3(z)=-C\sin\mu(z)$ with a constant $C\neq 0$. Furthermore, we have $c_1^2-c_2^2-1=0$. In particular, when we take ${\bs B}=:(0,0,0,1)$, there is a hyperplane $H(z_0)$ perpendicular to ${\bs B}$ and $z$-dependent $(x,y)$-surfaces $\tilde f(x,y,z)$ in $H(z_0)$ such that $f$ is expressed as 
$f(x,y,z)=(\tilde f(x,y,z), f^4(z))\subset H(z_0)\times \{t{\bs B}\}_{t\in \R}$, where $f^4(z)$ satisfies $(f^4)'(z)=P/\sqrt{1+\mu_z^2}=(Cc_1\cos^2\mu)^{-1}$ (here, we assumed $c_1\cos\mu>0$).

(2) \underline{In the case that $f(x,y,z)$ has the structure (b)}. In the description below, $L$ and $L^{\perp}$ are planes of $\R^4$ perpendicular to each other:

If $\mu\equiv 0$, then $\kappa_3$ and $P$ are independent of $z$: $\kappa_3(x,y,z)=:\kappa_3(x,y)$ and $P(x,y,z)=:P(x,y)$, and $\kappa_3(x,y)\neq 0$. Furthermore, we have $P^2(a^2+b^2+\kappa_3^2)=c_1^2-1>0$ and $a^2+b^2\neq 0$. In particular, we can take $X_{\gamma}(z)$ as 
$X_{\gamma}(z)=:(-\sin(\lower2pt\hbox{\small{$\sqrt{\smash[b]{c_1^2-1}}$}}(z-z_0)),\cos(\lower2pt\hbox{\small{$\sqrt{\smash[b]{c_1^2-1}}$}}(z-z_0)),0,0)$ for a fixed $z_0$.
When we define a function $f^1(x,y)$ by $f^1(x,y):=P(x,y)/\lower2pt\hbox{\small{$\sqrt{\smash[b]{c_1^2-1}}$}}~(>0)$, there is a surface $\tilde f(x,y):=(f^1(x,y),0,f^3(x,y),f^4(x,y))$ in the upper half hyperplane $H(z_0)^{+}\subset L\times L^{\perp}=\R^4$ such that $f(x,y,z)$ is expressed as 
$$f(x,y,z)=f^1(x,y)(\cos (\lower2pt\hbox{\small{$\sqrt{\smash[b]{c_1^2-1}}$}}(z-z_0)), \sin (\lower2pt\hbox{\small{$\sqrt{\smash[b]{c_1^2-1}}$}}(z-z_0)),0,0)+(0,0, f^3(x,y),f^4(x,y)).$$

If $\mu_z\neq 0$, then $P$ and $\kappa_3$ are expressed as 
$$P^{-1}(x,y,z)=C(x,y)\cos\mu(z)+D(x,y)\sin\mu(z), $$ 
$$\kappa_3(x,y,z)=-C(x,y)\sin\mu(z)+D(x,y)\cos\mu(z).$$ 
Here, the functions $C(x,y)$ and $D(x,y)$ of $(x,y)$ satisfy $C^2(c_1^2-c_2^2-1)-D^2(c_2^2+1)=a^2+b^2$ and $c_1^2-c_2^2-1>0$. The function $(a^2+b^2)(x,y)$ is positive on $U$. In particular, we can take ${\bs B}(z)$ as ${\bs B}(z)=:(-\sin\rho(z),\cos\rho(z),0,0)$, where the function $\rho(z)$ satisfies $\rho'(z)=\sqrt{(1+c_2^2)(c_1^2-c_2^2-1)}/(1+\mu_z^2(z))$ and $\rho(z_0)=0$ at a fixed $z_0$. When we define a function $f^1(x,y,z)$ by
$$f^1(x,y,z):=P(x,y,z)\sqrt{\textstyle\frac{1+\mu_z^2(z)}{(1+c_2^2)(c_1^2-c_2^2-1)}} \ (>0),$$ 
there are $z$-dependent $(x,y)$-surfaces $\tilde f(x,y,z)=(f^1(x,y,z),0,f^3(x,y,z),f^4(x,y,z))$
in the upper half hyperplane $H(z_0)^+\subset L\times L^{\perp}$ 
such that 
$f(x,y,z)$ is expressed as 
$$f(x,y,z)=f^1(x,y,z)(\cos\rho(z),\sin\rho(z),0,0)+(0,0,f^3(x,y,z),f^4(x,y,z)).$$
Regardless of $\mu\equiv 0$ and $\mu_z\neq 0$, the case $\kappa_1(P^{-1})_x=\kappa_2(P^{-1})_y\equiv 0$ on $U$ does not occur.

(3) \underline{In the case that $f(x,y,z)$ has the structure (c)}: 

If $\mu\equiv 0$, then $\kappa_3(x,y,z)\equiv 0$ and $X_{\gamma}(x,y,z)$ is independent of $z$: $X_{\gamma}(x,y,z)=:X_{\gamma}(x,y)$. Furthermore, $P$ depends only on $z$: $P(x,y,z)=:P(z)$, and $1-c_1^2>0$. Let $c_3:=\sqrt{1-c_1^2}$. Then, the function $P$ and the hypersurface $f$ are expressed as $P(z)=C\exp(c_3z)$ with a constant $C(> 0)$ and $f(x,y,z)=X_{\gamma}(x,y)/K(z)=(1/c_3)P(z)X_{\gamma}(x,y)$, respectively.

If $\mu_z\neq 0$, then $P$ and $\kappa_3$ depend only on $z$: $P(x,y,z)=:P(z)$ and $\kappa_3(x,y,z)=:\kappa_3(z)$, and $(\kappa_3)_z\neq 0$. Furthermore, we have $1+c_2^2-c_1^2>0$. Let $c_4:=\sqrt{(1+c_2^2)(1+c_2^2-c_1^2)}$. Then, the functions $P$, $\kappa_3$ and the hypersurface $f$ are expressed as, with a constant $C(>0)$ 
$$\begin{array}{ll}
P(z)=(C/\sqrt{1+\mu_z^2}~)\exp\left(c_4\int_0^z(1+\mu_z^2)^{-1}~dz\right), \ \ \kappa_3(z)=
\left((c_4\mu_z+\mu_{zz})/(1+\mu_z^2)\right)P^{-1}, \\[2mm] 
f(x,y,z)=-A(x,y,z)/K(z)=-(1/c_4)P(z)A(x,y,z).
\end{array}$$ }\\

{\it{Proof}.} \ Let a Guichard net $g$ be defined by a function $\varphi(x,y,z)=\lambda(x,y)+\mu(z)$ on $U\subset \R^3$. Let $f(x,y,z)$ be a normal form in $\R^4$ with the first fundamental form $I_f=P^2g$. Then, the function $K(x,y,z)$ satisfies either $K\equiv 0$ or $K(x,y,z)=K(z)\neq 0$ at all points $(x,y,z)\in U$ and the functions $a$ and $b$ in (4.5) are independent of $z$: $a(x,y,z)=:a(x,y)$ and $b(x,y,z)=:b(x,y)$, by Lemma 4.1. We have $(\kappa_3)_z=-P^{-1}\mu_z$ by Lemma 2.2. Hence, if $\mu\equiv 0$, then $\kappa_3$ is independent of $z$; if $\mu_z\neq 0$, then $(\kappa_3)_z\neq 0$. 

We verify the facts (1) and (2) in two separate cases (i) $\mu\equiv 0$ and (ii) $\mu_z\neq 0$.

(i) \underline{In the case that $f$ satisfies $K=\mu\equiv 0$}. Firstly, note that $X_{\gamma}$ depends only on $z$ by $A=-X_{\gamma}$ and Lemma 4.1. The condition $K\equiv 0$ implies that every $(x,y)$-surface $f(x,y,z)$ with fixed $z$ lies in the hyperplane $H(z)$ perpendicular to ${\bs B}(z)=X_{\gamma}(x,y,z)$, by Lemma 4.1.

Next, since $(P^{-1})_z=-K\equiv 0$ and $(\kappa_3)_z\equiv 0$ by $\mu\equiv 0$, the functions $P^{-1}$ and $\kappa_3$ are independent of $z$: $P(x,y,z)=:P(x,y)$ and $\kappa_3(x,y,z)=:\kappa_3(x,y)$. Since $\partial^2 X_{\gamma}/\partial z^2=-P^2(a^2+b^2+\kappa_3^2)X_{\gamma}$ and $P^2(a^2+b^2+\kappa_3^2)$ is a function of $(x,y)$ by $a=a(x,y)$ and $b=b(x,y)$, the function $P^2(a^2+b^2+\kappa_3^2)$ is constant. Moreover, we have $P^2(a^2+b^2+\kappa_3^2)=c_1^2-1$ by $\zeta=\kappa_3^2$ and $\zeta=(c_1^2-1)P^{-2}-(a^2+b^2)$.
 
Now, this case (i) is divided into the two cases $\kappa_3\equiv 0$ and $\kappa_3\neq 0$. 

(i-1) \underline{Suppose $\kappa_3\equiv 0$}. Then, every $z$-curve $f(x,y,z)$ with fixed $(x,y)$ is a line in $\R^4$ and ${\bs B}(z)=X_{\gamma}(x,y,z)$ is constant as mentioned in Section 2: ${\bs B}(z)={\bs B}$.
Hence, we have $a^2+b^2+\kappa_3^2\equiv 0$ and $c_1^2=1$. Furthermore, $P$ is a constant function by $a^2+b^2\equiv 0$ and the definitions of $a$ and $b$. Hence, by $f_z=PX_{\gamma}$, $f(x,y,z)$ has the form given when $\mu\equiv 0$ in (1).   

(i-2) \underline{Suppose $\kappa_3(x,y)\neq0$}. Then, $P^2(a^2+b^2+\kappa_3^2)=c_1^2-1$ is a non-zero constant on $U$, i.e., $c_1^2>1$. 
Hence, ${\bs B}(z)=X_{\gamma}(x,y,z)$ moves on a unit circle of a linear plane $L$ and every $z$-curve $f(x,y,z)$ with fixed $(x,y)$ is a circle in an affine plane parallel to $L$. Furthermore, since every $(x,y)$-surface $f(x,y,z)$ with fixed $z$ lies in a hyperplane $H(z)$ with normal vector ${\bs B}(z)$, the hypersurface $f(x,y,z)$ has the form given when $\mu\equiv 0$ in (2) by $f_z=PX_{\gamma}$. 
In particular, we have $a^2+b^2\neq 0$. In fact, if $P(x,y)$ would be constant, then $P^{-1}=c_1^2P^{-1}$ holds by (4.1) and (4.2), which contracts $c_1^2>1$. 

Next, we prove that the case $\kappa_1(P^{-1})_x=\kappa_2(P^{-1})_y\equiv 0$ on $U$ does not occur. If $\kappa_1(P^{-1})_x=\kappa_2(P^{-1})_y\equiv 0$ holds on $U$, then we have at least one of the following two cases (A) and (B), since $P(x,y)$ is not constant: 

(A) There is an open set $U'(\subset U)$ such that $(P^{-1})_x\neq 0$ at each point of $U'$. 

(B) There is an open set $\tilde{U}(\subset U)$ such that $(P^{-1})_y\neq 0$ at each point of $\tilde{U}$. 

\noindent
In the case (A), we have $\kappa_1=(P^{-1})_y\equiv 0$ on $U'$. In the case (B), we have $\kappa_2=(P^{-1})_x\equiv 0$ on $\tilde{U}$. Hence, it is sufficient to prove that the equations $\kappa_1=(P^{-1})_y\equiv 0$ and $(P^{-1})_x\neq 0$ on $U$ does not occur. Suppose that $\kappa_1=(P^{-1})_y\equiv 0$ and $(P^{-1})_x\neq 0$ hold on $U$. Then, we have $\lambda_x=\lambda_y\equiv 0$. In fact, $\kappa_3=-P^{-1}\tan\lambda$ holds by $\kappa_1\equiv 0$, and then we have $(\kappa_3)_x=-(P^{-1})_x\tan\lambda-P^{-1}(\lambda_x/\cos^2\lambda)$ and $(\kappa_3)_y=-(P^{-1})_y\tan\lambda-P^{-1}(\lambda_y/\cos^2\lambda)=-P^{-1}(\lambda_y/\cos^2\lambda)$. On the other hand, we have $(\kappa_3)_x=-(P^{-1})_x\tan\lambda$ and 
$(\kappa_3)_y=(P^{-1})_y\cot\lambda\equiv 0$ by Lemma 2.2, which show $\lambda_x=\lambda_y\equiv 0$.  
This conclusion that $\lambda$ is constant contradicts (4.1), by $c_1^2>0$. Therefore, the case $\kappa_1=(P^{-1})_y\equiv 0$ and $(P^{-1})_x\neq 0$ on $U$ does not occur. Note that the above proof of this result does not depend on $\mu\equiv 0$ and $\mu_z\neq 0$.

In consequence, the case (i) is divided into the two subcases (i-1) and (i-2). Then, the subcase (i-1) (resp. (i-2)) is the case that $f$ has the structure (a) (resp. the structure (b)), by the property of ${\bs B}(z)$ in Lemma 4.1.

(ii) \underline{In the case that $f$ satisfies $K \equiv 0$ and $\mu_z\neq 0$}. The vector $A/\sqrt{1+\mu_z^2}$ depends only on $z$, and the condition $K\equiv 0$ implies that every $(x,y)$-surface $f(x,y,z)$ with fixed $z$ lies in the hyperplane $H(z)$ perpendicular to ${\bs B}(z)=-(A/\sqrt{1+\mu_z^2})(x,y,z)$, by Lemma 4.1. 

Now, by $(P^{-1})_z=\kappa_3\mu_z$ and $(\kappa_3)_z=-P^{-1}\mu_z$, there are two functions $C=C(x,y)$ and $D=D(x,y)$ of $(x,y)$ such that $P^{-1}$ and $\kappa_3$ are expressed as  
$$P^{-1}=C\cos\mu+D\sin\mu, \ \ \ \kappa_3=-C\sin\mu+D\cos\mu.$$
Since the functions $\mu$ and $\psi_{zz}$ in (4.2) satisfy  
$$\mu_{zz}=-c_1^2\sin\mu\cos\mu, \ \ \mu_{zzz}=-c_1^2\mu_z\cos(2\mu), \ \ 2\psi_{zz}=c_1^2\cos(2\mu)$$
by (4.1), we obtain the fact that the equation $\zeta=\kappa_3^2$ is equivalent to $a^2+b^2=C^2(c_1^2-c_2^2-1)-D^2(c_2^2+1)$, by direct calculation.  

Next, for the vector ${\bs B}(z)$ of $z$: $\bs{B}=-{A}/{\sqrt{1+\mu_z^2}}={(X_{\gamma}-\mu_zN)}/{\sqrt{1+\mu_z^2}}$, 
we define the unit vectors $\bs{C}=\bs{C}(x,y,z)$ and $\bs{D}=\bs{D}(x,y,z)$ by $$\bs{C}:=(\mu_zX_{\gamma}+N)/{\sqrt{1+\mu_z^2}}, \ \ \ 
\bs{D}:=(aX_{\alpha}+bX_{\beta})/{\sqrt{a^2+b^2}}.$$
These vectors ${\bs B}$, $\bs{C}$ and $\bs{D}$ are perpendicular to each other and we have $X_{\gamma}=(1+\mu_z^2)^{-1/2}(\bs{B}+\mu_z\bs{C})$. Furthermore, we define the functions $l=l(x,y,z)$ and $m=m(x,y,z)$ by
$$l:=P\left(-\kappa_3+P^{-1}{\mu_{zz}}/(1+\mu_z^2)\right), \ \ \ 
m:=P\sqrt{(a^2+b^2)/(1+\mu_z^2)}.  \eqno{(4.6)}$$
Then, $m$ satisfies $m_z=-m\mu_zP\left[\kappa_3+P^{-1}{\mu_{zz}}/(1+\mu_z^2)\right]$. Furthermore, we have 
$$\bs{B}_z=-l\bs{C}+m\bs{D}, \ \ \bs{C}_z=l\bs{B}+\mu_zm\bs{D}, \ \ \bs{D}_z=-m\bs{B}-\mu_zm\bs{C}$$ 
and, at the points $(x,y,z)\in U$ such that $(l^2+m^2)(x,y,z)\neq 0$, we have 
$$-\left[\bs{B}_z/\sqrt{l^2+m^2}\right]_z=\sqrt{l^2+m^2}\bs{B}$$
$$+\left[(l_zm-lm_z)+\mu_zm(l^2+m^2)\right](l^2+m^2)^{-3/2}(m\bs{C}+l\bs{D}). \eqno{(4.7)}$$
Now, the case (ii) is divided into the subcases (ii-1) $m\equiv 0$ and (ii-2) $m\neq 0$.

(ii-1) \underline{Suppose $m(x,y,z)\equiv 0$}. Then, we have $a(x,y)=b(x,y)\equiv 0$, and hence the functions $C(x,y)$ and $D(x,y)$ are constant by $a\cos\varphi=P_x$ and $b\sin\varphi=P_y$. Furthermore, we have $C(c_1^2-c_2^2-1)\sin\mu(z)+D(c_2^2+1)\cos\mu(z)\equiv 0$ from the formula of $\kappa_3$ in (4.3). Hence, we have $D=0$ and $c_1^2-c_2^2-1=0$ by $P(z)\neq 0$. Note that $1+\mu_z^2=(c_1\cos\mu)^2$ and $C^2(c_1^2-c_2^2-1)-D^2(c_2^2+1)=a^2+b^2$ are satisfied. 

Now, we have obtained $P^{-1}=C\cos\mu$ and $\kappa_3=-C\sin\mu$ with a constant $C\neq 0$. 
Furthermore, we have $l\equiv 0$ by $c_1^2-c_2^2-1=0$, which implies that $\bs{B}(z)$ is a constant vector: $\bs{B}(z)=\bs{B}$. 

(ii-2) \underline{Suppose $m\neq 0$}. We take an open subset $U'$ of $U$ such that $(l^2+m^2)(x,y,z)\neq 0$ holds at every point $(x,y,z)\in U'$. For a while, our argument will go on $U'$ and we shall verify that the equation $\zeta=\kappa_3^2$: $a^2+b^2=C^2(c_1^2-c_2^2-1)-D^2(c_2^2+1)$, regulates the movement of the vector ${\bs B}(z)$ on $U'$. Now, for the coefficient of $m(l^2+m^2)^{-3/2}(m\bs{C}+l\bs{D})$ in the right hand side of (4.7), we have
$$(l_zm-lm_z)m^{-1}+\mu_z(l^2+m^2)=\left[{\mu_{zz}}/(1+\mu_z^2)\right]_z+\mu_z\left[1+\left(\mu_{zz}/(1+\mu_z^2)\right)^2\right]+\mu_z(l^2+m^2)$$
$$=\mu_z\left[-(1+c_2^2)(c_1^2-c_2^2-1)(1+\mu_z^2)^{-2}+(l^2+m^2)\right]. \eqno{(4.8)}$$
The equation (4.8) vanishes identically by $\zeta=\kappa_3^2$. In fact, we firstly have
$$(1+\mu_z^2)^{2}(l^2+m^2)=\left[(P^{-1}\mu_{zz}-\kappa_3(1+\mu_z^2))^2+(a^2+b^2)(1+\mu_z^2)\right]/P^{-2} \eqno{(4.9)}$$
by (4.6). Next, we study the condition under which the right hand side of the equation (4.9) is independent of $z$.
For the numerator and the denominator of the right hand side of (4.9), we have 
$$\bullet \ {\rm the \ numerator}=(P^{-1}\mu_{zz}-\kappa_3(1+\mu_z^2))^2+(a^2+b^2)(1+\mu_z^2)$$
$$=[(a^2+b^2)(1+c_2^2)+\{C^2(c_1^2-c_2^2-1)^2-c_1^2(a^2+b^2)\}/2+D^2(1+c_2^2)^2/2]$$
$$+\textstyle{\frac{1}{2}}\left[-C^2(c_1^2-c_2^2-1)^2+c_1^2(a^2+b^2)+D^2(1+c_2^2)^2\right]\cos(2\mu)
+CD(1+c_2^2)(c_1^2-c_2^2-1)\sin(2\mu).$$
$$\bullet \  {\rm the \ denomirator}=P^{-2}=(1/2)(C^2+D^2)+(1/2)(C^2-D^2)\cos(2\mu)+CD\sin(2\mu).$$
Comparing the coefficients of $\sin(2\mu)$ in the two equations above, the equation (4.9) should be equal to $(1+c_2^2)(c_1^2-c_2^2-1)$ if (4.9) is independent of $z$. Under the assumption that (4.9) is equal to $(1+c_2^2)(c_1^2-c_2^2-1)$, we further compare, respectively, the terms of $(x,y)$-function and the coefficients of $\cos(2\mu)$ in the two equations. Then, we obtain the equation $a^2+b^2=C^2(c_1^2-c_2^2-1)-D^2(c_2^2+1)$, which implies $\zeta=\kappa_3^2$. On the other hand, $\zeta=\kappa_3^2$ always holds. Hence, the equation $\zeta=\kappa_3^2$ implies that (4.9) is equal to the constant $(1+c_2^2)(c_1^2-c_2^2-1)$, i.e., we have $(1+\mu_z^2)^{2}(l^2+m^2)=(1+c_2^2)(c_1^2-c_2^2-1)(>0)$ on $U'$. 

Now, by the continuity of $(1+\mu_z^2)^{2}(l^2+m^2)$ and the connectivity of $U$, the equation $(1+\mu_z^2)^{2}(l^2+m^2)=(1+c_2^2)(c_1^2-c_2^2-1)$ is extended on the whole $U$. In particular, we have verified $(l^2+m^2)(x,y,z)\neq 0$ for any $(x,y,z)\in U$. Furthermore, the equation $(1+\mu_z^2)^{2}(l^2+m^2)=(1+c_2^2)(c_1^2-c_2^2-1)$ on $U$ implies that the equation (4.8) vanishes identically.

Since (4.8) vanishes identically, the equation $[\bs{B}_z/\sqrt{l^2+m^2}]_{z}=-\sqrt{l^2+m^2}\bs{B}$ holds on $U$ by (4.7), that is, $\bs{B}(z)$ moves on a unit circle in a plane. Here, we have $\sqrt{l^2+m^2}=\sqrt{(1+c_2^2)(c_1^2-c_2^2-1)}/(1+\mu_z^2)$.

In consequence, in the case (ii), the subcase (ii-1) (resp. (ii-2)) is the case that $f$ has the structure (a) (resp. the structure (b)) by the property of ${\bs B}(z)=-A(x,y,z)/\sqrt{1+\mu_z^2}$ in Lemma 4.1, and we know that $f$ in each case has the expression given, by the above argument.

(3) \underline{Suppose that all $(x,y)$-surfaces $f(x,y,z)$ with fixed $z$ lie on concentric $3$-spheres $H(z)$}. 

Then, we have $K(x,y,z)=K(z)\neq 0$ for all $(x,y,z)\in U$ and $f(x,y,z)=-(A/K)(x,y,z)$ by Lemma 4.1. Now, we take the derivative of both sides of $f(x,y,z)=-(A/K)(x,y,z)$ by $z$. Then, we have 
$$P=(1/K)\left[-(K_z/K)+P\kappa_3\mu_z\right], \ \ a=b\equiv 0, \ \ (K_z/K)\mu_z+(P\kappa_3-\mu_{zz})\equiv 0 \eqno{(4.10)}$$      
by $\partial X_{\gamma}/\partial z=P(aX_{\alpha}+bX_{\beta}+\kappa_3N)$. Hence, the function $P(x,y,z)$ depends only on $z$ by $a=b\equiv0$: $P(x,y,z)=:P(z)$. Furthermore, the first equation of (4.10) is equivalent to $K_z/K=(P^{-1})_z/P^{-1}$ by $K=-(P^{-1})_z+\kappa_3\mu_z$.  

In the case that $f$ has a Guichard net determined from $\mu\equiv 0$, we have $\kappa_3\equiv 0$ by the last equation in (4.10), and hence $X_{\gamma}$ is independent of $z$: $X_{\gamma}(x,y,z)=:X_{\gamma}(x,y)$. Furthermore, we have $(P^{-1})_z=-c_3P^{-1}$ with a constant $c_3\neq 0$ by $K=-(P^{-1})_z$ and $K_z/K=(P^{-1})_z/P^{-1}$. Then, we can determine $c_3$ as $c_3=\sqrt{1-c_1^2}~(>0)$ by $(P^{-1})_{zz}=(1-c_1^2)P^{-1}$ in (4.2). Hence, $P^{-1}(z)=(1/C)\exp(-c_3z)$ with a constant $C(>0)$ and $f(x,y,z)=X_{\gamma}(x,y)/K(z)=(1/c_3)P(z)X_{\gamma}(x,y)$ hold. 

In the case that $f$ has a Guichard net determined from $\mu_z\neq 0$, we have $K=c_4P^{-1}$ with a constant $c_4\neq 0$ by $K_z/K=(P^{-1})_z/P^{-1}$. Now, we determine the constant $c_4$. We firstly have $P^{-1}=(1/C)\sqrt{1+\mu_z^2}\exp(-\int_0^zc_4/(1+\mu_z^2)~dz)$ with a constant $C(>0)$ by $K=-(P^{-1})_z+\kappa_3\mu_z$ and $\kappa_3=P^{-1}\mu_{zz}-(P^{-1})_z\mu_z$ in (4.10). Differentiating this $P^{-1}$ yields
$$(P^{-1})_{zz}/P^{-1}=(c_4^2+\mu_{zz}^2)/(1+\mu_z^2)^2+\mu_z\mu_{zzz}/(1+\mu_z^2).$$  
On the other hand, we have $(P^{-1})_{zz}/P^{-1}=1-c_1^2\cos(2\mu)$ by (4.2). 
From these two equations for $(P^{-1})_{zz}/P^{-1}$, we have $c_4^2=(1+c_2^2)(1+c_2^2-c_1^2)~(>0)$ by direct calculation. Hence, we can determine $c_4$ as $c_4=\sqrt{(1+c_2^2)(1+c_2^2-c_1^2)}$. Then, we have 
$f(x,y,z)=-(1/K(z))A(x,y,z)=-(1/c_4)P(z)A(x,y,z)$.

Finally, we show $\kappa_3=\left((c_4\mu_z+\mu_{zz})/(1+\mu_z^2)\right)P^{-1}$. In fact, we have $\kappa_3\mu_z=c_4P^{-1}+(P^{-1})_z$ and $\kappa_3\mu_z=P^{-1}\mu_z\mu_{zz}-(P^{-1})_z\mu_z^2$ as above. Adding two of these equations gives the formula $\kappa_3$. 
In consequence, the case (3) has been verified.
\hspace{\fill}$\Box$\\

In the following Theorem 4.3 and its proof, we say that $f^*$ is a normal form with the structure (a), (b) or (c) if $f^*$ has such a structure by a parallel translation. In Theorem 4.3, we add the following condition (d) for $\varphi$ determining the Guichard net $g$ to clearly distinguish between the coordinates $(x,y)$ and $z$ of $f$:

(d) In $\varphi(x,y,z)=\lambda(x,y)+\mu(z)$, the $\lambda(x,y)$ is not a single variable function. 

\noindent
The condition (d) induces $\kappa_1\kappa_2\neq 0$. In fact, suppose $\kappa_1\equiv 0$, then $\kappa_3=-P^{-1}\tan\varphi$ holds. Hence, we have $(\kappa_3)_x=-(P^{-1})_x\tan\varphi-P^{-1}(\varphi_x/\cos^2\varphi)$. On the other hand, we have $(\kappa_3)_x=-(P^{-1})_x\tan\varphi$ by Lemma 2.2, which shows $\varphi_x\equiv 0$. In the same way, $\kappa_2\equiv0$ induces $\varphi_y\equiv0$.

In the proof, the inversion invariant $Inv^I(f^*)=K^*f^*+A^*$ and the dual invariant $Inv^D(f)_3=f^*+\sigma_3f+\kappa_3N$ play an important role. We use notations in Proposition 4.2. The description $\lambda^*(x,y)=\lambda(x,y)\pm\pi$ is settled on one of them under the condition that $\lambda^*(x,y)$ is contained in the interval $(-\pi,\pi)$. By the way, in the case that $f$ has the structure (a) and $\mu_z\neq 0$, the dual $f^*$ was also given in \cite[Section 2]{hsuy}: the dual pair $(f^+,f^-)(x,y,z)$ in that paper corresponds to our pair $(f, f^*)(x,y,z)$ with $C=1$ as $f^+=f$ and $f^-=-(1/2)f^*$ (cf. \cite[(2.10)]{su3}); hence, the difference between $f^+$ and $f$ (or $f^-$ and $f^*$) stems from homotheties.  \\

{\thm 4.3}. \ \ {\it Let a function $\varphi(x,y,z)=\lambda(x,y)+\mu(z)$ satisfy the condition (d) and determine a Guichard net $g$ in (1.1) on a connected open set $U\subset \R^3$: hence, $\lambda$ and $\mu$ satisfy (4.1). Suppose that $f(x,y,z)$ is a normal form in $\R^4$ with the first fundamental form $I_f=P^2g$: hence, $f(x,y,z)$ satisfies the conditions of Proposition 4.2 according to the structures (a), (b) and (c).  

Then, the dual $f^*(x,y,z)$ of $f(x,y,z)$ is a normal form with the same structure as $f(x,y,z)$. 
Let $\varphi^*(x,y,z)=\lambda^*(x,y)+\mu^*(z)$ be the function determining the Guichard net $g^*$ of $f^*(x,y,z)$. Then, we have $\mu^*\equiv 0$ if $\mu\equiv 0$, and $\mu^*_z\neq 0$ if $\mu_z\neq 0$.

If $f$ has the structure (a) or (c) and $\mu\equiv 0$, then $\varphi^*(x,y)=\lambda(x,y)\pm\pi$, and $f^*(x,y,z)$ and $f(x,y,z)$ are conformally equivalent: in the case (a), $f^*$ is expressed as $f^*(x,y,z)=(P^{-2}/2)(-f^1(x,y),Pz))$ for $f(x,y,z)=(f^1(x,y),Pz)$ in Proposition 4.2; in the case (c), $f^*$ is expressed as $f^*(x,y,z)=-\sigma_3(z)f(x,y,z)=-(1/(2c_3^2))\iota_0(f(x,y,z))$, where $\iota_0$ is the inversion $\iota_0(x):=x/\parallel x\parallel^2$. 

In the other cases, the dual $f^*(x,y,z)$ are not conformally equivalent to $f(x,y,z)$.  
In particular, we have the following facts in each case:

\underline{In the case that $f$ has the structure (a) and $\mu_z\neq 0$}: We have $\varphi^*(x,y,z)=(\lambda(x,y)\pm\pi)-\mu(z)$ and two $(x,y)$-surfaces of $f(x,y,z)$ and $f^*(x,y,z)$ with fixed $z$, respectively, lie in the parallel hyperplanes $H(z)$ and $H^*(z)$ perpendicular to the unit constant vector ${\bs B}=-A/\sqrt{1+\mu_z^2}$. In particular, for the constant $C$: $P^{-1}=C\cos\mu$ and $\kappa_3=-C\sin\mu$, and the unit normal vector field $N(x,y,z)$ of $f(x,y,z)$, the dual $f^*(x,y,z)$ is expressed as 
$$f^*(x,y,z)=-(C^2/2)f(x,y,z)-\kappa_3(z) N(x,y,z)+q(z){\bs B},$$
where $q(z)$ is a function satisfying $q'(z)=(P^{-1}\sqrt{1+\mu_z^2})(z)=Cc_1\cos^2\mu(z)$.

\underline{In the case that $f$ has the structure (b)}: Regardless of whether $\mu\equiv 0$ and $\mu_z\neq 0$, two $(x,y)$-surfaces of $f(x,y,z)$ and $f^*(x,y,z)$ with fixed $z$ lie in the same linear hyperplane $H(z)$ and the following two cases ($b_1$) and ($b_2$) do not occur: ($b_1$) $(\lambda^*)_x= \lambda_x$ and $(\lambda^*)_y=\lambda_y$ on $U$; ($b_2$) $(\lambda^*)_x=- \lambda_x$ and $(\lambda^*)_y=-\lambda_y$ on $U$ (the exact formulae of $\lambda^*$ are given in the proof). If $\mu_z\neq 0$, then $\mu^*(z)=-\mu(z)$ holds. In particular, for $f(x,y,z)$ given in Proposition 4.2, there is a $\R^2$-valued function $W(x,y)=(0,0,w^3(x,y),w^4(x,y))$ independent of $z$ such that $f^*(x,y,z)$ is given by 
$$f^*(x,y,z)=\sigma_3(x,y)f(x,y,z) +W(x,y).$$

\underline{In the case that $f$ has the structure (c) and $\mu_z\neq 0$}: We have $\lambda^*(x,y)= \lambda(x,y)$ and $\mu^*_z\neq \pm\mu_z$. In particular, we have $f^*(x,y,z)=-A^*(x,y,z)/K^*(z)=-(\sigma_3f+\kappa_3N)(x,y,z)$. } \\

{\it{Proof}.} \ Firstly, we give a few facts that are readily apparent. The function $K(x,y,z)$ for $f$ depends only on $z$ in every case, by Lemma 4.1. If $K(z)\equiv 0$, then $K^*(x,y,z)\equiv 0$ and $\varphi^*_z(x,y,z)=-\mu_z(z)$ by $K^*=-K/\sigma_3$ and Lemma 2.3. Hence, if $f$ has the structure (a) or (b), then we have $K^*\equiv 0$ and $\varphi^*_z(z)=-\mu_z(z)$. Furthermore, when $f$ with the Guichard net $g$ determined by $\mu\equiv 0$ has the structure (a) or (c), the hypersurfaces $f^*(x,y,z)$ and $f(x,y,z)$ are conformally equivalent by $\kappa_3\equiv 0$. In the argument below, we check carefully where the condition (d) is required.

Now, we study the dual $f^*$ of each hypersurface $f$ with the structure (a), (b) or (c), using Lemma 4.1 and Proposition 4.2. To do so, we separate the cases of (1) $\mu\equiv 0$ and (2) $\mu_z\neq 0$. 

(1) In the case that $g$ is determined by $\mu\equiv 0$: \underline{Suppose that $f$ has the structure (a)}. Then, since $\varphi^*_z=-\mu_z\equiv 0$ as mentioned above, we have $\varphi^*(x,y,z)=\lambda^*(x,y)$ by $\varphi^*_{xz}=\varphi^*_{yz}\equiv 0$. Furthermore, since $\sigma_1=\sigma_2=-\sigma_3=(\kappa_1\kappa_2)/2<0$ by (2.2) and $\kappa_3\equiv 0$, we have $\lambda^*(x,y)=\lambda(x,y)\pm\pi$ by $\cos\varphi^*=-\cos\varphi$ and $\sin\varphi^*=-\sin\varphi$ in (2.6). 
   
Since $X^*_{\gamma}=X_{\gamma}$ is the constant vector ${\bs B}$ and $K^*=-K/\sigma_3\equiv 0$, any $(x,y)$-surface $f^*(x,y,z)$ with $z$ lies in a hyperplane $H^*(z)$ perpendicular to $\bs B$, and hence all hyperplanes $H(z)$ and $H^*(z)$ are parallel, where $H(z)$ is a hyperplane containing an $(x,y)$-surface $f(x,y,z)$ with fixed $z$. Furthermore, since $P$ is constant and $\kappa_3\equiv 0$, the function $P^*$ is also constant by $P^*=\sigma_3P=[(P^{-2}+\kappa_3^2)/2]P=P^{-1}/2$ in (2.6). Now, all $z$-curves of both hypersurfaces $f(x,y,z)$ and $f^*(x,y,z)$ are lines in $\R^4$ by $\kappa_3=\kappa_3^*\equiv 0$ and $f$ is given by $f(x,y,z)=(f^1(x,y),Pz)\in H(z_0)\times \{t{\bs B}\}_{t\in \R}$. 
Hence, $f^*$ is expressed as $(P^{-2}/2)(-f^1(x,y),Pz)$ by $f^*_x=\sigma_1f_x=-(P^{-2}/2)f_x$, $f^*_y=\sigma_2f_y=-(P^{-2}/2)f_y$ and $f^*_z=\sigma_3f_z=(P^{-2}/2)f_z$.

\underline{Suppose that $f$ has the structure (b)}. Then, we have $\varphi^*(x,y,z)=\lambda^*(x,y)$ by $\varphi^*_z=-\mu_z\equiv 0$, as in the above case. The function $\lambda^*(x,y)$ satisfies the equation for $\lambda(x,y)$ in (4.1) by replacing the constant $c_1$ with $c_1^*$. Then, we have $\lambda^*_x=\lambda_x-(\kappa_1/\sigma_3)(P^{-1})_x$ and $\lambda^*_y=\lambda_y-(\kappa_2/\sigma_3)(P^{-1})_y$ by Lemma 2.3. In these equations, note that $\sigma_3=(P^{-2}+\kappa_3^2)/2$ and $\kappa_i \ (i=1,2)$ are independent of $z$, because $P$ and $\kappa_3$ are independent of $z$ by Proposition 4.2. In particular, since the case $\kappa_1(P^{-1})_x=\kappa_2(P^{-1})_y\equiv 0$ does not occur, at least either $(\lambda^*)_x\neq \lambda_x$ or $(\lambda^*)_y\neq \lambda_y$ is possible, that is, the case ($b_1$) does not occur. 

Now, we show that the case ($b_2$): $(\lambda^*)_x=- \lambda_x$ and $(\lambda^*)_y=- \lambda_y$ on $U$, also does not occur. Here, we use the condition (d). If $(\lambda^*)_x=- \lambda_x$ and $(\lambda^*)_y=- \lambda_y$ hold on $U$, then we have either $\lambda^*=-\lambda$ or $\lambda^*=-\lambda\pm \pi$ by (4.1). Hence, $f$ and $f^*$ are conformally equivalent, which implies either $\kappa_1\equiv0$ or $\kappa_2\equiv0$ by $\kappa_3\neq 0$ and contradicts (d). (Precisely, when we define a function $\theta_1=\theta_1(x,y)$ of $x$ and $y$ by $\cos\theta_1:=(\kappa_3^2-P^{-2})/(\kappa_3^2+P^{-2})$ and $\sin\theta_1:=2\kappa_3P^{-1}/(\kappa_3^2+P^{-2})$, we have $\lambda^*=\lambda-\theta_1$.)

Next, by Lemma 4.1, the vector ${\bs B}(z)=X_{\gamma}(x,y,z)$ depends only on $z$ and ${\bs B}(z)=X_{\gamma}(x,y,z)=X_{\gamma}^*(x,y,z)$ moves on the unit circle in a plane $L$. Hence, all $z$-curves $f^*(x,y,z)$ with fixed $(x,y)$ are circles, and  
every $(x,y)$-surface $f^*(x,y,z)$ with fixed $z$ lies in a hyperplane $H^*(z)$ perpendicular to ${\bs B}(z)=X_{\gamma}^*(x,y,z)$ by $K^*=0$. 
In consequence, we can choose $f^*(x,y,z)$ so that every two $(x,y)$-surfaces $f(x,y,z)$ and $f^*(x,y,z)$ with fixed $z$ lie in the same linear hyperplane $H(z)$ perpendicular to ${\bs B}(z)$. In particular, $f^*$ has the structure (b).
 
Now, let $f(x,y,z)$ be the hypersurface given in Proposition 4.2. Since $\sigma_3$ is independent of $z$, we have $f^*_z=\sigma_3f_z=(\sigma_3f)_z$. Hence, there is a $\R^4$-valued function ${\bs W}(x,y)$ independent of $z$ such that $f^*(x,y,z)=\sigma_3(x,y)f(x,y,z)+{\bs W}(x,y)$. Then, ${\bs W}(x,y)$ takes values in the plane $L^{\perp}$ perpendicular to $L$. In fact, we have ${\bs W}_x=(\sigma_1-\sigma_3)f_x-(\sigma_3)_xf$ and ${\bs W}_y=(\sigma_2-\sigma_3)f_y-(\sigma_3)_yf$ from $f^*=\sigma_3f+{\bs W}$. Furthermore, by Lemmata 2.2 and 2.4, we have 
$${\bs W}_x=\kappa_2(\kappa_1-\kappa_3)f_x-\kappa_2(\kappa_3)_xf=\kappa_2\tan\varphi(P^{-1}f_x+(P^{-1})_xf)=\kappa_2\tan\varphi[P^{-1}f]_x,$$
$${\bs W}_y=\kappa_1(\kappa_2-\kappa_3)f_y-\kappa_1(\kappa_3)_yf=-\kappa_1\cot\varphi(P^{-1}f_y+(P^{-1})_yf)=-\kappa_1\cot\varphi[P^{-1}f]_y.$$ 
In these two equations, the planar $L$-component of $P^{-1}f$ depends only on $z$. Hence, we can take ${\bs W}$ so that ${\bs W}(x,y)\in L^{\perp}$. Note that the integrability condition ${\bs W}_{xy}={\bs W}_{yx}$ is guaranteed by the existence of $f^*$.  

\underline{Suppose that $f^*$ has the structure (c)}. Then, we have $\varphi^*_z=-\mu_z=0$ by $\kappa_3\equiv 0$ and Lemma 2.3. Furthermore, we have $\varphi^*(x,y,z)=\varphi(x,y,z)\pm\pi$ by $\kappa_3\equiv 0$, as in the case (a). Thus, we have $\varphi^*(x,y,z)=\lambda(x,y)\pm\pi$.

Now, $\sigma_3=P^{-2}/2=(2C^2)^{-1}\exp(-2c_3z)$ depends only on $z$ by $P=C\exp(c_3z)$ with a constant $C$ and $\kappa_3\equiv 0$. Hence, $f^*$ can be taken as $f^*(x,y,z)=-\sigma_3(z)(X_{\gamma}(x,y)/K(z))$ by $f(x,y,z)=X_{\gamma}(x,y)/K(z)(=-A(x,y)/K(z))$ and Theorem 3.5-(1), which is rewritten as $f^*(x,y,z)=-A^*(x,y)/K^*(z)$ by $X_{\gamma}^*(x,y)=X_{\gamma}(x,y)$. In particular, $f^*$ has the structure (c). Furthermore, since $K=-(P^{-1})_z$, we also have the expression $f^*=-(2c_3^2)^{-1}(\iota_0\circ f)$ by $f=(C/c_3)\exp(c_3z)X_{\gamma}$.

(2) In the case that $g$ is determined by $\mu_z\neq 0$: \underline{Suppose that $f$ has the structure (a)}. Then, $\varphi^*_z=-\mu_z$ holds as mentioned at the beginning of this proof. Hence, we have $\varphi^*(x,y,z)=\lambda^*(x,y)+\mu^*(z)$, where $\sin^2\mu^*=\sin^2\mu.$ 
By $P^{-1}(z)=C\cos\mu(z)$ and $\kappa_3(z)=-C\sin\mu(z)$ with a non-zero constant $C$, we have $\varphi^*_x=\varphi_x$, $\varphi^*_y=\varphi_y$ by Lemma 2.3. 
More precisely, we have $\varphi^*(x,y,z)=(\lambda(x,y)\pm\pi)-\mu(z)$ by (2.6), since 
$$\sigma_1=-\frac{C^2}{2}\frac{\cos(\lambda-\mu)}{\cos\varphi}, \ \ \sigma_2=-\frac{C^2}{2}\frac{\sin(\lambda-\mu)}{\sin\varphi}, \ \  \sigma_3=\frac{P^{-2}+\kappa_3^2}{2}=\frac{C^2}{2}$$
by (2.2). 

Since $A^*=-X_{\gamma}+\mu_zN=A$ by $X_{\gamma}^*=X_{\gamma}$ and $N^*=-N$, we have ${\bs B}^*=-A^*/\sqrt{1+\mu_z^2}=-A/\sqrt{1+\mu_z^2}={\bs B}$, which is a constant vector. Furthermore, since $X_{\alpha}^*=X_{\alpha}$ and $X_{\beta}^*=X_{\beta}$ are perpendicular to ${\bs B}$, any $(x,y)$-surfaces $f^*(x,y,z)$ with fixed $z$ lie in a hyperplane $H^*(z)$ perpendicular to ${\bs B}$. Hence, all hyperplanes $H^*(z)$ and $H(z)$ are parallel. Then, if $C\neq \sqrt{2}$, then $H^*(z)\neq H(z)$ by $f^*_z=\sigma_3f_z=(C^2/2)f_z$.
In particular, $f^*$ has the structure (a). 

Now, we verify that $f^*(x,y,z)$ is expressed in the given form. We have $[f^*+\sigma_3f+\kappa_3N]_x=[f^*+\sigma_3f+\kappa_3N]_y={\bs 0}$ and $[f^*+\sigma_3f+\kappa_3N]_z=-P^{-1}A$ by $(\kappa_3)_x=(\kappa_3)_y=0$ and Theorem 3.5, which shows our assertion.

\underline{Suppose that $f$ has the structure (b)}. Then, we have $\varphi^*_z=-\mu_z$. Hence, we have $\varphi^*(x,y,z)=\lambda^*(x,y)+\mu^*(z)$. Since $P^{-1}(x,y,z)=C(x,y)\cos\mu(z)+D(x,y)\sin\mu(z)$ and $\kappa_3(x,y,z)=-C(x,y)\sin\mu(z)+D(x,y)\cos\mu(z)$, we have $\sigma_3=(1/2)(C^2+D^2)$ and 
$$\frac{\sigma_1}{\sigma_3}=\frac{1}{\cos\varphi}\left[\frac{-C^2+D^2}{C^2+D^2}\cos(\lambda-\mu)+2\frac{CD}{C^2+D^2}\sin(\lambda-\mu)\right], $$ 
$$\frac{\sigma_2}{\sigma_3}=\frac{1}{\sin\varphi}\left[\frac{-C^2+D^2}{C^2+D^2}\sin(\lambda-\mu)-2\frac{CD}{C^2+D^2}\cos(\lambda-\mu)\right].$$ 
Thus, we have $\lambda^*(x,y)=\lambda(x,y)-\theta_2(x,y)$ and $\mu^*(z)=-\mu(z)$, where $\theta_2(x,y)$ is determined by 
$\cos\theta_2:=(-C^2+D^2)/(C^2+D^2)$ and $\sin\theta_2:=2CD/(C^2+D^2)$. Furthermore, in the same way as the case $\mu\equiv 0$, we can prove the fact that the cases ($b_1$) and ($b_2$) do not occur by the condition (d).

Since $X^*_{\alpha}=X_{\alpha}$ and $X^*_{\beta}=X_{\beta}$ are perpendicular to ${\bs B}(z)=-(A/\sqrt{1+\mu_z^2})(x,y,z)$, every $(x,y)$-surface $f^*(x,y,z)$ 
with fixed $z$ lies in a hyperplane $H^*(z)$ perpendicular to ${\bs B}(z)$. Furthermore, since $f(x,y,z)$ has the structure (b), the vector ${\bs B}(z)$ moves on the unit circle in a plane $L$ and every $(x,y)$-surface $f(x,y,z)$ with fixed $z$ lie in a linear hyperplane $H(z)$ perpendicular to ${\bs B}(z)$. When we note that $f^*$ is determined uniquely up to a parallel translation, we can choose $f^*(x,y,z)$ so that $H^*(z)=H(z)$, since ${\bs B}(z)$ moves on the unit circle of a plane and we have ${\bs B}^*(z)={\bs B}(z)$ and $X^*_{\gamma}=X_{\gamma}$.

In consequence, two $(x,y)$-surfaces of $f(x,y,z)$ and $f^*(x,y,z)$ with fixed $z$ lie in the same linear hyperplane $H(z)$ perpendicular to ${\bs B}(z)={\bs B}^*(z)$. In particular, $f^*$ has the structure (b). Furthermore, since $\sigma_3$ is independent of $z$, $f^*(x,y,z)$ is given in the same form, with respect to $f(x,y,z)$ in this case, as for $\mu\equiv 0$.

\underline{Suppose that $f$ has the structure (c)}. In this case, $P^{-1}$ and $\kappa_3(\neq 0)$ are functions of only $z$, and hence $\sigma_3$ and $K(\neq 0)$ also depend only on $z$. Then, by Lemma 2.3, the function $\varphi^*_z=-\mu_z+(\kappa_3/\sigma_3)K~(\neq -\mu_z)$ also depend only on $z$. 
Furthermore, we have $\varphi^*_x=\lambda_x$ and $\varphi^*_y=\lambda_y$ by $(P^{-1})_x=(P^{-1})_y\equiv 0$. Hence, we have either $\lambda^*(x,y)=\lambda(x,y)$ or $\lambda^*(x,y)=\lambda(x,y)\pm\pi$. Then, we also have $\mu^*_z(z)\neq \mu_z(z)$ (note that $\mu^*_z=\mu_z$ implies either $\mu^*= \mu$ or $\mu^*= \mu\pm\pi$ by $\mu^*_z\neq -\mu_z$). In fact, if $\mu_z^*(z)=\mu_z(z)$ holds, then we have either $\varphi^*(x,y,z)=\varphi(x,y,z)$ or $\varphi^*(x,y,z)=\varphi(x,y,z)\pm \pi$. 
Hence, $f$ and $f^*$ are conformally equivalent and then we have either $\kappa_1\equiv0$ or $\kappa_2\equiv 0$ by $\kappa_3\neq 0$, which contradicts the fact that every $(x,y)$-surface $f(x,y,z)$ with fixed $z$ lies on a $3$-sphere. Thus, we have $\mu^*_z(z)\neq\mu_z(z)$. Therefore, $\varphi^*$ can be written as $\varphi^*(x,y,z)=\lambda(x,y)+(\mu-\theta_3)(z)$ with the following function $\theta_3(z)$ of $z$ satisfying $(\theta_3)_z\neq0$ and $(\theta_3)_z\neq 2\mu_z$: the $\theta_3=\theta_3(z)$ is determined by $\cos\theta_3:=(\kappa_3^2-P^{-2})/(\kappa_3^2+P^{-2})$ and $\sin\theta_3:=(2\kappa_3P^{-1})/(\kappa_3^2+P^{-2})$.  
  
Now, we have $[f^*+\sigma_3f+\kappa_3N]_x=[f^*+\sigma_3f+\kappa_3N]_y=[f^*+\sigma_3f+\kappa_3N]_z={\bf 0}$ by Theorem 3.5-(1) and $Kf+A={\bf 0}$. Hence, we can determine $f^*$ by $f^*(x,y,z)=-(\sigma_3f+\kappa_3N)(x,y,z)$. Then, by $K^*=-K/\sigma_3$ and
$$A^*=-X_{\gamma}+\mu^*_zN^*=A-(\kappa_3/\sigma_3)KN=-K[f+(\kappa_3/\sigma_3)N]=K^*[\sigma_3f+\kappa_3N],$$ 
$f^*$ is also expressed as $f^*(x,y,z)=-A^*(x,y,z)/K^*(z)$, and hence $f^*$ has the structure (c). 
\hspace{\fill}$\Box$\\

By Proposition 4.2 and Theorem 4.3, the dual invariants $Inv^D(f)_3$ of normal forms $f$ satisfy the following facts: 

$\bullet$ If $f$ has the structure (a), then $Inv^D(f)_3(x,y,z)=q(z){\bs B}$ with a unit constant ${\bs B}$. In particular, when $f$ has a Guichard net determined by $\mu\equiv 0$, we have $q(z)=P^{-1}z$ with a positive constant $P$.

$\bullet$ If $f$ has the structure (b), then $[Inv^D(f)_3]_z(x,y,z)=r(x,y,z){\bs B}(z)$ with a vector ${\bs B}(z)$ of $z$ moving on a unit circle in a plane. Here, $r(x,y,z)=:P^{-1}(x,y)$ if $\mu\equiv 0$, and $r(x,y,z)=:(P^{-1}\sqrt{1+\mu_z^2})(x,y,z)$ if $\mu_z\neq 0$. 

$\bullet$ If $f$ has the structure (c), then $Inv^D(f)_3(x,y,z)= {\bs 0}$.

Let $g$ be a Guichard net determined by a function $\varphi(x,y,z)$ satisfying $\varphi_{xz}=\varphi_{yz}\equiv 0$. Let $f$ be a normal form $f$ with the first fundamental form $I_f=P^2g$. 
In Theorem 4.3, we have obtained an exact dual $f^*$ (under the condition (d)) by paying attention to the structures (a)-(c) of $f$. However, even the conformal transformations $\iota_q\circ f$ of those $f$ would be difficult to get explicit duals in any form, because the functions $(\varphi^q)^*$ for $(\iota_q\circ f)^*$ do not satisfy $(\varphi^q)^*_{xz}=(\varphi^q)^*_{yz}\equiv 0$. \\

{\bf 5. Two kind of Approximations of dual hypersurfaces} 

In this section, let $g=\cos^2\varphi(dx)^2+\sin^2\varphi(dy)^2+(dz)^2$ be a Guichard net defined on a connected open set $U(\subset \R^3)$ (not limited to the Guichard nets considered in the previous section). 
Let $f(x,y,z)$ be a generic conformally flat hypersurface in $\R^4$ with the first fundamental form $I_f=P^2g$. Same as previous sections, let $\kappa_i(x,y,z)$ and $\sigma_i(x,y,z)$ ($i=1,2,3$) be the principal curvatures and the eigenvalues of the Schouten $(1,1)$-tensor $S$ of $f$, respectively. The dual $f^*(x,y,z)$ of $f(x,y,z)$ is defined by $df^*=df\circ S$ on $U$, and then it is regular only on $U^*=\{(x,y,z)\in U| \ (\sigma_1\sigma_2\sigma_3)(x,y,z)\neq 0\}$. The dual invariant $Inv^D(f)_3=f^*+\sigma_3f+\kappa_3N$ of $f$ satisfies the following equations by Theorem 3.5:
$$
[f^*+\sigma_3f+\kappa_3N]_x=(\kappa_3)_x(\kappa_2f+N), \ \  
[f^*+\sigma_3f+\kappa_3N]_y=(\kappa_3)_y(\kappa_1f+N).\eqno{(5.1)}$$
Recall that these equations are obtained from $(\sigma_3)_x=(\kappa_3)_x\kappa_2$ and $(\sigma_3)_y=(\kappa_3)_y\kappa_1$. Hence, the two equations in (5.1) are rewrote as $f^*_x=-\sigma_3f_x-\kappa_3N_x$ and $f^*_y=-\sigma_3f_y-\kappa_3N_y$.
In this section, using these equations, we define two kinds of approximations of $(x,y)$-surfaces $f^*(x,y,z)$ with fixed $z$ and further, using the following lemma for the dual invariants $Inv^D(f)_i=f^*+\sigma_if+\kappa_iN \ (i=1,2)$, we approximate $z$-curves of the dual $f^*(x,y,z)$ with fixed $(x,y)$. Note that the formulae of $Inv^D(f)_1$ are also rewrote as $f^*_y=-\sigma_1f_y-\kappa_1N_y$ and $f^*_z=-\sigma_1f_z-\kappa_1N_z$, for example. \\

{\lem 5.1}. \ \ {\it We have the following equations: $(\sigma_1)_y=(\kappa_1)_y\kappa_3$, $(\sigma_1)_z=(\kappa_1)_z\kappa_2$, $(\sigma_2)_x=(\kappa_2)_x\kappa_3$, $(\sigma_2)_z=(\kappa_2)_z\kappa_1$ and 
$$[f^*+\sigma_1f+\kappa_1N]_y=(\kappa_1)_y(\kappa_3f+N), \ \ [f^*+\sigma_1f+\kappa_1N]_z=(\kappa_1)_z(\kappa_2f+N),$$
$$[f^*+\sigma_2f+\kappa_2N]_x=(\kappa_2)_x(\kappa_3f+N), \ \ [f^*+\sigma_2f+\kappa_2N]_z=(\kappa_2)_z(\kappa_1f+N).$$
Here, we have
$$(\kappa_1)_y=[(P^{-1})_y+P^{-1}\varphi_y\tan\varphi](\sin\varphi\cos\varphi)^{-1}, \ \ (\kappa_1)_z=[(P^{-1})_z+P^{-1}\varphi_z\tan\varphi]\tan\varphi,$$ 
$$(\kappa_2)_x=[-(P^{-1})_x+P^{-1}\varphi_x\cot\varphi](\sin\varphi\cos\varphi)^{-1}, \ \ (\kappa_2)_z=[-(P^{-1})_z+P^{-1}\varphi_z\cot\varphi]\cot\varphi.$$}

{\it{Proof}.} \ In this proof, only some formulae related to $[f^*+\sigma_1f+\kappa_1N]_z$ are verified, since the other formulae are verified in the same way.   

The equation $(\sigma_1)_z=(\kappa_1)_z\kappa_2$ follows from the definition of $\sigma_1$ in (2.1) and the equation $(\kappa_1)_z(\kappa_3-\kappa_2)+(\kappa_2)_z(\kappa_1-\kappa_3)+(\kappa_3)_z(\kappa_1-\kappa_2)=0$ in (\cite{ca}, \cite{la}, \cite[Equations (2.2.9)]{su4}). By $\kappa_1=P^{-1}\tan\varphi+\kappa_3$ in (1.3) and $(\kappa_3)_z=-P^{-1}\varphi_z$ in Lemma 2.2, we have the formula for $(\kappa_1)_z$. Finally, the formula for $[f^*+\sigma_1f+\kappa_1N]_z$ is obtained from $\sigma_1+\sigma_3-\kappa_3\kappa_1=0$. 
\hspace{\fill}$\Box$\\

We construct two kinds of approximations $(f^*)^{\bar{\delta}_n}(x,y,z)$ and $(f^*)^{\underline{\delta}_n}(x,y,z)$ of the dual $f^*(x,y,z)$ for each positive integer $n$, focusing on the following two properties (P1) and (P2) of $f^*(x,y,z)$: 

(P1) The hypersurface $f^*(x,y,z)$ is a one parameter family of $(x,y)$-curvature surfaces with parameter $z$. That is, $f^*(x,y,z)$ is made of $(x,y)$-surfaces. 

(P2) Every $(x,y)$-curvature surface $f^*(x,y,z)$ with fixed $z$ is a one parameter family of principal curvature $x$-lines (or $y$-lines). Hence, each $(x,y)$-surface $f^*(x,y,z)$ can be regarded as made of $x$- (or $y$-) curves only.

Now, let $Q:=[x_0,x_e]\times[y_0,y_e]\times[z_0,z_e]$ be a compact cube in $U^*$ and $a:=x_e-x_0=y_e-y_0=z_e-z_0>0$. For an arbitrarily fixed $n>0$, let 
$$x_0<x_1<\cdots <x_{n}=x_e, \ \ y_0<y_1<\cdots <y_{n}=y_e, \ \ z_0<z_1<\cdots <z_{n}=z_e$$
be the divisions of three intervals of equal length $\delta_n:=a/n$. Let $L_n$ be the $3$-dimensional lattice in $Q$ consisting of the coordinate line segments passing through all points $(x_i,y_j,z_k)$. Note that $f(L_n)$ is a $3$-dimensional net in the hypersurface $f$. As described in Steps 1 and 2 below, for each $(x,y)$-surface $f^*(x,y,c_0)$ with an arbitrarily fixed $c_0:=z_k$, we firstly focus on one of the $x$-curves and the $y$-curves on $f^*(x,y,c_0)$ by (P2) and construct two kind of approximations defined on certain subsets of $L_n\cap (\R^2\times\{c_0\})$. The construction will proceed in the following steps.

{\bf Step 1}. In this Step, the $x$-curves ${\bs u}^{\delta_n}_i(x,y_j,z_k)$ and ${\bs u}^{\delta_n}(x,y_j,z_k)$ given in (1.5) and (1.6) of the introduction are defined again and the equations associated with them are verified. 
That is, for an $(x,y)$-curvature surface $f^*(x,y,c_0)$ with fixed $c_0=z_k$, we approximate its principal curvature $x$-lines $f^*(x,y_j,c_0)-f^*(x_0,y_j,c_0)$ on $[x_0,x_e]\times\{(y_j,c_0)\}$ (resp. its principal curvature $y$-lines $f^*(x_i,y,c_0)-f^*(x_i,y_0,c_0)$ on $\{x_i\}\times[y_0,y_e]\times\{c_0\}$) by the curves $\bs{u}^{\delta_n}(x, y_j,c_0)$ satisfying $\bs{u}^{\delta_n}(x_0, y_j,c_0)={\bf 0}$ (resp. curves $\bs{v}^{\delta_n}(x_i, y,c_0)$ satisfying $\bs{v}^{\delta_n}(x_i, y_0,c_0)={\bf 0}$). 
In the construction of these approximate $x$-curves $\bs{u}^{\delta_n}(x, y_j,c_0)$ and $y$-curves $\bs{v}^{\delta_n}(x_i, y,c_0)$, the equations $f^*_x=-\sigma_3f_x-\kappa_3N_x$ and $f^*_y=-\sigma_3f_y-\kappa_3N_y$ in (5.1) are used, respectively. 

Here, we say that the curve $\bs{u}^{\delta_n}(x, y_j,c_0)$ with fixed $(y_j,c_0)$ approximates the curve $f^*(x,y_j,c_0)-f^*(x_0,y_j,c_0)$, if there is a constant $C$ depending only on $f$ and $Q$ such that $\parallel(f^*-\bs{u}^{\delta_n})(x,y_j,c_0)-f^*(x_0,y_j,c_0)\parallel<C/n$ holds for $x\in [x_0,x_e]$. 

In this step, since $c_0=z_k$ is fixed, we denote these curves $\bs{u}^{\delta_n}(x, y_j,c_0)$ and $\bs{v}^{\delta_n}(x_i, y,c_0)$ by $\bs{u}^{\delta_n}(x, y_j)$ and $\bs{v}^{\delta_n}(x_i, y)$, respectively.

{\bf Step 2}. Under the consideration of (P2), we approximate each $(x,y)$-surface $f^*(x,y,c_0)$ with fixed $c_0=z_k$, by two surfaces $S^{\bar{\delta}_n}(x,y,c_0)$ and $S^{\underline{\delta}_n}(x,y,c_0)$. 
The surface $S^{\bar{\delta}_n}(x,y,c_0)$ is created from the $y$-curve $\bs{v}^{\delta_n}(x_0,y,c_0)$ and translated $x$-curves $\bs{u}^{\delta_n}(x,y_j,c_0)$ ($0\leq j\leq n$), by attaching to the $y$-curve $\bs{v}^{\delta_n}(x_0,y,c_0)$ these $x$-curves under the conditions $\bs{v}^{\delta_n}(x_0,y_j,c_0)=\bs{u}^{\delta_n}(x_0,y_j,c_0)$. Hence, the domain of $S^{\bar{\delta}_n}(x,y,c_0)$ is the set $(\{x_0\}\times[y_0,y_e]\times\{c_0\})\cup_{j=0}^{n} ([x_0,x_e]\times\{(y_j,c_0)\})$. 

On the other hand, the surface $S^{\underline{\delta}_n}(x,y,c_0)$ is created from the $x$-curve $\bs{u}^{\delta_n}(x,y_0,c_0)$ and translated $y$-curves $\bs{v}^{\delta_n}(x_i,y,c_0)$ ($0\leq i\leq n$), by attaching to the $x$-curve $\bs{u}^{\delta_n}(x,y_0,c_0)$ these $y$-curves under the conditions $\bs{u}^{\delta_n}(x_i,y_0,c_0)=\bs{v}^{\delta_n}(x_i,y_0,c_0)$. Hence, the domain of $S^{\underline{\delta}_n}(x,y,c_0)$ is the set $([x_0,x_e]\times\{(y_0,c_0)\})\cup_{i=0}^{n} (\{x_i\}\times[y_0,y_e]\times\{c_0\})$.       
Then, the surfaces $S^{\bar{\delta}_n}(x,y,c_0)$ and $S^{\underline{\delta}_n}(x,y,c_0)$ will approximate the curvature $(x,y)$-surfaces $f^*(x,y,c_0)-f^*(x_0,y_0,c_0)$ on the respective defined domains.

In this step, since $c_0=z_k$ is fixed, we denote these surfaces $S^{\bar{\delta}_n}(x,y,c_0)$ and $S^{\underline{\delta}_n}(x,y,c_0)$ by $S^{\bar{\delta}_n}(x,y)$ and $S^{\underline{\delta}_n}(x,y)$, respectively. In Step 3, we construct the approximate hypersurfaces $(f^*)^{\bar\delta_n}(x,y,z)$ and $(f^*)^{\underline{\delta}_n}(x,y,z)$ of $f^*(x,y,z)$ for each $n$.

{\bf Step 3}. In order to construct the approximation $(f^*)^{\bar\delta_n}(x,y,z)$, we use the equation $f^*_z=-\sigma_2f_z-\kappa_2N_z$ in Lemma 5.1. Using the equation, in the same way as the construction of ${{\bs u}^{\delta_n}}(x,y_j,c_0)$ in Step 1, 
we firstly construct a $z$-curve $\bs{w}^{\delta_n}(x_0,y_0,z)$ on $\{(x_0,y_0)\}\times [z_0,z_e]$ under the condition $\bs{w}^{\delta_n}(x_0,y_0,z_0)={\bf 0}$ as an approximation of $f^*(x_0,y_0,z)-f^*(x_0,y_0,z_0)$. 
Next, by the conditions $\bs{w}^{\delta_n}(x_0,y_0,z_k)=S^{\bar{\delta}_n}(x_0,y_0,z_k)$, we attach to the curve $\bs{w}^{\delta_n}(x_0,y_0,z)$ the family of surfaces $\{S^{\bar{\delta}_n}(x,y,z_k)\}_k$ by parallel translation. Thus, we obtain the hypersurface $(f^*)^{\bar{\delta}_n}(x,y,z)$ defined on subset of the lattice $L_n$. It will be the approximations of $f^*(x,y,z)-f^*(x_0,y_0,z_0)$ on the defined domain, by (P1).

On the other hand, when constructing $(f^*)^{\underline{\delta}_n}(x,y,z)$, instead of $\bs{w}^{\delta_n}(x_0,y_0,z)$ above, a $z$-curve $\bs{z}^{\delta_n}(x_0,y_0,z)$ on $\{(x_0,y_0)\}\times [z_0,z_e]$ constructed from the equation $f^*_z=-\sigma_1f_z-\kappa_1N_z$ is used as an approximation of the $z$-curve $f^*(x_0,y_0,z)-f^*(x_0,y_0,z_0)$. \\[-3mm]

In the next section 6, 
for every $(i,j,k)$, we shall construct an approximate $y$-curve connecting adjacent $x$-curves $S^{\bar{\delta}_n}(x,y_j,z_k)$ and  $S^{\bar{\delta}_n}(x,y_{j+1},z_k)$, at $x=x_i$. Then, $(f^*)^{\bar\delta_n}(x,y,z)$ and $(f^*)^{\underline{\delta}_n}(x,y,z)$ will be approximate discrete hypersurfaces of $f^*(x,y,z)$ on $Q$. 
In fact, the definition domains of the approximations $(f^*)^{\bar\delta_n}(x,y,z)$ and $(f^*)^{\underline{\delta}_n}(x,y,z)$ are extended to the whole lattice $L_n$, since the remaining approximate $x$-curves in $S^{\underline{\delta}_n}(x,y,z_k)$ and $z$-curves in $(f^*)^{\bar{\delta}_n}(x,y,z)$ and $(f^*)^{\underline{\delta}_n}(x,y,z)$ can be constructed in the same way: so that this process would work, we have taken the different $z$-approximations in the constructions of $f^{\bar\delta_n}(x,y,z)$ and $f^{\underline{\delta}_n}(x,y,z)$. That is, in the construction of $z$-curves $(f^*)^{\bar\delta_n}(x_i,y_j,z)$ with fixed $(x_i,y_j)$, we focus on $(x,z)$-surfaces $f^*(x,y_j,z)$ with fixed $y_j$ and these $z$-curves are defined so that each $(x,z)$-surface $(f^*)^{\bar\delta_n}(x,y_j,z)$ is an approximate discrete surface of $f^*(x,y_j,z)$. On the other hand, in the construction of $z$-curves $(f^*)^{\underline{\delta}_n}(x_i,y_j,z)$ with fixed $(x_i,y_j)$, we focus on $(y,z)$-surfaces $f^*(x_i,y,z)$ with fixed $x_i$.\\[-3mm]

Now, since the the first fundamental form $I_f=P^2(\cos^2\varphi(dx)^2+\sin^2\varphi(dy)^2+(dz)^2)$ of $f$ is positive definite on the compact cube $Q$, the functions $1/\sin\varphi$ and $1/\cos\varphi$ are bounded on $Q$. 
Adding this fact, there are positive constants $C_1$ and $C_2$ such that the following inequalities hold on $Q$: 
$$\begin{array}{ll}
C_1^{-1}<| P(x,y,z)|<C_1, \ \ |\kappa_i(x,y,z)|<C_1, \ \ |(\kappa_i)_x(x,y,z)|<C_2, \\[2mm] 
 \ \ |(\kappa_i)_y(x,y,z)|<C_2, \ \ |(\kappa_i)_z(x,y,z)|<C_2 \ \ for \ i=1,2,3.  
\end{array} \eqno{(5.2)}$$  
With these constants $C_i$, we further have $\parallel f_x(x,y,z)\parallel<C_1, \ \parallel f_y(x,y,z)\parallel<C_1, \ \parallel f_z(x,y,z)\parallel <C_1$ and
$$\begin{array}{ll}
|(\sigma_i)_x(x,y,z)|<C_1C_2 \ for \ i=2,3, \ \ |(\sigma_j)_y(x,y,z)|<C_1C_2 \ for \ j=1,3,\\[2mm] 
|(\sigma_k)_z(x,y,z)|<C_1C_2 \ for \ k=1,2
\end{array} \eqno{(5.3)}$$ 
by $\parallel f_x\parallel=|P\cos\varphi|, \ \parallel f_y\parallel=|P\sin\varphi|$, $\parallel f_z\parallel=|P|$ and $(\sigma_2)_x=(\kappa_2)_x\kappa_3, \ (\sigma_3)_x=(\kappa_3)_x\kappa_2$, $(\sigma_1)_y=(\kappa_1)_y\kappa_3$, $(\sigma_3)_y=(\kappa_3)_y\kappa_1$, $(\sigma_1)_z=(\kappa_1)_z\kappa_2$, $(\sigma_2)_z=(\kappa_2)_z\kappa_1$.   

In the Steps 1 and 2 below, we fix a positive integer $n$ and $c_0\in\{z_0,z_1,\cdots,z_n\}$ arbitrarily. 
We denote 
$$E:=[x_0,x_e]\times[y_0,y_e]\times\{c_0\}, \ \ f^0(x,y):=f(x,y,c_0), \ \ (f^*)^0(x,y):=f^*(x,y,c_0),$$ 
$$N^0(x,y):=N(x,y,c_0), \  \bar P(x,y):=P(x,y,c_0), \  \bar\kappa_i(x,y):=\kappa_i(x,y,c_0), \  \bar\sigma_i(x,y):=\sigma_i(x,y,c_0).$$
Hence, the surfaces $f^0(x,y)$ and $(f^*)^0(x,y)$ are the curvature surfaces of $f(x,y,z)$ and $f^*(x,y,z)$, respectively, and then we identify $E$ with the square $[x_0,x_e]\times[y_0,y_e]\subset \R^2$. 
We construct an approximation of $(f^*)^0(x,y)$ on $E$, for each $n$. Here, $x_0<x_1<\cdots <x_n=x_e$ and $y_0<y_1<\cdots <y_n=y_e$ are the divisions of the intervals $[x_0,x_e]$ and $[y_0,y_e]$ of equal length $\delta_n=a/n$.

{For Step 1}. For each $y=y_j$ ($0\leq j\leq n$), we define a curve $\bs{u}^{\delta_n}_i(x,y_j)(:=\bs{u}^{\delta_n}_i(x,y_j,c_0))$ on $[x_i,x_{i+1}]\times\{y_j\}$ ($0\leq i\leq n-1$) by 
$$\bs{u}^{\delta_n}_i(x,y_j):=-\bar\sigma_3(x_{i},y_j)\left[f^0(x,y_j)\right]_{x_i}^x-\bar\kappa_3(x_i,y_j)\left[N^0(x,y_j)\right]_{x_i}^x, \eqno{(5.4)}$$
where $\left[f^0(x,y_j)\right]_{x=x_i}^x=f^0(x,y_j)-f^0(x_i,y_j)$.
Note that $\bs{u}^{\delta_n}_i(x,y_j)$ on $[x_i,x_{i+1}]$ is a kind of parallel curve of the curve $f^0(x,y_j)-f^0(x_i,y_j)$.     
Now, we have $\bs{u}^{\delta_n}_i(x_i,y_j)={\bf 0}$ and the tangent vectors of two curves $\bs{u}^{\delta_n}_i(x,y_j)$ and $(f^*)^0(x,y_j)$ are parallel on $[x_i,x_{i+1}]$ by $(\nabla'_{\partial/\partial x}N^0)(x,y_j)=-(\bar\kappa_1f^0_x)(x,y_j)$ and $(f^*)^0_x(x,y_j)=\sigma_1f^0_x(x,y_j)$. 
In particular, we have the following equations:
$$(\partial\bs{u}^{\delta_n}_i/\partial x)(x_i{\small +0},y_j)=-\left[\bar\sigma_3-\bar\kappa_3\bar\kappa_1\right](x_i,y_j)f^0_x(x_i,y_j)=\left(\bar\sigma_1f^0_x\right)(x_i,y_j)=(f^*)^0_x(x_i,y_j),\leqno{ \ \ \ \bullet}$$
$$(\partial\bs{u}^{\delta_n}_{i}/\partial x)(x_{i+1}{\small -0},y_j)=-\left[\bar\sigma_3(x_{i},y_j)-\bar\kappa_3(x_{i},y_j)\bar\kappa_1(x_{i+1},y_j)\right]f^0_x(x_{i+1},y_j), \leqno{ \ \ \ \bullet}$$ 
$$\parallel \textstyle{\frac{\partial}{\partial x}}[\bs{u}^{\delta_n}_i-(f^*)^0](x,y_j)\parallel \leqno{ \ \ \ \bullet}$$
$$=\left|(\bar\sigma_3(x,y_j)-\bar\sigma_3(x_i,y_j))-(\bar\kappa_3(x,y_j)-\bar\kappa_3(x_i,y_j)\kappa_1(x,y_j)\right|~\parallel f^0_x(x_i,y_j)\parallel$$
$$\leq 2C_1^2C_2(x-x_0)\leq 2C_1^2C_2(a/n) \ \ \ \ for \ x\in[x_i,x_{i+1}].$$
Next, we define the curve $\bs{u}^{\delta_n}(x,y_j)(:=\bs{u}^{\delta_n}(x,y_j,c_0))$ on $[x_0,x_n]\times\{y_j\}$ from $\bs{u}^{\delta_n}_i(x,y_j)$  by 
$$\bs{u}^{\delta_n}(x,y_j)=\Sigma_{l=0}^{i-1}\bs{u}^{\delta_n}_{l}(x_{l+1},y_j)+\bs{u}^{\delta_n}_{i}(x,y_j) \ \ \ for \  x\in [x_{i},x_{i+1}], \ \ i\in\{0,1,\cdots, n-1\}, \eqno{(5.5)}$$
where $\Sigma_{l=0}^{i-1}\bs{u}^{\delta_n}_{l}(x_{l+1},y_j)=:{\bf 0}$ for $i=0$. Then, $\bs{u}^{\delta_n}(x_0,y_j)={\bf 0}$ holds and $\bs{u}^{\delta_n}(x,y_j)$ is continuous on $[x_0,x_e]$. 
Furthermore, we have 
$$\parallel \textstyle\frac{\partial}{\partial x}\left[\bs{u}^{\delta_n}-(f^*)^0\right] (x,y_j)\parallel\leq 2C_1^2C_2(a/n) \ \ \ for \ x\in [x_0,x_e] \eqno{(5.6)}$$
by the above equation for $\parallel ({\partial}/{\partial x})[\bs{u}^{\delta_n}_i-(f^*)^0](x,y_j)\parallel$.

Similarly, for each $x=x_i$ ($0\leq i\leq n$), we define a curve $\bs{v}^{\delta_n}_j(x_i,y)(:=\bs{v}^{\delta_n}_j(x_i,y,c_0))$ on $\{x_i\}\times[y_j,y_{j+1}]$ ($0\leq j\leq n-1$) by 
$$\bs{v}^{\delta_n}_j(x_i,y):=-\bar\sigma_3(x_{i},y_j)\left[f^0(x_i,y)\right]_{y=y_j}^{y}-\bar\kappa_3(x_i,y_j)\left[N^0(x_i,y)\right]_{y=y_j}^{y}.  \eqno{(5.7)}$$
Then, $\bs{v}^{\delta_n}_j(x_i,y_j)={\bf 0}$ and the tangent vectors of two curves $\bs{v}^{\delta_n}_j(x_i,y)$ and $(f^*)^0(x_i,y)$ are parallel on $[y_j,y_{j+1}]$. 
In particular, we have the following equations:
$$(\partial\bs{v}^{\delta_n}_j/\partial y)(x_i,y_j{\small +0})=-\left[\bar\sigma_3-\bar\kappa_3\bar\kappa_2\right](x_i,y_j)f^0_y(x_i,y_j)=\left(\bar\sigma_2f^0_y\right)(x_i,y_j)=(f^*)^0_y(x_i,y_j),   \leqno{ \ \bullet}$$ 
$$(\partial\bs{v}^{\delta_n}_{j}/\partial y)(x_i,y_{j+1}{\small -0})=-\left[\bar\sigma_3(x_i,y_{j})-\bar\kappa_3(x_i,y_{j})\bar\kappa_2(x_i,y_{j+1})\right]f^0_y(x_i,y_{j+1}),  \leqno{ \ \bullet}$$ 
$$\parallel \textstyle\frac{\partial}{\partial y}[\bs{v}^{\delta_n}_j)-(f^*)^0](x_i,y)\parallel\leq 2C_1^2C_2(a/n) \ \ \ \ for \ y\in[y_j,y_{j+1}].   \leqno{ \ \bullet}$$
Next, we define the curve $\bs{v}^{\delta_n}(x_i,y)(:=\bs{v}^{\delta_n}(x_i,y,c_0))$ on $\{x_i\}\times[y_0,y_e]$ by
$$\bs{v}^{\delta_n}(x_i,y):=\Sigma_{l=0}^{j-1}\bs{v}^{\delta_n}_{l}(x_i,y_{l+1})+\bs{v}^{\delta_n}_j(x_i,y) \ \ \ \ {\it for} \  y\in [y_j,y_{j+1}], \ \ j\in\{0,1,\cdots, n-1\},  \eqno{(5.8)}$$
where $\Sigma_{l=0}^{j-1}\bs{v}^{\delta_n}_{l}(x_i,y_{l+1})=:{\bf 0}$ for $j=0$. Then, $\bs{v}^{\delta_n}(x_i,y_0)={\bf 0}$ and $\bs{v}^{\delta_n}(x_i,y)$ is continuous on $[y_0,y_e]$. Furthermore, we have
$$\parallel \textstyle\frac{\partial}{\partial y}\left[\bs{v}^{\delta_n}-(f^*)^0\right](x_i,y)\parallel\leq 2C_1^2C_2(a/n) \ \ \ for \ y\in[y_0,y_e]. \eqno{(5.9)}$$

{For Step 2}: In this step, we assume $(f^*)^0(x_0,y_0)=f^*(x_0,y_0,c_0)={\bf 0}$ for the sake of simplicity. Under the condition, we define two surfaces ${S}^{\overline{\delta}_n}(x,y)(:=S^{\overline{\delta}_n}(x,y,c_0))$ and ${S}^{\underline{\delta}_n}(x,y)(:=S^{\underline{\delta}_n}(x,y,c_0))$, which are approximations of the $(x,y)$-surface $(f^*)^0(x,y)=f^*(x,y,c_0)$. Then, the condition $(f^*)^0(x_0,y_0)={\bf 0}$ corresponds to ${S}^{\overline{\delta}_n}(x_0,y_0)={S}^{\underline{\delta}_n}(x_0,y_0)={\bf 0}$. 

The surface $S^{\overline{\delta}_n}(x,y)$ is defined on $(\{x_0\}\times[y_0,y_e])\cup_{j=0}^n([x_0,x_e]\times \{y_j\})$ by
$$S^{\bar{\delta}_n}(x,y):=
\left[
\begin{array}{ll}
\bs{v}^{\delta_n}(x_0,y) \ \ \ for \ x=x_0, \ \ y\in[y_0,y_e] \\[2mm]
\bs{v}^{\delta_n}(x_0,y_j)+\bs{u}^{\delta_n}(x,y_j) \ \ \ for \ x\in[x_0,x_e], \ \ y=y_j \ (j=0,1,\cdots,n).
\end{array} \right. 
\eqno{(5.10)}$$
The surface ${S}^{\overline{\delta}_n}(x,y)$ is continuous and the tangent vectors of two $x$-curves ${S}^{\overline{\delta}_n}(x,y_j)$ and $(f^*)^0(x,y_j)$ with fixed $y=y_j$ are parallel on the definition domain of ${S}^{\overline{\delta}_n}(x,y)$.

Similarly, the surface ${ S}^{\underline{\delta}_n}(x,y)$ is defined on $([x_0,x_e]\times\{y_0\})\cup_{i=0}^n(\{x_i\}\times[y_0,y_e])$ by  
$${ S}^{\underline{\delta}_n}(x,y):=\left[
\begin{array}{ll}
\bs{u}^{\delta_n}(x,y_0) \ \ \ for \ x\in[x_0,x_e], \ \ y=y_0 \\[2mm]
\bs{u}^{\delta_n}(x_i,y_0)+\bs{v}^{\delta_n}(x_i,y) \ \ \ for \ x=x_i \ (i=0,1,\cdots,n), \ \ y\in[y_0,y_e]
\end{array}\right. \eqno{(5.11)}$$
The surface ${ S}^{\underline{\delta}_n}(x,y)$ is continuous and the tangent vectors of two $y$-curves ${ S}^{\underline{\delta}_n}(x_i,y)$ and $(f^*)^0(x_i,y)$ with fixed $x=x_i$ are parallel on the definition domain of ${ S}^{\underline{\delta}_n}(x,y)$. Note that both surfaces ${ S}^{\overline{\delta}_n}(x,y)$ and ${ S}^{\underline{\delta}_n}(x,y)$ are constructed directly from $2$-dimensional subnets of $f^0(L_n\cap E)$. 

The following Theorem 5.2 gives the relation between the corresponding $(x,y)$-curvature surfaces of $f(x,y,z)$ and $f^*(x,y,z)$, via the approximations ${ S}^{\overline{\delta}_n}(x,y)$ and ${ S}^{\underline{\delta}_n}(x,y)$ of $(f^*)^0(x,y)=f^*(x,y,c_0)$. In the theorem, the lattice $L_n(\subset Q)$ for $n$ are defined as above. Let ${ S}^{\overline{\delta}_n}(x,y)$ and ${ S}^{\underline{\delta}_n}(x,y)$ be the surfaces defined in (5.10) and (5.11), respectively.\\

{\thm 5.2} \ \ {\it Let $f(x,y,z)$ be a generic conformally flat hypersurface on a connected open set $U$ and $f^*(x,y,z)$ be the dual of $f(x,y,z)$. Let $n$ be a positive integer. For an arbitrarily fixed $c_0:=z_k \ (k=0,1,\cdots,n)$,
let $E:=[x_0,x_e]\times[y_0,y_e]\times\{c_0\}$ be a compact square in $U^*$ with $a:=x_e-x_0=y_e-y_0>0$.  
Let $f^0(x,y):=f(x,y,c_0)$ and $(f^*)^0(x,y):=f^*(x,y,c_0)$ be the curvature surfaces on $E$ of $f(x,y,z)$ and $f^*(x,y,z)$, respectively. 
Suppose $(f^*)^0(x_0,y_0)={\bf 0}$. 
Then, the surfaces ${ S}^{\overline{\delta}_n}(x,y)$ and ${ S}^{\underline{\delta}_n}(x,y)$ defined for each $n$ are approximations of $(f^*)^0(x,y)$ constructed by focusing on the principal curvature $x$- and $y$-lines of $(f^*)^0(x,y)$, respectively.
That is, we have the following inequalities:
$$\parallel \left[(f^*)^0-{ S}^{\overline{\delta}_n}\right](x,y)\parallel<4C_1^2C_2({a^2}/{n}) \ \ \ for \ (x,y)\in (\{x_0\}\times[y_0,y_e])\cup_{j=0}^n([x_0,x_e]\times\{y_j\}),$$
$$\parallel \left[(f^*)^0-{ S}^{\underline{\delta}_n}\right](x,y)\parallel<4C_1^2C_2({a^2}/{n}) \ \ \ for \ (x,y)\in([x_0,x_e]\times\{y_0\})\cup_{i=0}^n(\{x_i\}\times[y_0,y_e]).$$}

{\it{Proof}.} \ We only show the inequality for ${ S}^{\overline{\delta}_n}(x,y)$, because the inequality for ${ S}^{\underline{\delta}_n}(x,y)$ can be shown in the same way. We denote by $(x,y)$ the point $(x,y,c_0)\in E$. Now, on the definition domain of ${ S}^{\overline{\delta}_n}(x,y)$, the path connecting two points $(x_0,y_0)$ and $(x,y_j)$ ($x\in[x_0,x_e]$) is only one, which is given by $(x_0,y_0)\rightarrow (x_0,y_j)\rightarrow (x,y_j)$. Then, we have 
$$\parallel\left[\bs{v}^{\delta_n}-(f^*)^0\right](x_0,y_i)\parallel\leq \int_{y_0}^{y_j}\parallel \textstyle\frac{\partial}{\partial y}\left[\bs{v}^{\delta_n}-(f^*)^0\right](x_0,y)\parallel dy$$
$$\leq 2C_1^2C_2(a/n)(y_j-y_0)\leq 2C_1^2C_2(a^2/n)$$
by $[\bs{v}^{\delta_n}-(f^*)^0](x_0,y_0)={\bf 0}$ and (5.2). In the same way, we have 
$$\parallel\left[(\bs{u}^{\delta_n}-(f^*)^0)(x,y_j)\right]_{x_0}^{x}\parallel
\leq \int_{x_0}^{x}\parallel \textstyle\frac{\partial}{\partial x}[\bs{u}^{\delta_n}-(f^*)^0](x,y_j)\parallel dx\leq 2C_1^2C_2(a^2/n)$$
by (5.1). Hence, by the definition of $S^{\bar{\delta}_n}(x,y)$, we have 
$$\parallel [S^{\bar{\delta}_n}-(f^*)^0](x,y)\parallel\leq 4C_1^2C_2(a^2/n)$$ 
for all $(x,y)\in (\{x_0\}\times[y_0,y_e])\cup_{j=0}^n([x_0,x_e]\times\{y_j\})$.
\hspace{\fill}$\Box$\\

For Step 3: We construct an approximation $(f^*)^{\bar{\delta}_n}(x,y,z)$ of $f^*(x,y,z)$ on $Q=[x_0,x_e]\times[y_0,y_e]\times [z_0,z_e]$ for each $n$, by using (P1) and $f^*_z=-\sigma_2f_z-\kappa_2N_z$ of Lemma 5.1. In this Step, we assume $f^*(x_0,y_0,z_0)={\bf 0}$ for the sake of simplicity. Now, the interval $[z_0,z_e]$ has been divided into $n$ sub-intervals $z_0<z_1<\cdots <z_n=z_e$ of equal length $\delta_n=a/n$. We firstly define a curve $\bs{w}^{\delta_n}_k(x_0,y_0,z)$ for $z\in [z_k,z_{k+1}]$ by
$$\bs{w}^{\delta_n}_k(x_0,y_0,z):=-\sigma_2((x_0,y_0,z_k)\left[f(x_0,y_0,z)\right]_{z_k}^z-\kappa_2(x_0,y_0,z_k)\left[N(x_0,y_0,z)\right]_{z_k}^z. \eqno{(5.12)}$$
As same as in the cases $\bs{u}^{\delta_n}_i$ and $\bs{v}^{\delta_n}_j$,  the tangent vectors of two $z$-curves $\bs{w}^{\delta_n}_k(x_0,y_0,z)$ and $f^*(x_0,y_0,z)$ are parallel at any $z\in[z_k,z_{k+1}]$ and we have
$$\parallel \textstyle\frac{\partial}{\partial z}[\bs{w}^{\delta_n}_k-f^*](x_0,y_0,z)\parallel<2C_1^2C_2(a/n) \ \ \ for \ z\in[z_k,z_{k+1}].   \eqno{(5.13)}$$
Next, we define a curve $\bs{w}^{\delta_n}(x_0,y_0,z)$ for $z\in[z_0,z_e]$ by
$$\bs{w}^{\delta_n}(x_0,y_0,z):=\Sigma_{l=0}^{k-1}\bs{w}^{\delta_n}_{l}(x_0,y_0,z_{l+1})+\bs{w}^{\delta_n}_k(x_0,y_0,z) \ \ \  for \ z\in[z_k,z_{k+1}], \ k\in\{0,1,\cdots,n-1\},$$
where $\Sigma_{l=0}^{k-1}\bs{w}^{\delta_n}_{l}(x_0,y_0,z_{l+1})=:{\bf 0}$ for $k=0$.

Now, let $S^{\bar{\delta}_n}(x,y,z_k) \ (k=1,2,\cdots n)$ be the $(x,y)$-surfaces $S^{\bar{\delta}_n}(x,y)$ with $c_0=z_k$ in Theorem 5.2. Under the condition $\bs{w}^{\delta_n}(x_0,y_0,z_k)=S^{\bar{\delta}_n}(x_0,y_0,z_k)$, we attach to the curve $\bs{w}^{\delta_n}(x_0,y_0,z)$ the surfaces $\{S^{\bar{\delta}_n}(x,y,z_k)\}_k$ by parallel translation. Thus, we have obtained the hypersurface $(f^*)^{\bar{\delta}_n}(x,y,z)$ from the families $\{S^{\bar{\delta}_n}(x,y,z_k)\}$. 
Then, the definition domain $Q^{\bar{\delta}_n}$ of $(f^*)^{\bar{\delta}_n}(x,y,z)$ 
is given by
$$Q^{\bar{\delta}_n}:= \ \left(\{(x_0,y_0)\}\times[z_0,z_e]\right)\cup_k\left(\{x_0\}\times[y_0,y_e]\times\{z_k\}\right)\cup_{j,k}\left([x_0,x_e]\times\{(y_j,z_k)\}\right).$$

On the other hand, as the approximate $z$-curve in the definition of $(f^*)^{\underline{\delta}_n}(x,y,z)$, we adopt the following $z$-curve ${\bs z}^{\delta_n}(x_0,y_0,z)$ determined by the equation $f^*_z=-\sigma_1f_z-\kappa_1N_z$. That is, instead of the curve $\bs{w}^{\delta_n}_k(x_0,y_0,z)$, we firstly take a $z$-curve $\bs{z}^{\delta_n}_k(x_0,y_0,z)$ on $[z_k,z_{k+1}]$, where $\bs{z}^{\delta_n}_k(x_0,y_0,z)$ is defined by replacing $(\sigma_2,\kappa_2)$ with $(\sigma_1,\kappa_1)$ in (5.12). Then, note that the equation (5.13) with $\bs{w}^{\delta_n}_k(x_0,y_0,z)$ replaced by $\bs{z}^{\delta_n}_k(x_0,y_0,z)$ holds.  Next, just as $\bs{w}^{\delta_n}(x_0,y_0,z)$ is obtained from $\bs{w}^{\delta_n}_k(x_0,y_0,z)$, we obtain the curve $\bs{z}^{\delta_n}(x_0,y_0,z)$ for $z\in[z_0,z_e]$ from $\bs{z}^{\delta_n}_k(x_0,y_0,z)$. Then, from the $z$-curve $\bs{z}^{\delta_n}(x_0,y_0,z)$ and the surfaces $\{S^{\underline{\delta}_n}(x,y,z_k)\}_k$, the hypersurface $(f^*)^{\underline{\delta}_n}(x,y,z)$ is defined in the same way that $(f^*)^{\bar{\delta}_n}(x,y,z)$ is defined.   
Therefore, the definition domain $Q^{\underline{\delta}_n}$ of $(f^*)^{\underline{\delta}_n}(x,y,z)$ is given by 
$$Q^{\underline{\delta}_n}:= \ \left(\{(x_0,y_0)\}\times[z_0,z_e]\right)\cup_k\left([x_0,x_e]\times\{(y_0,z_k\}\right)\cup_{i,k}\left(\{x_i\}\times[y_0,y_e]\times\{z_k\}\right).$$

Now, from the method of constructing $(f^*)^{\bar{\delta}_n}(x,y,z)$ and $(f^*)^{\underline{\delta}_n}(x,y,z)$, we obtain the following corollary by (5.13) in the same way as the proof of Theorem 5.2. \\

{\cor 5.3}. \ \ {\it The hypersurfaces $(f^*)^{\bar{\delta}_n}(x,y,z)$ and $(f^*)^{\underline{\delta}_n}(x,y,z)$ approximate the dual $f^*(x,y,z)$ on the respective defined domains $Q^{\bar{\delta}_n}$ and $Q^{\underline{\delta}_n}$, that is, we have the following inequalities:
$$ \parallel [f^*-(f^*)^{\bar{\delta}_n}](x,y,z)\parallel\leq 6C_1^2C_2(a^3/n) \ \ \ for \ (x,y,z)\in Q^{\bar{\delta}_n},$$
$$ \parallel [f^*-(f^*)^{\underline{\delta}_n}](x,y,z)\parallel\leq 6C_1^2C_2(a^3/n) \ \ \ for \ (x,y,z)\in Q^{\underline{\delta}_n}.$$
} \\

{\bf 6. Hypersurfaces $(f^*)^{\bar\delta_n}$ and $(f^*)^{\underline{\delta}_n}$ are approximate discrete hypersurfaces of $f^*$} 

In this section, we use the same notations as in the previous section: $f(x,y,z)$ is a generic conformally flat hypersurface on $U\subset \R^3$ and $f^*(x,y,z)$ is the dual of $f(x,y,z)$ on $U^*$; the set $Q=[x_0,x_e]\times[y_0,y_e]\times[z_0,z_e]$ is a compact cube in $U^*$ and $a=x_e-x_0=y_e-y_0=z_e-z_0>0$; for each positive integer $n$, the divisions of three sides of $Q$ of equal length $\delta_n=a/n$ are given by $x_0<x_1<\cdots<x_n=x_e$,  $y_0<y_1<\cdots<y_n=y_e$ and $z_0<z_1<\cdots<z_n=z_e$, and the set $L_n$ is the lattice in $Q$ determined from the divisions. 

Now, we fix $c_0\in \{1,2,\cdots, n\}$ arbitrarily, and then the compact square $E=[x_0,x_e]\times[y_0,y_e]\times\{c_0\}$ is identified with $[x_0,x_e]\times[y_0,y_e]$. 
In Theorem 5.2, for each positive integer $n$, we have constructed the approximation $S^{\overline{\delta}_n}(x,y)$ (resp. $S^{\underline{\delta}_n}(x,y)$) of the surface $(f^*)^0(x,y):=f^*(x,y,c_0)$ on $E$ by regarding $(f^*)^0(x,y)$ as the one parameter family of $x$-curves (resp. of $y$-curves). In this section, for any $1\leq i\leq n$ and $0\leq j\leq n-1$, we firstly construct a $y$-curve $l^{\delta_n}_{i,j}(y)$ on $[y_j,y_{j+1}]$, which connects two points $S^{\overline{\delta}_n}(x_i,y_j)$ and $S^{\overline{\delta}_n}(x_i,y_{j+1})$ and approximates $(f^*)^0(x_i,y)$ on $[y_j,y_{j+1}]$. 
Here, the fact that the $y$-curve $l^{\delta_n}_{i,j}(y)$ on $[y_j,y_{j+1}]$ connects $S^{\overline{\delta}_n}(x_i,y_j)$ and $S^{\overline{\delta}_n}(x_i,y_{j+1})$ means that $l^{\delta_n}_{i,j}(y_j)=S^{\overline{\delta}_n}(x_i,y_j)$ and $l^{\bar\delta_n}_{i,j}(y_{j+1})=S^{\overline{\delta}_n}(x_i,y_{j+1})$, and the fact that the curve $l^{\delta_n}_{i,j}(y)$ approximates $(f^*)^0(x_i,y)$ implies that there is a constant $M$ depending only on $f$ and $Q$ such that 
$$\parallel  \textstyle\frac{\partial}{\partial y}\left[l^{\bar\delta_n}_{i,j}(y)-(f^*)^0(x_i,y)\right]\parallel<M/n \ \ \ \ for \ \ y\in [y_j,y_{j+1}] \eqno{(6.1)}$$ 
is satisfied. Because, by (6.1), we have the inequality $\parallel l^{\bar\delta_n}_{i,j}(y)-(f^*)^0(x_i,y)\parallel<\parallel S^{\bar\delta_n}_{i,j}-(f^*)^0\parallel(x_i,y_j)+aM/n^2$ for $y\in[y_j,y_{j+1}]$. For simplicity of description, we denote the equation (6.1) by $(l^{\bar\delta_n}_{i,j})_y(y)=(f^*)^0_y(x_i,y)+O(1/n)$.

In the same way, for the surface $S^{\underline{\delta}_n}(x,y)$, an approximate $x$-curve of $(f^*)^0(x,y_j)$ connecting two points $S^{\underline{\delta}_n}(x_i,y_j)$ and $S^{\underline{\delta}_n}(x_{i+1},y_j)$ can be constructed. Hence, when the $y$-curves $l^{\delta_n}_{i,j}(y)$ for $S^{\overline{\delta}_n}(x,y)$ are made, the surfaces $S^{\overline{\delta}_n}(x,y)$ and $S^{\underline{\delta}_n}(x,y)$ are approximate discrete surfaces of $(f^*)^0(x,y)$, defined on $L_n\cap(R^2_{(x,y)}\times\{c_0\})$.

Now, let us define $N^0(x,y):=N(x,y,c_0)$ for the unit normal vector field $N$ of $f(x,y,z)$ and $\bar k(x,y):=k(x,y,c_0)$ for a function $k(x,y,z)$ on $U$, as in the previous section. Let $\bs{v}^{\delta_n}_j(x_i,y)$ be the $y$-curve on $\{x_i\}\times[y_j,y_{j+1}]$ in (5.7). As shown in Step 1 of the previous section, the curve $S^{\overline{\delta}_n}(x_i,y_j)+\bs{v}^{\delta_n}_j(x_i,y)$ ($y_j\leq y\leq y_{j+1}$) is an approximation of the curve $(f^*)^0(x_i,y)$ ($y_j\leq y\leq y_{j+1}$), but the two points $S^{\overline{\delta}_n}(x_i,y_j)+\bs{v}^{\delta_n}_j(x_i,y_{j+1})$ and $S^{\overline{\delta}_n}(x_i,y_{j+1})$ are distinct from each other by the definition of $S^{\overline{\delta}_n}(x,y)$. Therefore, we must modify the curve $S^{\overline{\delta}_n}(x_i,y_j)+\bs{v}^{\delta_n}_j(x_i,y)$ on $[y_j,y_{j+1}]$ so that its value at $(x_i,y_{j+1})$ is $S^{\overline{\delta}_n}(x_i,y_{j+1})$.

We firstly define a constant vector $\bs{a}^{\delta_n}_{i,j}$ and a $y$-curve $\hat{\bs{v}}^{\delta_n}_j(x_i,y)$ on $[y_j,y_{j+1}]$ by 
$$\bs{a}^{\delta_n}_{i,j}:=S^{\overline{\delta}_n}(x_i,y_{j+1})-S^{\overline{\delta}_n}(x_i,y_j),$$
$$\hat{\bs{v}}^{\delta_n}_j(x_i,y):=\bar\sigma_3(x_i,y_j)\left[f^0(x_i,y)\right]^{y_{j+1}}_{y=y}+\bar\kappa_3(x_i,y_j)\left[N^0(x_i,y)\right]^{y_{j+1}}_{y=y}.$$Then, we have $\hat{\bs{v}}^{\delta_n}_j(x_i,y_j)=-\bs{v}^{\delta_n}_j(x_i,y_{j+1})$, \ $\hat{\bs{v}}^{\delta_n}_j(x_i,y_{j+1})={\bf 0}$ and 
$$(\partial\hat{\bs{v}}^{\delta_n}_j/\partial y)(x_i,y)
=\left(-\bar\sigma_3(x_i,y_j)+\bar\kappa_3(x_i,y_j)\bar\kappa_2(x_i,y)\right)f^0_y(x_i,y)=
(\partial\bs{v}^{\delta_n}_j/\partial y)(x_i,y).$$
Hence, the tangent vector of the curve $\hat{\bs{v}}^{\delta_n}_j(x_i,y)$ is parallel to that of $(f^*)^0(x_i,y)$ as well as $\bs{v}^{\delta_n}_j(x_i,y)$.
Next, using the equation $(n/a)[(y_{j+1}-y)+(y-y_j)]=1$, we define a $y$-curve $\tilde{\bs{v}}^{\delta_n}_j(x_i,y)$ ($y_j\leq y\leq y_{j+1}$) by
$$\tilde{\bs{v}}^{\delta_n}_j(x_i,y):=(n/a)\left[(y_{j+1}-y)\bs{v}^{\delta_n}_j(x_i,y)
+(y-y_j)(\bs{a}^{\delta_n}_{i,j}+\hat{\bs{v}}^{\delta_n}_j(x_i,y))\right].$$
Then, $\tilde{\bs{v}}^{\delta_n}_j(x_i,y_j)={\bf 0}$ and $\tilde{\bs{v}}^{\delta_n}_j(x_i,y_{j+1})=\bs{a}^{\delta_n}_{i,j}$ hold. Hence, $S^{\overline{\delta}_n}(x_i,y_j)+\tilde{\bs{v}}^{\delta_n}_j(x_i,y)(=l^{\bar\delta_n}_{i,j}(y))$ is a curve connecting $S^{\overline{\delta}_n}(x_i,y_j)$ and $S^{\overline{\delta}_n}(x_i,y_{j+1})$.

Now, we shall verify the fact that the curve $S^{\overline{\delta}_n}(x_i,y_j)+\tilde{\bs{v}}^{\delta_n}_j(x_i,y)$ is an approximation of the $y$-curve $(f^*)^0(x_i,y)$ on $[y_j,y_{j+1}]$. 
Then, since we have 
$$(\partial\tilde{\bs{v}}^{\delta_n}_j/\partial y)(x_i,y)=(n/a)\left[
\bs{a}^{\delta_n}_{i,j}-\bs{v}^{\delta_n}_j(x_i,y)+\hat{\bs{v}}^{\delta_n}_j(x_i,y)\right]+(\partial\bs{v}^{\delta_n}_j/\partial y)(x_i,y)$$
and $(\partial\bs{v}^{\delta_n}_j/\partial y)(x_i,y)=(f^*)^0_y(x_i,y)+O(a/n)$ by Step 1 of the previous section, it is sufficient to verify the following equation: with $1\leq i\leq n$ and $0\leq j\leq n-1$, 
$$\parallel \bs{a}^{\delta_n}_{i,j}-\bs{v}^{\delta_n}_j(x_i,y)+\hat{\bs{v}}^{\delta_n}_j(x_i,y)\parallel=O(a^3/n^2) \ \ \ for \ \ y\in[y_j,y_{j+1}].\eqno{(6.2)}$$

The equation (6.2) will be verified in Theorem 6.4, after three Lemmata 6.1-6.3 below. 
\\

{\lem 6.1}. \ \ {\it With $1\leq i\leq n$ and $0\leq j\leq n-1$, the following equation holds:
$$-\bs{v}^{\delta_n}_j(x_i,y)+\hat{\bs{v}}^{\delta_n}_j(x_i,y)=-\bs{v}^{\delta_n}_j(x_i,y_{j+1}) \ \ \ for \ y  \in[y_j,y_{j+1}].$$ 
In particular, the vector $\bs{a}^{\delta_n}_{i,j}-\bs{v}^{\delta_n}_j(x_i,y)+\hat{\bs{v}}^{\delta_n}_j(x_i,y)$ of $y$ in (6.2) is the constant vector $\bs{a}^{\delta_n}_{i,j}-\bs{v}^{\delta_n}_j(x_i,y_{j+1})$ on $[y_j,y_{j+1}]$.
}\\

{\it{Proof}.} \ We have
$$-\bs{v}^{\delta_n}_j(x_i,y)+\hat{\bs{v}}^{\delta_n}_j(x_i,y)$$
$$=\bar\sigma_3(x_i,y_j)\left(f^0(x_i,y)-f^0(x_i,y_j)\right)+\bar\kappa_3(x_i,y_j)\left(N^0(x_i,y)-N^0(x_i,y_j)\right)$$
$$+\bar\sigma_3(x_i,y_j)\left(f^0(x_i,y_{j+1})-f^0(x_i,y)\right)+\bar\kappa_3(x_i,y_j)\left(N^0(x_i,y_{j+1})-N^0(x_i,y)\right)$$
$$=\bar\sigma_3(x_i,y_j)\left[f^0(x_i,y)\right]^{y_{j+1}}_{y=y_j}+\bar\kappa_3(x_i,y_j)\left[N^0(x_i,y)\right]^{y_{j+1}}_{y=y_j}=-\bs{v}^{\delta_n}_j(x_i,y_{j+1}).$$
\hspace{\fill}$\Box$\\

In the following lemma, $\bs{v}^{\delta_n}_j(x_k,y_{j+1})+\bs{u}^{\delta_n}_k(x_{k+1},y_{j+1})$ and $\bs{u}^{\delta_n}_k(x_{k+1},y_j)+\bs{v}^{\delta_n}_j(x_{k+1},y_{j+1})$ are points in $\bs R^4$ obtained by the mapping along paths $(x_k,y_j)\rightarrow(x_k,y_{j+1})\rightarrow(x_{k+1},y_{j+1})$ and $(x_k,y_j)\rightarrow(x_{k+1},y_j)\rightarrow(x_{k+1},y_{j+1})$, respectively, under the condition that $(x_k,y_j)$ maps to the origin. \\

{\lem 6.2}. \ \ {\it Let $j\in\{0,1,\cdots, n-1\}$. For the vector in Lemma 6.1, we have the following equations:
$$\bs{a}^{\delta_n}_{1,j}-\bs{v}^{\delta_n}_j(x_1,y_{j+1})=
\left(\bs{v}^{\delta_n}_j(x_0,y_{j+1})+\bs{u}^{\delta_n}_0(x_1,y_{j+1})\right)-\left(\bs{u}^{\delta_n}_0(x_1,y_j)+\bs{v}^{\delta_n}_j(x_1,y_{j+1})\right). \leqno{ \ \ \ \bullet}$$
$$\bs{a}^{\delta_n}_{k+1,j}-\bs{v}^{\delta_n}_j(x_{k+1},y_{j+1})=
\bs{a}^{\delta_n}_{k,j}-\bs{v}^{\delta_n}_j(x_{k},y_{j+1}) \leqno{ \ \ \ \bullet}$$
$$+\left(\bs{v}^{\delta_n}_j(x_k,y_{j+1})+\bs{u}^{\delta_n}_k(x_{k+1},y_{j+1})\right)-\left(\bs{u}^{\delta_n}_k(x_{k+1},y_j)+\bs{v}^{\delta_n}_j(x_{k+1},y_{j+1})\right) \ \ for \ \ 1\leq k\leq n-1.$$}

{\it{Proof}.} \ For the first equation, we have the following:
$$\bs{a}^{\delta_n}_{1,j}-\bs{v}^{\delta_n}_j(x_1,y_{j+1})$$
$$=S^{\bar\delta_n}(x_1,y_{j+1})-S^{\bar\delta_n}(x_1,y_{j})-\bs{v}^{\delta_n}_j(x_1,y_{j+1})$$
$$=\left(S^{\bar\delta_n}(x_1,y_{j+1})-S^{\bar\delta_n}(x_0,y_{j})\right)-\left(S^{\bar\delta_n}(x_1,y_{j})-S^{\bar\delta_n}(x_0,y_{j})+\bs{v}^{\delta_n}_j(x_1,y_{j+1})\right)$$
$$=\left(\bs{v}^{\delta_n}_j(x_0,y_{j+1})+\bs{u}^{\delta_n}_0(x_1,y_{j+1})\right)-\left(\bs{u}^{\delta_n}_0(x_1,y_j)+\bs{v}^{\delta_n}_j(x_1,y_{j+1})\right),$$
which shows the first equation. Next, we verify the second equation. For $1\leq k\leq n-1$, we have
$$\bs{a}^{\delta_n}_{k+1,j}-\bs{v}^{\delta_n}_j(x_{k+1},y_{j+1})$$
$$=S^{\bar\delta_n}(x_{k+1},y_{j+1})-S^{\bar\delta_n}(x_{k+1},y_{j})-\bs{v}^{\delta_n}_j(x_{k+1},y_{j+1})$$
$$=\left(S^{\bar\delta_n}(x_{k},y_{j+1})+\bs{u}^{\delta_n}_k(x_{k+1},y_{j+1})\right)-\left(S^{\bar\delta_n}(x_{k},y_{j})+\bs{u}^{\delta_n}_k(x_{k+1},y_{j})\right)-\bs{v}^{\delta_n}_j(x_{k+1},y_{j+1})$$
$$=\left(S^{\bar\delta_n}(x_{k},y_{j+1})-S^{\bar\delta_n}(x_{k},y_{j})-\bs{v}^{\delta_n}_j(x_{k},y_{j+1})\right)+\left(\bs{v}^{\delta_n}_j(x_{k},y_{j+1})+\bs{u}^{\delta_n}_k(x_{k+1},y_{j+1})\right)$$
$$-\left(\bs{u}^{\delta_n}_k(x_{k+1},y_{j})+\bs{v}^{\delta_n}_j(x_{k+1},y_{j+1})\right),$$
which shows the second equation.
\hspace{\fill}$\Box$\\

{\lem 6.3}. \ \ {\it For $0\leq i,j\leq n-1$, the following equation holds:
$$\parallel\left(\bs{v}^{\delta_n}_j(x_k,y_{j+1})+\bs{u}^{\delta_n}_k(x_{k+1},y_{j+1})\right)-\left(\bs{u}^{\delta_n}_k(x_{k+1},y_j)+\bs{v}^{\delta_n}_j(x_{k+1},y_{j+1})\right)\parallel=O((a/n)^3).$$}\\

{\it{Proof}.} \ The equation is verified by the Taylor expansion for functions and the compactness of $Q$. Firstly, note that we have $x_{i+1}-x_i=y_{j+1}-y_j=a/n$, $(\bar\sigma_3)_x=\bar\kappa_2(\bar\kappa_3)_x$, $(\bar\sigma_3)_y=\bar\kappa_1(\bar\kappa_3)_y$, $N^0_x=-\bar\kappa_1f^0_x$ and $N^0_y=-\bar\kappa_2f^0_y$. 
Now, we have
$$\bs{v}^{\delta_n}_j(x_k,y_{j+1})-\bs{v}^{\delta_n}_j(x_{k+1},y_{j+1})$$
$$=-\bar\sigma_3(x_k,y_j)(f^0(x_k,y_{j+1})-f^0(x_k,y_j))-\bar\kappa_3(x_k,y_j)(N^0(x_k,y_{j+1})-N^0(x_k,y_j))$$
$$+\bar\sigma_3(x_{k+1},y_j)(f^0(x_{k+1},y_{j+1})-f^0(x_{k+1},y_j))+\bar\kappa_3(x_{k+1},y_j)(N^0(x_{k+1},y_{j+1})-N^0(x_{k+1},y_j))$$
$$=(\bar\sigma_3(x_{k+1},y_j)-\bar\sigma_3(x_{k},y_j))(f^0(x_{k+1},y_{j+1})-f^0(x_{k+1},y_j))$$
$$+\bar\sigma_3(x_{k},y_j)\left[(f^0(x_{k+1},y_{j+1})-f^0(x_{k+1},y_j))-(f^0(x_k,y_{j+1})-f^0(x_k,y_j))\right]$$
$$+(\bar\kappa_3(x_{k+1},y_j)-\bar\kappa_3(x_{k},y_j))(N^0(x_{k+1},y_{j+1})-N^0(x_{k+1},y_j))$$
$$+\bar\kappa_3(x_{k},y_j)\left[(N^0(x_{k+1},y_{j+1})-N^0(x_{k+1},y_j))-(N^0(x_k,y_{j+1})-N^0(x_k,y_j))\right]$$
$$=(a/n)^2~[(\bar\sigma_3)_x(x_k,y_j)f^0_y(x_{k+1},y_j)+(\bar\kappa_3)_x(x_k,y_j)N^0_y(x_{k+1},y_j)]$$
$$+(a/n)^2~\left[\bar\sigma_3f^0_{yx}+\bar\kappa_3(N^0)_{yx}\right](x_k,y_j)+O((a/n)^3)$$
$$=(a/n)^2~[(\bar\sigma_3)_xf^0_y+(\bar\kappa_3)_xN^0_y](x_k,y_j)
+(a/n)^2~\left[\bar\sigma_3f^0_{yx}+\bar\kappa_3(N^0)_{yx}\right](x_k,y_j)+O((a/n)^3)$$
$$=(a/n)^2~\left[\bar\sigma_3f^0_{yx}+\bar\kappa_3(N^0)_{yx}\right](x_k,y_j)+O((a/n)^3).$$

In the same way, we have 
$$\bs{u}^{\delta_n}_k(x_{k+1},y_{j+1})-\bs{u}^{\delta_n}_k(x_{k+1},y_j)
=-(a/n)^2~\left[\bar\sigma_3f^0_{xy}+\bar\kappa_3N^0_{xy}\right](x_k,y_j)+O((a/n)^3).$$
In consequence, we have completed the proof. 
\hspace{\fill}$\Box$\\

{\thm 6.4}. \ {\it We fix $c_0\in\{z_0,z_1,\cdots,z_n\}$ arbitrarily. Let ${S}^{\bar \delta_n}(x,y)$ be the approximation of the surface $(f^*)^0(x,y):=f^*(x,y,c_0)$ on $E:=[x_0,x_e]\times[y_0,y_e]\times\{c_0\}$, in Theorem 5.2. Then, with $1\leq i\leq n$ and $0\leq j\leq n-1$, the curve ${S}^{\bar\delta_n}(x_i,y_j)+\tilde{\bs{v}}^{\delta_n}_j(x_i,y)$ for $y_j\leq y\leq y_{j+1}$ connects two points ${S}^{\bar \delta_n}(x_i,y_j)$ and ${S}^{\bar \delta_n}(x_i,y_{j+1})$, and it approximates the curve $(f^*)^0(x_i,y)$ for $y_j\leq y\leq y_{j+1}$. That is, the following three equations are satisfied for all $(i,j)$: 
$$\tilde{\bs{v}}^{\delta_n}_j(x_i,y_j)={\bf 0}, \ \ \ \tilde{\bs{v}}^{\delta_n}_j(x_i,y_{j+1})={S}^{\bar \delta_n}(x_i,y_{j+1})-{S}^{\bar \delta_n}(x_i,y_j), $$  $$(\partial\tilde{\bs{v}}^{\delta_n}_j/\partial y)(x_i,y)=(f^*)^0_y(x_i,y)+O(1/n) \ \ \ for \ y\in[y_j,y_{j+1}]. $$ 
These facts implies that ${S}^{\bar \delta_n}(x,y)$ is an approximate discrete surface of $(f^*)^0(x,y)$, defined on $E\cap L_n$.   
}\\

{\it{Proof}.} \ Only the proof of the equation (6.2) is reminded. 
Let us fix $j \ (0\leq j\leq n-1)$ arbitrarily. 
By Lemma 6.1, we need to verify the equations $\bs{a}_{i,j}^{\delta_n}-\bs{v}^{\delta_n}_j(x_i,y_{j+1})=O(a^3/n^2)$ for all $i\in\{1,2,\cdots,n\}$ to get (6.2). However, as the following proof process shows, it is sufficient to verify the equation $\bs{a}_{n,j}^{\delta_n}-\bs{v}^{\delta_n}_j(x_n,y_{j+1})=O(a^3/n^2)$ only for $i=n$. 

Now, by Lemmata 6.2 and 6.3, we obtain 
$$\bs{a}_{n,j}^{\delta_n}-\bs{v}^{\delta_n}_j(x_n,y_{j+1})=\bs{a}_{n-1,j}^{\delta_n}-\bs{v}^{\delta_n}_j(x_{n-1},y_{j+1})+O((a/n)^3)$$
$$=\bs{a}_{n-2,j}^{\delta_n}-\bs{v}^{\delta_n}_j(x_{n-2},y_{j+1})+2\times O((a/n)^3)=\cdots$$
$$=\bs{a}_{1,j}^{\delta_n}-\bs{v}^{\delta_n}_j(x_{1},y_{j+1})+(n-1)\times O((a/n)^3)=O(a^3/n^2).$$
In consequence, the equation (6.2) has been verified for $1\leq i\leq n$, $0\leq j\leq n-1$.
\hspace{\fill}$\Box$\\

Now, by Theorem 6.4, two hypersurfaces $(f^*)^{\bar\delta_n}(x,y,z)$ and $(f^*)^{\underline{\delta}_n}(x,y,z)$ in Step 3 of the previous section have been extended to the approximations of the dual $f^*(x,y,z)$, defined on 
$$\left[L_n\cap(\cup_{k=0}^n (\R^2_{(x,y)}\times\{z_k\})\right]\cup (\{(x_0,y_0)\}\times[z_0,z_e]).$$
 
Next, we verify that $(f^*)^{\bar\delta_n}(x,y,z)$ is an approximate discrete hypersurface of $f^*(x,y,z)$, defined on $L_n$. That is,  
for every $(i,j,k)$ we define an approximate $z$-curve connecting two points $(f^*)^{\bar\delta_n}(x_i,y_j,z_k)$ and $(f^*)^{\bar\delta_n}(x_i,y_j,z_{k+1})$, by using $f^*_z=-\sigma_2f_z-\kappa_2N_z$ in Lemma 5.1. 
%Then, the $z$-curve $f^*(x_i,y_j,z)$ on $[z_k,z_{k+1}]$, which is verified by using the other equation $f^*_x=-\sigma_2f_x-\kappa_2N_x$.    

Let us define the vectors ${\bs b}^{\delta_n}_{i,j,k}$ and the $z$-curves $\bs{z}_{i,j,k}(z)$ and $\hat{\bs z}_{i,j,k}(z)$ on $[z_k,z_{k+1}]$ by 
$${\bs b}^{\delta_n}_{i,j,k}:=(f^*)^{\bar{\delta}_n}(x_i,y_j,z_{k+1})-(f^*)^{\bar{\delta}_n}(x_i,y_j,z_k),$$
$$\bs{z}_{i,j,k}(z):=-\sigma_2(x_i,y_j,z_k)[f(x_i,y_j,z)]_{z=z_k}^z-
\kappa_2(x_i,y_j,z_k)[N(x_i,y_j,z)]_{z=z_k}^z,$$
$$\hat{\bs{z}}_{i,j,k}(z):=\sigma_2(x_i,y_j,z_k)[f(x_i,y_j,z)]_{z=z}^{z_{k+1}}+
\kappa_2(x_i,y_j,z_k)[N(x_i,y_j,z)]_{z=z}^{z_{k+1}}.$$
Then, we have $\bs{z}_{i,j,k}(z_k)=\hat{\bs{z}}_{i,j,k}(z_{k+1})=\bs{0}$ and 
$$d\bs{z}_{i,j,k}/dz=f^*_z(x_i,y_j,z)+O(1/n), \ \ \ d\hat{\bs{z}}_{i,j,k}/dz=f^*_z(x_i,y_j,z)+O(1/n).$$
Hence, for the following $z$-curve $\tilde{\bs{z}}_{i,j,k}(z)$ on $[z_k,z_{k+1}]$,
$$\tilde{\bs{z}}_{i,j,k}(z):=(n/a)[(z_{k+1}-z)\bs{z}_{i,j,k}(z)+(z-z_k)[{\bs b}^{\delta_n}_{i,j,k}+\hat{\bs{z}}_{i,j,k}(z)],$$
the $z$-curve $(f^*)^{\bar\delta_n}(x_i,y_j,z_k)+\tilde{\bs{z}}_{i,j,k}(z)$ connects two points $(f^*)^{\bar\delta_n}(x_i,y_j,z_k)$ and $(f^*)^{\bar\delta_n}(x_i,y_j,z_{k+1})$, and we have 
$$d\tilde{\bs z}_{i,j,k}/dz=(n/a)[{\bs b}^{\delta_n}_{i,j,k}-{\bs z}_{i,j,k}(z)+\hat{\bs z}_{i,j,k}(z)]+f^*_z(x_i,y_j,z)+O(1/n).$$

Now, we turn our attention to an $(x,z)$-surface of $f^*(x,y_j,z)$ with fixed $y=y_j$. If the $x$-curves $\bs{u}(x,y_j,z_k)$ of $(f^*)^{\bar\delta_n}(x_i,y_j,z)$, which are $\bs{u}(x,y_j)$ with $c_0=z_k$ in Step 1, would be constructed from the equation $f^*_x=-\sigma_2f_x-\kappa_2N_x$, then each $z$-curve $(f^*)^{\bar\delta_n}(x_i,y_j,z_k)+\tilde{\bs{z}}_{i,j,k}(z)$ approximates $(f^*)^{\bar\delta_n}(x_i,y_j,z)$ on $[z_k,z_{k+1}]$, because $\left[(d/dz)(\tilde{\bs{z}}_{i,j,k}-f^*)\right](x_i,y_j,z)=O(1/n)$ for $z\in[z_k,z_{k+1}]$ holds as in Theorem 6.4. Even in our case, each $z$-curve $(f^*)^{\bar\delta_n}(x_i,y_j,z_k)+\tilde{\bs{z}}_{i,j,k}(z)$ also approximates $(f^*)^{\bar\delta_n}(x_i,y_j,z)$ on $[z_k,z_{k+1}]$. In fact, we define $x$-curves $\bar{\bs{u}}(x,y_j,z_k)$ from $f^*_x=-\sigma_2f_x-\kappa_2N_x$ in the same way that $\bs{u}(x,y_j,z_k)$ were defined from $f^*_x=-\sigma_3f_x-\kappa_3N_x$. Then, the equation (5.6) with $\bs{u}^{\delta_n}(x,y_j)$ replaced by $\bar{\bs{u}}^{\delta_n}(x,y_j,z_k)$ holds. Therefore, since $\parallel(\partial/\partial x)[\bs{u}^{\delta_n}-\bar{\bs{u}}^{\delta_n}](x,y_j,z_k)\parallel=O(1/n)$ and $\parallel(\bs{u}^{\delta_n}-\bar{\bs{u}}^{\delta_n})(x,y_j,z_k)\parallel=O(1/n)$ hold for $x\in[x_0,x_e]$, we have $\left[(d/dz)(\tilde{\bs{z}}_{i,j,k}-f^*)\right](x_i,y_j,z)=O(1/n)$ for $z\in[z_k,z_{k+1}]$.

Furthermore, the approximate $z$-curves connecting two points $(f^*)^{\underline{\delta}_n}(x_i,y_j,z_k)$ and $(f^*)^{\underline{\delta}_n}(x_i,y_j,z_{k+1})$ are constructed from $f^*_z=-\sigma_1f_z-\kappa_1N_z$ in the same way as above. In consequence, both hypersurfaces $(f^*)^{\bar\delta_n}(x,y,z)$ and $(f^*)^{\underline{\delta}_n}(x,y,z)$ are approximate discrete hypersurfaces of $f^*(x,y,z)$, defined on $L_n$.\\

{\thm 6.5}. \ {\it Both hypersurfaces $(f^*)^{\bar\delta_n}(x,y,z)$ and $(f^*)^{\underline{\delta}_n}(x,y,z)$ are approximate discrete hypersurfaces of $f^*(x,y,z)$, defined on $L_n$.}\\

\end{document}